\documentclass{amsart}
\usepackage{amsmath} 
\usepackage{amssymb}
\usepackage{mathrsfs}

\newtheorem{theorem}{Theorem}[section] 
\newtheorem{claim}[theorem]{Claim}
\newtheorem{conclusion}[theorem]{Conclusion}
\newtheorem{lemma}[theorem]{Lemma} 
 
\newtheorem{corollary}[theorem]{Corollary} 

\theoremstyle{definition}
\newtheorem{definition}[theorem]{Definition}

\newtheorem{example}[theorem]{Example}
\newtheorem{explanation}[theorem]{Explanation}
\newtheorem{fact}[theorem]{Fact}
\newtheorem{problem}[theorem]{Problem}

\newtheorem{observation}[theorem]{Observation} 
\newtheorem{conjecture}[theorem]{Conjecture}

\newtheorem{exercise}[theorem]{Exercise}
\newtheorem{discussion}[theorem]{Discussion}

\newtheorem{hypothesis}[theorem]{Hypothesis}

\theoremstyle{remark}
\newtheorem{remark}[theorem]{Remark}
\newtheorem{notation}[theorem]{Notation}
\newtheorem{question}[theorem]{Question}
\newtheorem{context}[theorem]{Context}

\newcommand{\rest}{{\restriction}}
 
\newcommand{\Dom}{{\rm Dom}}
\newcommand{\AP}{{\rm AP}}
\newcommand{\BP}{{\rm BP}}
\newcommand{\EC}{{\rm EC}}
\newcommand{\Av}{{\rm Av}}
\newcommand{\Th}{{\rm Th}}
\newcommand{\tp}{{\rm tp}}  
\newcommand{\aut}{{\rm aut}}  
\newcommand{\eq}{{\rm eq}}  
\newcommand{\id}{{\rm id}}  
\newcommand{\nsp}{{\rm nsp}}  
\newcommand{\cf}{{\rm cf}} 

\newcommand{\dcl}{{\rm dcl}}
\newcommand{\mxK}{{\rm mxK}}

\newcommand{\rang}{{\rm rang}}
\newcommand{\Seq}{{\rm Seq}}
\newcommand{\Min}{{\rm Min}}
\newcommand{\Rang}{{\rm Rang}}
\newcommand{\wilog}{{\rm without loss of generality}}

\newcommand{\then}{{\underline{then}}}
\newcommand{\when}{{\underline{when}}}
\newcommand{\Then}{{\underline{Then}}}

\newcommand{\If}{{\underline{if}}}
\newcommand{\Iff}{{\underline{iff}}}
\newcommand{\mn}{{\medskip\noindent}}
\newcommand{\sn}{{\smallskip\noindent}}

\newcommand{\cA}{{\mathscr A}}

\newcommand{\bd}{{\bf d}}
\newcommand{\bc}{{\bf c}}

\newcommand{\gC}{{\mathfrak C}}

\newcommand{\cD}{{\mathscr D}}

\newcommand{\cF}{{\mathscr F}}

\newcommand{\bbL}{{\mathbb L}}

\newcommand{\cP}{{\mathscr P}}

\newcommand{\bbQ}{{\mathbb Q}}

\newcommand{\cU}{{\mathscr U}}

\newcommand{\cY}{{\mathscr Y}}

\newcount\skewfactor
\def\mathunderaccent#1#2 {\let\theaccent#1\skewfactor#2
\mathpalette\putaccentunder}
\def\putaccentunder#1#2{\oalign{$#1#2$\crcr\hidewidth
\vbox to.2ex{\hbox{$#1\skew\skewfactor\theaccent{}$}\vss}\hidewidth}}

\newenvironment{PROOF}[2][\proofname.]
   {\begin{proof}[#1]\renewcommand{\qedsymbol}{\textsquare$_{\rm #2}$}}
   {\end{proof}}

\begin{document}

\title {Dependent theories and the generic pair conjecture}
\author {Saharon Shelah}
\address{Einstein Institute of Mathematics\\
Edmond J. Safra Campus, Givat Ram\\
The Hebrew University of Jerusalem\\
Jerusalem, 91904, Israel\\
 and \\
 Department of Mathematics\\
 Hill Center - Busch Campus \\ 
 Rutgers, The State University of New Jersey \\
 110 Frelinghuysen Road \\
 Piscataway, NJ 08854-8019 USA}
\email{shelah@math.huji.ac.il}
\urladdr{http://shelah.logic.at}
\thanks{The author would like to thank the Israel Science Foundation for
partial support of this research (Grant No. 242/03).
I would like to thank 
Alice Leonhardt for the beautiful typing. 
First version as F756 was 2006/Feb/14; and 900 - 2007/Jan/10. Publication 900}
 




\date{November 15, 2013}

\begin{abstract}
We try to understand complete types over a
somewhat saturated model of a complete first order theory which is
dependent (previously called NIP), 
by ``decomposition theorems for such types".  Our thesis
is that the picture of dependent theory is the combination of the one for
stable theories and the one for the theory of dense linear order or trees (and
first we should try to understand the quite saturated case).
As a measure of our progress, we give several
applications considering some test questions; in particular we try 
to prove the generic pair conjecture and do it for measurable cardinals.   
\end{abstract}

\maketitle
\numberwithin{equation}{section}
\setcounter{section}{-1}
\newpage

\centerline {Annotated Content}
\bigskip

\noindent
\S0 \quad Introduction, pg. \pageref{s:introduction}
\bigskip

\noindent
\S1 \quad Non-splitting construction, pg. \pageref{s:non}
\mn
\begin{enumerate}
\item[${{}}$]   [For $\kappa$-saturated $M$ and $N$ such that $M
\prec N$ we try to analyze $N$ over $M$ by finding $M_1,N_1$ such
that $M \prec M_1 \prec N_1,M \prec N \prec N_1$ 
and both $M_1/M,N_1/N$ are understood but in opposite ways.  
The first similar in some sense to the stable
situation, the second to the situation for order.]
\end{enumerate}
\bigskip

\noindent
\S2 \quad The type decomposition theorem, pg. \pageref{s:the}
\mn
\begin{enumerate}
\item[${{}}$]  [For $\kappa$-saturated $M \prec \gC$ and $\bar d \in {\frak
C}$ of length $< \theta^+$ we try to analyze the type tp$(\bar d,M)$
in two steps - pseudo stable and tree-like one. 
This is the main aim of the section (and a major part of the paper).
It is done by looking at $K_\ell$ and
mxK$^\ell_{\lambda,\kappa,\theta}$.  A consequence which fulfilled to some
extent the aim is the Type Decomposition Theorem (\ref{tp25.43}).  As
a second consequence we
give a characterization of ``$M$ is exactly $\kappa$-saturated, $\kappa >
\text{ cf}(\kappa) > |T|$", see \ref{tp16.14}.  
In fact, we deal a little with singular exact saturation per
se. ``Unfortunately" there are independent (complete first order theories) $T$ 
which has no model with singular exact saturation, see
\ref{tp16.15}.  But the
existence of an indiscernible set for dependent $T$ suffice
(see \ref{tp16.17} under instances of GCH) and has a neat characterization.  
Also, if $p$ is a complete 1-type over a model $M$ of $T$ which is quite 
saturated \then \, $p$ has a spectrum in a suitable sense, see \ref{tp35.46}.]
\end{enumerate}
\bigskip

\noindent
\S3 \quad Existence of strict decomposition, pg. \pageref{s:existence}
\mn
\begin{enumerate}
\item[${{}}$]  [E.g. here complete types over 
a saturated model $M$ of cardinality $\kappa$, a
measurable cardinal, is analyzed.  What we get is a better
decomposition theorem (the strict one).]
\end{enumerate}
\bigskip

\noindent
\S4 \quad Consequences of strict decomposition, pg. \pageref{s:consequences}
\mn
\begin{enumerate}
\item[${{}}$]   [We start by sufficient conditions for a sequence
being indiscernible.  For a measurable cardinal $\kappa$ ($> |T|$) we
confirm the structure half of the generic pair conjecture.
Toward this, if we have the consequences of \S3 we can analyze
generic pairs of models of $T$ in $\kappa$.  In a slightly simplified
formulation this means: if $2^\kappa = \kappa^+,(\kappa =
\kappa^{< \kappa} > |T|),M_\alpha$, a model of $T$ of cardinality
$\kappa$ for $\alpha < \kappa^+$ is
$\prec$-increasing continuous, $M = \cup\{M_\alpha:\alpha <
\kappa^+\}$ is $\kappa^+$-saturated, \then \, for a club $E$ of 
$\kappa^+$ for all $\alpha < \beta$ belonging to $\{\delta \in
E:\delta$ has cofinality $\kappa\}$ the pair
$(M_\beta,M_\alpha)$ has the same isomorphism type.  In fact for $\kappa_1
< \kappa_2$ we get $\Bbb L_{\infty,\kappa_1}(\tau_T)$-equivalence, so
we have a derived first order theory. 
For the proof we show that an increasing (short) sequence of so called
strict $(\kappa,\theta)$-decompositions has a limit.]
\end{enumerate}
\newpage

\section{Introduction} \label{s:introduction}

We first give a page of introductory remarks for non-logicians to clarify
notions and to motivate working on dependent theories.  The classical
center of model theory is investigating elementary classes, i.e. we
fix a vocabulary $\tau$ (i.e. a set of predicates and function
symbols), for a $\tau$-structure $M$ let Th$(M)$ be the set of first
order sentences which $M$ satisfies,  a complete first order theory
$T$ is Th$(M)$ for some $\tau$-model $M$.  We fix $T$ and $\tau =
\tau_T$ and investigate $\tau$-models of $T$, 
i.e. $\tau$-structures $M$ such that $T = \text{ Th}(M)$; about 
other contexts, see e.g. \cite{Sh:E53}.

Let $M,N$ denote such structures and they are called models (of $T$).  
Let $\bar a,\bar b,\bar c,\bar d$ denote sequences of elements of such models
and $\varphi(\bar x)$ or $\varphi(\bar x,\bar y)$ denote members of
$\bbL(\tau)$, i.e. the set of first order formulas in this vocabulary
but we allow $\bar x$ to be infinite though the formula is finite so
only finitely many $x_\ell,y_j$ are relevant.  

Let $M \models \varphi[\bar a]$ mean that the model $M$ satisfies the
formula $\varphi(\bar x)$ under the substitution $\bar x \mapsto \bar
a$ (so $\bar a,\bar x$ have the same length).

The right notion of sub-models is $\prec$, being elementary submodel
where $M \prec N$ iff $M \subseteq N$ and for every $\varphi(\bar x)
\in \bbL(\tau)$ and $\bar a \in {}^{\ell g(\bar x)}M$ we have $M
\models \varphi[\bar a]$ iff $N \models \varphi[\bar a]$.                

Recall that an ordinal is the isomorphism type of a well ordering
(which is a linear order for which every non-empty set has a first
member).  But we identify an ordinal with the 
set of smaller ordinals.  Also a cardinal is an
ordinal $\lambda$ with no smaller ordinal of the
same power.  Here saying ``$x$ is a cardinal" means ``$x$ is an 
infinite cardinal" if not said
otherwise.  Let $\aleph_\alpha$ be the $\alpha$-th infinite cardinal 
and the cardinality $|\cU|$ of a set $\cU$ 
is the minimal ordinal of the same power.

Let the successor $\lambda^+$ of a cardinal $\lambda$ be $\aleph_{\alpha +1}$
when $\lambda = \aleph_\alpha$.

We say $E$ is a closed subset of the limit cardinal $\gamma$ when $E
\subseteq \gamma$ and $\delta < \gamma \wedge \delta =
\sup(\delta \cap E) \Rightarrow \delta \in E$ and $E$ is called
unbounded when $(\forall \alpha < \gamma)(\exists \beta)(\alpha
\le \beta \in E)$, ``$E$ is a club of $\gamma$" is the shorthand for
``$E$ is a closed unbounded subset of $\gamma$".  

For an ordinal $\alpha$ let cf$(\alpha) = \text{ min}\{|C|:C$ an
unbounded subset of $\alpha\} = \text{min}\{\text{otp}(C):C$ a closed unbounded
subset of $\alpha\}$; we say $\alpha$ is regular if $\alpha = \text{
cf}(\alpha)$ is infinite (hence is a cardinal), now recall (see
e.g. \cite{J}) that if $\alpha$ is a limit ordinal (e.g. a cardinal) 
then cf$(\alpha)$ is regular, and every cardinal of the form
$\lambda^+$ is regular.  When $\cf(\delta) > \aleph_0$ we say ``$S
\subseteq \delta$ is stationary" when $S \cap E \ne \emptyset$ for
every club $E$ of $\delta$.

A central notion is type; for $A \subseteq M$ and $\bar a$ a
sequence from $M$ let tp$(\bar a,A,M)$ be the set $\{\varphi(\bar
x,\bar b):\bar\varphi(\bar x,\bar y) \in \bbL(\tau),\bar b$ a sequence
from $A$ and $M \models \varphi[\bar a,\bar b]\}$.  We may write $a$
instead of $\langle a \rangle$.

Let 

\begin{equation*}
\begin{array}{clcr}
\bold S^\alpha(A,M) = \{\text{tp}(\bar a,A,N):&\text{ for some }
N,\bar a \text{ we have} \\
  &M \prec N,\bar a \text{ a sequence of length } \alpha \text{ from } N\}
\end{array}
\end{equation*}

\[
\bold S^\alpha(M) = \bold S^\alpha(M,M).
\]
\mn
By this we can define another central notion.  $M$ is
$\kappa$-saturated \Iff \, for every $A \subseteq M,|A| < \kappa$ and
$p \in \bold S^1(A,M)$ some $a \in M$ realizes $p$ in $M$ which means
$p = \text{ tp}(a,A,M)$.  We say the model $M$ is saturated when it is
$\kappa$-saturated and of cardinality $\kappa$.  Let EC$_\lambda(T)$ be
the class of models of $T$ of cardinality $\lambda$.

It is classically known that for $\lambda \ge |T|$, 
(assuming $2^\lambda = \lambda^+$, mostly done here for transparency) 
there is a saturated member of
EC$_{\lambda^+}(T)$, it is unique up to isomorphism, and the union of
an $\prec$-increasing chain of saturated members of
EC$_{\lambda^+}(T)$ of length $\lambda^+$ is a saturated member of
EC$_{\lambda^+}(T)$.
On the background so far, see e.g. Chang-Keisler \cite{CK73}.

\centerline {$* \qquad * \qquad *$}

A major theme of the author's work is trying to find natural dividing
lines (i.e. properties) in the family of first order complete $T$, a
criterion for natural is having both 
``inside definition" by formulas and ``outside definition"
by properties of the class of its models.  That is, such a property is
interesting as a dividing line when we have consequences for
those with the property and for those without it; see e.g. 
\cite[\S(1A)]{Sh:E53}.

A major such dividing line is ``$T$ is stable" which holds iff
$(*)^1_T$ iff $(*)^2_T$ where
\mn
\begin{enumerate}
\item[$(*)^1_T$]  for some $\varphi(\bar x,\bar y) \in \bbL(\tau_T)$,
model $M$ of $T$ and $\bar a_n \in {}^{(\ell g(\bar x))}M,\bar b_n \in
{}^{(\ell g(\bar y))}M$ for $n < \omega$ we have $n < m
\Leftrightarrow M \models \varphi[\bar a_n,\bar b_m]$
\sn 
\item[$(*)^2_T$]  for every $\lambda \ge |T|$ and limit ordinal $\delta
\le \lambda$ of cofinality $> |T|$, the union of any
$\prec$-increasing chain of length $\delta$
of saturated models of $T$ of cardinality $\lambda$ is saturated.
\end{enumerate}
\mn
Another major dividing line is ``$T$ is superstable" which holds iff
\mn
\begin{enumerate}
\item[$(*)^3_T$]  like $(*)^2_T$ allowing any limit ordinal $\delta$.
\end{enumerate}
\mn
On this and the relevant history, see e.g. \cite{Sh:c}.

The property we deal with here is ``$T$ is dependent", also called ``$T$ is
NIP", where its negation,
``$T$ is independent" or ``$T$ has the independence property" means
\mn
\begin{enumerate}
\item[$(*)^4_T$]  there are $\varphi(\bar x,\bar y) \in \bbL(\tau_T)$, a
model $M$ of $T$ and $\bar a_u \in {}^{(\ell g(\bar x))}M,\bar b_n \in
{}^{(\ell g(\bar y))}M$ for $u \subseteq \omega,n < \omega$ such that
$n \in u \Leftrightarrow M \models \varphi[\bar a_u,\bar b_n]$.
\end{enumerate}
\mn
What is the motivation to investigate this dividing line?  First, it
has a nice, simple definition, parallel to the one for stable
theories.  Second, it is a much wider class than that of the stable
theories; also, extremely important for many, whereas infinite fields with 
stable first order complete
theory are hard to come by (algebraically closed and separably closed
are the only known ones), 
there are many important fields with dependent first order complete
theory (the $p$-adics and many of the power series fields).  Third, there are
some results on it indicating it is not unreasonable to hope there is
a rich theory on it to be discovered.

On background on dependent theories, see \cite{Sh:715}, \cite{Sh:783}.
\bigskip

\noindent
\centerline {$* \qquad * \qquad *$}
\bigskip

Let $T$ be a fixed first order complete theory.  For transparency,
till \ref{0n.1}, we
assume G.C.H., i.e. $2^\kappa = \kappa^+$ for every infinite cardinal $\kappa$
and consider only $\lambda$ regular $> |T|$.  Let
$\bar M = \langle M_\alpha:\alpha < \lambda^+\rangle$ be an $\prec$-increasing
continuous sequence of models of $T$ of cardinality
$\lambda$ with $M$ being saturated where $M := \cup
\{M_\alpha:\alpha < \lambda^+\}$.  Now $M$ is unique
(up to isomorphism, for 
each $\lambda$) and though for a given $M,\bar M$ is not 
unique, for any two such
sequences $\bar M',\bar M''$ there is a closed unbounded subset 
$E$ of $\lambda^+$ and
isomorphism $f$ from
$M' = \cup\{M'_\alpha:\alpha < \lambda^+\}$ onto $M'' =
\cup\{M''_\alpha:\alpha < \lambda\}$ such that $f$ maps $M'_\delta$
onto $M''_\delta$ for every $\delta \in E$.

So it is natural to ask ($\lambda > |T|$ regular and $E$ varies on closed
unbounded subsets of $\lambda^+$)
\mn
\begin{enumerate}
\item[$\odot_1$]   what\footnote{We can present the problem
differently, about the existence of (variations of)
$(\lambda,\kappa)$-limit models (so $2^\lambda = \lambda^+$ is no
longer necessary, by forcing this is equivalent).  Also, instead of the
function $\bold n$ getting the value $\lambda^+$ we can consider
saying for some club no two relevant cases are isomorphic.  This 
does not make a real
difference but we find the present choice has more transparent presentation.}
is $\bold n_\lambda(T) := \text{ Min}_E|\{M_\delta/\cong:\delta \in
E\}|$? where $M_\delta/\cong$ is the isomorphism type of $M_\delta$.  
When is $\bold n_\lambda(T)$ equal to one?
\end{enumerate}
\mn
Now (see \cite{Sh:868}):
\mn
\begin{enumerate}
\item[$\odot_2$]  $\bold n_\lambda(T) = 1$ iff $T$ is superstable
\sn
\item[$\odot_3$]  for countable $T, \bold n_\lambda(T) = 2$ \underline{iff} $T$
is strictly stable (i.e. $T$ is stable, not superstable)
\sn
\item[$\odot_4$]  given an ordinal $\gamma$, for $\lambda$ large enough
$\bold n_\lambda(T) = |\gamma +1|$ \If \, $T$ is
stable and $\kappa(T) = \aleph_\gamma$ (recalling that for a stable
$T,\kappa(T)$ is cardinal $\le |T|^+$, so for countable $T$ it is
$\aleph_0$ or $\aleph_1$)
\sn
\item[$\odot_5$]   if $T$ is unstable, $\lambda = \aleph_\gamma$
then $\bold n_\lambda(T) \ge |\gamma +1|$.
\end{enumerate}
\mn
[Why?  Because for some closed unbounded subset $E$ of
$\lambda^+$, if $\delta \in E$ then $M_\delta$ is
$\cf(\delta)$-saturated but not $(\cf(\delta))^+$-saturated hence
$[\delta_1,\delta_2 \in E \wedge \text{ cf}(\delta_1) \ne
\text{ cf}(\delta_2) \Rightarrow M_{\delta_1} \ncong M_{\delta_2}$.]

\mn
Hence it is natural to replace $\bold n_\lambda(T)$ by:
\mn
\begin{enumerate}
\item[$\odot_6$]   let $\bold n_{\lambda,\kappa}(T) = 
\Min_E|\{M_\delta/\cong:\delta \in E$ and cf$(\delta) = \kappa\}|$ when
$\lambda > \kappa = \text{ cf}(\kappa)$ 
(as above $E$ varies on the clubs of $\lambda^+$).
\end{enumerate}
\mn
Below we use $\bold n_{\lambda,\kappa}(T)$ only when $\lambda =
\cf(\lambda) > |T| +
\kappa \wedge \kappa = \text{ cf}(\kappa)$ and remember that for
simplicity we are assuming G.C.H.

Now (see \cite{Sh:868}):
\mn
\begin{enumerate}
\item[$\odot_7$]  if $T$ is stable then $\bold n_{\lambda,\kappa}(T) = 1$.
\end{enumerate}
\mn
It is natural to ask whether this characterizes stable theories.  The
answer is no, in fact, by an example everyone knows
(by \cite[\S1]{Sh:877}):
\mn
\begin{enumerate}
\item[$\odot_8$]   $\bold n_{\lambda,\kappa}(T)=1$
for $T = \text{ Th}(\bbQ,<)$, the theory of dense linear orders with
neither first nor last element, so $\lambda
= \lambda^{< \lambda} > \kappa = \text{ cf}(\kappa)$.
\end{enumerate}
\mn
During the proof we analyze $p \in \bold S(M_\alpha),M_\alpha$
saturated, of course, only when $p \ne \tp(a,M_\alpha,M_\alpha)$ for
$a \in M_\alpha$.  So $M_\alpha$ is a linear order and $p$
induces a cut $(C^-_p,C^+_p)$ of $M_\alpha$, i.e. $C^-_p = \{a \in
M_\alpha:(a < x) \in p\}$ is an initial segment of $M_\alpha$ and its
compliment, $\{a \in M_\alpha:(a < x) \notin p\}$ is an end segment.  
This gives a pair of cofinalities, $(\mu^-_p,\mu^+_p),\mu^-_p$ the
cofinality of the linear order $C^-_p$ and $\mu^+_p$ the cofinality of
the inverse of $C^+_p$.  

Now
\mn
\begin{enumerate}
\item[$(*)_{8.1}$]  if $\mu_p := \text{ min}\{\mu^-_p,\mu^+_p\} <
\lambda$, \then \, the type is determined by any subset of $M_\alpha$ 
of cardinality $\mu_p$ such that:
\sn
\begin{enumerate}
\item[$\bullet$]  the set is unbounded in $C^-_p$ if 
$\mu_p = \mu^-_p$ and 
\sn
\item[$\bullet$]  the set is 
unbounded from below in $C^+_p$ if $\mu_p = \mu^+_p$.
\end{enumerate}
\sn
\item[$(*)_{8.2}$]  if $\mu_p = \lambda$, and we expand $M_\alpha$ by
the (unary) relation $C^-_p$, we still get a saturated model.
\end{enumerate}
\mn
Next considering $\odot_7 + \odot_8$ you may think that 
for every $T$ we get $\bold n_{\lambda,\kappa}(T)=1$,
 but (\cite[2.3(2)]{Sh:877} implies directly that):
\mn
\begin{enumerate}
\item[$\odot_9$]   $\bold n_{\lambda,\kappa}(T) = \lambda^+$ if $T$
is Peano arithmetic
\end{enumerate}
\mn
moreover, this holds for quite many theories $T$ (by \cite[\S2]{Sh:877}):
\mn
\begin{enumerate}
\item[$\odot_{10}$]   $\bold n_{\lambda,\kappa}(T) = \lambda^+$ if
$T$ has the strong independence property (see \cite{Sh:72}, i.e. 
for some first order formula
$\varphi(x,y),\langle \varphi(M,a):a \in M\rangle$ is an independent
sequence of subsets of $M$, see Definition \ref{0n.22}).
\end{enumerate}
\mn
For me this rings a bell and strengthens a suspicion - maybe the
dividing line is $T$ independent/$T$ dependent, indeed (by
\cite[\S2]{Sh:877}):
\mn
\begin{enumerate}
\item[$\odot_{11}$]   $\bold n_{\lambda,\kappa}(T) = \lambda^+$ if
$T$ is independent, $\lambda$ a successor cardinal.
\end{enumerate}
\mn
We try here to address the complement, the structure side.  This calls for
analyzing appropriate $\prec$-increasing continuous sequence $\bar M
= \langle M_i:i \le \kappa\rangle$ of models of $T$ of cardinality
$\lambda$.  Clearly in the relevant cases they ``increase fast enough" 
and $M_i$ is
saturated for $i$ non-limit.  Now among such sequences, is it not
reasonable to first deal with the case of length 2?

This leads to the generic pair conjecture which says that for $\lambda
 = \lambda^{< \lambda} > |T|$, we have $T$ is
independent iff $\bold n_{\lambda,2}(T) = \lambda^+$ where:  
\mn
\begin{enumerate}
\item[$\odot_{12}$]   $\bold n^*_{\lambda,2}(T) 
:= \text{ Min}_E|\{(M_\beta,M_\alpha)/\cong:\alpha < \beta$ belongs to $E$ and
cf$(\alpha) = \lambda = \text{ cf}(\beta)\}|$.
\end{enumerate}
\mn
Note that in defining $\bold n_{\lambda,\kappa}(T),\kappa \in
\text{ Reg } \cap [\aleph_0,\lambda]$ we speak on models of $T$,
i.e. $\delta \in E,\cf(\delta) = \kappa$ whereas here we deal 
with pairs of models.  However, to analyze $M_\delta$ for $\delta \in
E \wedge \cf(\delta) = \kappa,E$ small enough club of $\lambda^+$, it
is natural to assume $\delta = \sup\{\alpha \in E:\cf(\alpha) =
\lambda$ and $\alpha < \delta\}$ and choose $\bar \alpha \in
\Seq_{E,\kappa,\delta}$ which means $\bar\alpha$ is an 
increasing continuous sequence $\langle \alpha_i:i < \kappa\rangle$ of
ordinals with limit $\delta$ such that $i < \kappa$ non-limit
$\Rightarrow \cf(\alpha_i) = \lambda$.  So a sufficient condition for
$\bold n_{\lambda,\kappa}(T)=1$ is $\bold n^*_{\lambda,\kappa}(T)=1$
where $\bold n^*_{\lambda,\kappa} =
  \Min_E|\{M_{\bar\alpha}/\cong:\bar\alpha \in
\Seq_{E,\kappa,\delta}\}|,E$ varying on the clubs of $\lambda^+$.
  Now though it is not clear if this is also a necessary condition it
  seems more approachable and is natural.  Anyhow it seems reasonable
  to consider $\bold n^*_{\lambda,2}(T)=1$, i.e. the generic pair conjecture.

This connects us to the long term goal of classifying first order
theories by ``good" dividing lines, ones in which we find outside
properties (like here investigating $\bold n_{\lambda,\kappa}(T)$ or
just $\bold n_{\lambda,\lambda}(T)$,
trying to characterize it) with ``inside" definitions (like being
dependent), and developing an inside theory; here - looking at
decomposition (in \S1 decompositions of models, 
in \S2 decomposition of types, in \S3,\S4 strict decomposition of types).
More fully, for this we have to analyze types.  In \S1 we make a first
attempt; more exactly see \ref{3m.3} and \ref{3m.4}. 
We try to analyze a model 
$N := M_\beta$ over $M := M_\alpha$ by trying to find models 
$M_1,,N_1$ such that:
\mn
\begin{enumerate}
\item[$\boxplus_1$]   $M_\alpha = M \prec M_1 \prec N_1$ and $M_\alpha
= M \prec N = M_\beta \prec N_1$
\sn
\item[$\boxplus_2$]   for every $\bar a \in {}^{\omega >}(M_1)$ for some
$B_\alpha \in [M]^{< \lambda}$ the type 
tp$(\bar a,M_\beta,M_1)$ is definable over $B_\alpha$ in a weak
sense, i.e. does not split over $B_\alpha$,
this means that if $n < \omega$ and $\bar b,\bar c \in
{}^n(M_\alpha)$ realizes the same type over $B_\alpha$ \then \, so does
$\bar a \char 94 \bar b,\bar a \char 94 \bar c$ 
 (this is parallel to $(*)_{8.1}$ from $\odot_8$); it follows that for
any sequence $\bar a \in {}^{\kappa >}(M_1)$ a similar statement holds
\sn
\item[$\boxplus_3$]  tp$(N_1,M_1,N_1)$ is weakly orthogonal to every $q
\in \bold S^{< \omega}(M_1)$ which does not split over some $B \in
[M_1]^{< \lambda}$; the weakly orthogonal means that $q$ has a unique
extension in $\bold S^n(N_1)$ wherever $q \in \bold S^n(M_1)$.
\end{enumerate}
\mn
In \S2 we try to analyze a type rather than a pair of models, also we
find it better to deal with $\theta$-types, $\theta \ge |T|$, as during
the analysis we add more variables.  So for a $\kappa$-saturated model 
$M \prec \gC$ and sequence $\bar d$ of length $< \theta^+$ we try to 
analyze tp$(\bar d,M,\gC)$ in
two steps.  The first is to add $\bar c$ of length $< \theta^+$ such that
\mn
\begin{enumerate}
\item[$\boxplus_4$]   tp$(\bar c,M,\gC)$ does not split over some $B
\subseteq M \prec \gC$ of cardinality $< \kappa$.
\end{enumerate}
\mn
This corresponds to the stable type (``unfortunately" but unavoidably 
depending on $\kappa$),  so for the theory of dense linear
orders it corresponds to types $p \in \bold S(M)$ with $\mu_p <
\kappa$, see $(*)_{8.1}$ above.  
True, they are not really definable, but non-splitting is a
weak form of being definable.  The second step is
\mn
\begin{enumerate}
\item[$\boxplus_5$]   tp$(\bar d,M + \bar c,\gC)$ is tree like, i.e. if $A
\subseteq M \prec \gC$ and $|A| < \kappa$
\then \, for some $\bar e \in {}^{\theta^+ >} M$
we have tp$(\bar d,\bar c + \bar e) \vdash \text{ tp}(\bar d,A + c)$.
\end{enumerate}
\mn
This property holds for $T = \text{ Th}(\Bbb Q,<),p \in \bold S(M)$ when
$\mu_p \ge \kappa$!, i.e. when both cofinalities are $\ge \kappa$.  This
is the Type Decomposition Theorem (\ref{tp25.43}).

A consequence is some clarification of models of $M$ of a dependent
theory which are exactly $\kappa$-saturated for singular $\kappa$.  We
deal with this question to some extent per se.

In \S3 we get a better decomposition - strict decomposition.  But at
present with a price, assuming e.g. $\kappa = \|M\|$ is a measurable
cardinal.  The main
point appears in \S4, the existence of limits of increasing sequences of strict
decompositions.  

Using this we are able to prove the pair genericity 
conjecture, the structure side
for the case of a measurable cardinal.  The measurability assumption
seems undesirable.  Describing this to Udi Hrushovski he was more
concerned about also having the  non-structure side for independent
$T$.  Now at the time in \cite{Sh:877} it was remarked that a similar
proof should work for the strongly inaccessibles, but the author 
was not motivated enough to really look into it.  
Subsequently \cite{Sh:906} completes it.

The order of the sections is by their conceptions, 
so there are some repetitions.
In \cite{Sh:950} and Kaplan-Shelah \cite{KpSh:946} we start to continue this
work as well as in  \cite{Sh:F1127}.  Note that \cite{Sh:950}
concentrate on saturated models but it works just as well for special
models (in singular strong limit cardinals, see e.g. \cite{CK73}).

We thank the referee with thoroughness much above the call of duty
causing the paper to be much improved 
and John Baldwin for much helpful criticism.

\begin{context}
\label{0n.1}  
1) $T$ is complete first order theory.

\noindent
2) ${\frak C} = {\frak C}_T$ is a monster model for $T$, omitting $T$ when
   no confusion arises; i.e. $\bar\kappa$ is a large enough cardinal,
$\gC$ is a $\bar\kappa$-saturated model such that we deal only with
models $M \prec \gC$, sets $A \subseteq \gC$ of cardinality 
$< \bar\kappa$ and sequences $\bar a,\bar b,\bar c,\bar d,\bar e$ 
from ${}^\alpha \gC$ for some $\alpha < \bar\kappa$.  So tp$(\bar c,A)$ means 
tp$(\bar c,A,\gC)$.  

\noindent
3) We may not pedantically distinguish a model $M$ and its universe,
   the cardinality $\|M\|$ of $M$ is that of its universe.
\end{context}

\begin{notation}
\label{0n.4}  
1) For $M \prec {\frak C}$ and $\bar a \in {}^\alpha M$ or just $\bar
a \in {}^\alpha \gC$ let $M_{[\bar a]}$ be the expansion 
of $M$ by every relation $R_{\varphi(\bar x,\bar a)}
= \varphi(M,\bar a)$ where $\varphi(M,\bar a) :=
\{\bar b \in {}^{\ell g(\bar x)}M:{\frak C} \models \varphi[\bar
b,\bar a]\}$ for $\varphi(\bar x,\bar y) \in \Bbb L(\tau_T)$ such that
$\ell g(\bar y) = \alpha,\ell g(\bar x) < \omega$ or pedantically
$\varphi(\bar x,\bar y \restriction u)$ for $\bar x,\bar y$ as above,
$u \subseteq \alpha$ finite.  We define $M_{[A]}$ similarly, i.e. as
the expansion of $M$ by $R_{\varphi(\bar x,\bar a)} = \varphi(M,\bar
a)$ for every $\in\, {}^{\ell g(\bar y)}A$ and $\varphi(\bar x,\bar
y) \in \bbL(\tau_T)$.

\noindent
1A) For $p(\bar x) \in \bold S^\alpha(M)$ let $M_{[p]}$ be $M_{[\bar
    a]}$ whenever $\bar a \in {}^\alpha \gC$ realizes $p(\bar x)$.

\noindent
1B) We say the sequence $\langle \varphi_s(\bar x,\bar a_s):s \in
I\rangle$ of formulas from $\bbL(\tau_M)$ with $\bar a_s$ from $M$ is
independent in the model $M$ \when \, every finite non-trivial Boolean
combination of sets from $\varphi_s(M,\bar a_s)$ is non-empty.

\noindent
2) Writing $\varphi(\bar x,\bar y) \in \Bbb L(\tau_T),\varphi$ here is
always first order but $\bar x$ and $\bar y$ may be infinite, though
sometimes are finite (said or clear from the context).  Let $p(\bar
x),q(\bar x),r(\bar x)$ denote types over some $A \subseteq \gC$,
i.e. set of formulas of the form $\varphi(\bar x,\bar b),\bar b \in
{}^{(\ell g(\bar b))}A$.

\noindent
3) EC$_\lambda(T)$ is the class of models $M$ of $T$ (so $M \prec
{\frak C}$) of cardinality $\lambda$ and $\EC_{\lambda,\kappa}(T)$ is
the class of $\kappa$-saturated $M \in \EC_\lambda(T)$.

\noindent
4) $A + \bar c$ is $A \cup \text{ Rang}(\bar c)$, etc.

\noindent
5) Let $\tp(A,B)$ be $\tp(\bar a,B)$ where $\bar a$ is the identity
function on $A$.
\end{notation}

\begin{definition}
\label{0n.8}  
1)  If $\bar a_t \in {}^\gamma {\frak C}$ for $t
\in I$ and ${\cD}$ is a filter on $I$ and $\bar x = \langle x_i:i <
\gamma\rangle$ and $A \subseteq {\frak C}$ then Av$(\langle \bar a_t:t
\in I\rangle/{\cD},A) = 
\{\varphi(\bar x,\bar b):\bar b \in {}^{\omega >} A$
and the set $\{t \in I:{\frak C} \models \varphi[\bar a_t,\bar b]\}$
belongs to ${\cD}\}$. 
Note that if $T$ is dependent, $I$ is a linear order with no last
members and $\langle \bar a_t:t \in I\rangle$ is an indiscernible
sequence, see below \then \, the result $\in \bold S^\gamma(A)$.  Also
note that if $\cD$ is an ultrafilter on $I$ then Av$(\langle \bar
a_t:t \in I\rangle/\cD,A)$ belongs to $\bold S^\gamma(A)$.

\noindent
1A)  Recall that if ${\cD}$ is a filter on 
$\{\bar a_t:t \in I\} \subseteq {}^\alpha\gC$ and $A \subseteq \gC$ 
we define $\Av({\cD},A)$ similarly and if $I$ a linear order and
${\cD}$ is the filter of co-bounded subsets of $I$ we may omit it.

\noindent
2) If $p(\bar x),q(\bar y)$ are complete types over $A$ we say $p(\bar
x),q(\bar y)$ are weakly orthogonal \when \, for every $\bar a_1,\bar a_2$
realizing $p(\bar x)$ and $\bar b_1,\bar b_2$ realizing $q(\bar y)$ we
have tp$(\bar a_1 \char 94 \bar b_1,A) = \text{ tp}(\bar a_2 \char 94
\bar b_2,A)$.

\noindent
3) For a linear order $I,\langle \bar a_s:s \in I\rangle$ is an
   indiscernible sequence over $B$ \when \,: $\ell g(\bar a_s)$ is constant
   and if $s_0 <_I \ldots <_I s_{n-1}$ an $t_0 <_I \ldots <_I t_{n-1}$
then the sequences $\bar a_{s_0} \char 94 \ldots \char 94 
\bar a_{s_{n-1}}$ and $\bar a_{t_0} \char 94 \ldots 
\char 94 \bar a_{t_{n-1}}$ realize the same type over $B$.
\end{definition}

Recall also (see \cite[Ch.II,\S4]{Sh:c})
\begin{fact}
\label{0n.17}  If $T$ is dependent \then \, for any formula
$\varphi = \varphi(\bar x,\bar y,\bar z) \in \bbL(\tau_T)$ there is
$n = n_\varphi < \omega$ (depending on $T$) such that:
\mn
\begin{enumerate}
\item[$(a)$]   for no $\bar c \in {}^{\ell g(\bar z)} {\frak C}$ and $\bar
b_i \in {}^{\ell g(\bar y)} {\frak C}$ for $i < n$, is the sequence $\langle
\varphi(\bar x,\bar b_i,\bar c):i < n \rangle$ independent, i.e. every
non-trivial Boolean combination of the sets $\varphi(M,\bar b_i,\bar c) 
= \{\bar a \in {}^{\ell g(\bar x)} M:M \models 
\varphi[\bar a,\bar b_i,\bar c]\}$ for $i < n$ is non-empty.
\sn
\item[$(b)$]   If $\langle \bar b_i:i < n\rangle$ is an
indiscernible sequence over $C,\ell g(\bar b_i) = \ell g(\bar y),\bar c
\in {}^{\ell g(\bar z)} C$ (all in ${\frak C}$) then for no $\bar a
\in {}^{\ell g(\bar a)} M$ do we have ${\frak C} \models \varphi[\bar
a,\bar b_i,\bar c]^{\text{if}(\ell\text{ even})}$ for $\ell < n$.
\sn
\item[$(c)$]   Also there is a finite $\Delta_\varphi \subseteq 
\bbL(\tau_T)$ such that in clause (b) it is enough to demand that $\langle
 \bar b_i:i < n\rangle$ is a $\Delta$-indiscernible sequence.
\end{enumerate} 
\end{fact}

\noindent
Lastly, we quote Erd\"os-Rado \cite{ErRa69}.
\begin{fact}
\label{0n.19}
The $\Delta$-System Lemma for finite sets.

For every natural numbers $k,n$ there is a natural number $m$ such
that: if $u_i$ is a finite set with $\le k$ elements for $i < m$ \then \,
there are sets $w \subseteq \{0,\dotsc,m-1\}$ with $|w|=n$ and $u_*$
such that $\langle u_i:i \in w\rangle$ is a $\Delta$-system with heart
$u_*$, which means that $i \ne j \in w \Rightarrow u_i \cap u_j = u_*$.
\end{fact}

\begin{definition}
\label{0n.22}
1) A partial order $I$ is $\kappa$-directed \when \, every set $J
\subseteq I$ of cardinality $< \kappa$ has an upper bound $t \in I$
which means that $(\forall s)(s \in J \Rightarrow s \le_I t)$.

\noindent
2) A sequence $\langle A_s:s \in I\rangle$ is an independent sequence
   of subsets of $A_*$ \when \, ($A_s \subseteq A_*$ for $s \in I$
   and) $\bigcap\limits_{s \in u} A_s \backslash \bigcup\limits_{t \in
   v} A_t$ is non-empty for every disjoint finite $u,v \subseteq I$.
\end{definition}
\newpage

\section {Non-splitting constructions} \label{s:non}

On such constructions including $\bold F^{\text{nsp}}_\kappa$
see \cite[Ch.IV,\S1,\S3]{Sh:c} but 
$\bold F^{\text{nsp}}$ here is $\bold F^p$ there; and see
\cite[4.23-4.26]{Sh:715}, however this section is self-contained. 

We try here to analyze $\kappa$-saturated models $M \prec N$, e.g. by finding
$M_1,N_1$ such that $M \prec M_1 \prec N_1,N \prec N_1$ where $M_1$ is
$\bold F^{\text{nsp}}_\kappa$-constructible over $M$, see below and
tp$(N_1,M_1)$ is weakly orthogonal to any type over $M_1$ realized in
some $\bold F^{\text{nsp}}_\kappa$-construction over it, see Theorem
\ref{3m.4}, part (B) noting that $M,N,N_1,M_1$ here stands for
$A,A^+,M,N$ there.  We first
recall the definition of non-splitting and some of its properties.

\begin{definition}
\label{3k.0.4}  
We say $p(\bar x)$ does not
split over $A$ \when \,: if $\varphi(\bar x,\bar b),\neg \varphi(\bar
x,\bar c) \in p(\bar x)$ then $\tp(\bar b,A) \ne \text{ tp}(\bar c,A)$.
\end{definition}

\begin{fact}
\label{3k.0.7}
1) If $\circledast_{A,B,C}$ below holds and $p(\bar x) 
\in \bold S^m(B)$ does not split over $A$,
\then \, there is one and only one $q(\bar x) \in \bold S^m(C)$
extending $p(\bar x)$
and not splitting over $A$ (also called the non-splitting extension of
$p(\bar x)$ over $C$), where:
\mn
\begin{enumerate}
\item[$\circledast_{A,B,C}$]  $(a) \quad A \subseteq B \subseteq C$
\sn
\item[${{}}$]   $(b) \quad$ for every $\bar c \in {}^{\omega >}C$
there is $\bar b \in {}^{\ell g(\bar c)}B$ realizing tp$(\bar c,A)$.
\end{enumerate}
\mn
2) Let $I$ be a linear order.
If tp$(\bar a_t,B \cup \bigcup\{\bar a_s:s <_I t\})$ does not split
over $B$ and increases with $t \in I$ \then \, 
$\langle \bar a_t:t \in I\rangle$ is an indiscernible sequence over $B$.

\noindent
3) If $\tp(\bar a,B)$ does not split over $A$, the sequence $\langle
\bar b_t:t \in I\rangle$ is an indiscernible sequence over $A$
and
$\bar b_t \subseteq B$ for $t \in I$ \then \, $\langle \bar b_t:t \in
I\rangle$ is an indiscernible sequence over $A \cup \bar a$.

\noindent
4)  If $A \subseteq B$ \then \, the number of $p \in \bold S^\theta(B)$
which does not split over $A$ is $\le 2^{2^{|A|+|T|+\theta}}$,
moreover if $T$ is dependent the number is $\le 2^{|A|+|T|+\theta}$.

\noindent
5) If $A \subseteq B$ and $p \in \bold S^\alpha(B)$ is finitely
   satisfiable in $A$ \then \, $p$ does not split over $A$.
\end{fact}

\begin{PROOF}{\ref{3k.0.7}}
1) By \cite{Sh:3} or see \cite{Sh:715} or see
\cite[1.10]{Sh:300a} for uniqueness.

\noindent
2) By \cite[I]{Sh:c} or \cite[I]{Sh:300} or \cite[3.2]{Sh:300a}.

\noindent
3)  By the definitions.

\noindent
4) The first conclusion is easy and see \cite{Sh:3} or
\cite[\S1]{Sh:300a}, the second holds by \cite[5.26]{Sh:783}.

\noindent
5) Easy, too.
\end{PROOF}

\begin{fact}
\label{3k.2}
[Assume $T$ is dependent.]
  
If $p(\bar x)$ is an $\alpha$-type over $B \subseteq A$ \then \,
we can find $q(\bar x) \in \bold S^\alpha(A)$ 
extending $p(\bar x)$ such that for some $C \subseteq A$
of cardinality $\le |T| + |\alpha|$ 
the type $q(\bar x)$ does not split over $B \cup C$.
\end{fact}

\begin{PROOF}{\ref{3k.2}}
\cite[III,7.5,pg.140]{Sh:c} or see \cite[4.24]{Sh:715}.
\end{PROOF}

\begin{observation}
\label{3k.4}
For $\kappa$ regular.

\noindent
1) If $A \subseteq B,|A| < \kappa$ and $\bar a \in 
{}^{\kappa >}{\frak C}$ and tp$(\bar a,B)$ is finitely satisfiable in 
$A$ \then \, it does not split over $A$.

\noindent
2) If $A \subseteq B,\bar c \in {}^{\kappa >}{\frak C}$ and 
tp$(\bar c,B)$ does not split over $A$ and $i < \ell g(\bar c)$
\then \, tp$(c_i,B \cup \{c_j:j <i\})$ does not split over $A \cup
\{c_j:j<i\}$.  Similarly for $\langle \bar c_j:j < i \rangle$ when $j \le
i \Rightarrow \Rang(\bar c_j) \subseteq \Rang(\bar c)$.

\noindent
3) If tp$(\bar c_k,B + \bar c_0 + \ldots + \bar c_{k-1})$ does not
split over $A \subseteq B$ for $k < n$ \then \, tp$(\bar c_0 \char
94 \ldots \char 94 \bar c_{n-1},B)$ does not split over $A$.

\noindent
4) If $A \subseteq A_1 \subseteq B_1 \subseteq B$ and Rang$(\bar c_1)
\subseteq \text{ Rang}(\bar c)$ and tp$(\bar c,B)$ does not split
over $A$ \then \, tp$(\bar c_1,B_1)$ does not split over $A_1$.

\noindent
5) If $\bar c \in {}^{\kappa >}\gC$ and for every finite $u \subseteq
   \ell g(\bar c)$ and finite $B_1 \subseteq B$ the type $\tp(\bar c
\rest u,B_1)$ does not split over $A$ \then \, $\tp(\bar c,B)$ does
   not split over $A$.
\end{observation}

\begin{PROOF}{\ref{3k.4}}
Easy.
\end{PROOF}

\noindent
As in \cite[IV]{Sh:c}.
\begin{definition}
\label{3k.6} 
1) ${\cA}$ is an $\bold F^{\text{nsp}}_\kappa$-construction \when \,:
\mn
\begin{enumerate}
\item[$(a)$]   ${\cA} = (A,\bar a,\bar B,\bar A,\alpha) =
(A^{\cA},\bar a^{\cA},\bar B^{\cA},\bar A^{\cA},\alpha^{\cA})$,
\sn
\item[$(b)$]  $\bar a = \langle a_\beta:\beta < \alpha\rangle =
\langle a^{\cA}_\beta:\beta < \alpha\rangle$,
\sn
\item[$(c)$]   $\bar B = \langle B_\beta:\beta < \alpha\rangle =
\langle B^{\cA}_\beta:\beta < \alpha\rangle$,
\sn
\item[$(d)$]  $\bar A = \langle A_\beta:\beta \le \alpha\rangle =
\langle A^{\cA}_\beta:\beta \le \alpha\rangle$,
\sn
\item[$(e)$]   $A_\beta = A \cup\{a_\gamma:\gamma < \beta\}$,
\sn
\item[$(f)$]   $B_\beta \subseteq A_\beta$ and $|B_\beta| < \kappa$,
\sn
\item[$(g)$]  tp$(a_\beta,A_\beta)$ does not split over $B_\beta$.
\end{enumerate}
\mn
2) We let $\ell g({\cA}) = \alpha^{\cA}$ and writing ${\cA}$
we may omit $\bar A^{\cA},\alpha^{\cA}$ as they are
determined by the others so may write ${\cA} = (A,\bar a,\bar B)$ or
${\cA} = (A,\langle (a_\beta,B_\beta):\beta < \alpha\rangle)$.  We
may replace $a_\beta$ by a finite sequence $\bar a_\beta$ with no real change.

\noindent
3) We say the $\bold F^{\text{nsp}}_\kappa$-construction 
${\cA}$ is $\mu$-full \when \,
cf$(\ell g({\cA})) \ge \kappa$ and if $q \in 
\bold S(A^{\cA}_{\ell g({\cA})})$ does not split 
over $B$ where $B \subseteq A^{\cA}_{\ell g({\cA})}$ 
has cardinality $< \kappa$, \then \, $\{\beta
< \ell g({\cA}):a_\beta$ realizes $p \restriction A^{\cA}_\beta$
and $B \subseteq A^{\cA}_\beta\}$
is unbounded in $\alpha^{\cA}$ and has order type divisible by $\mu$.

\noindent
4) We say $C$ is $\bold F^{\text{nsp}}_\kappa$-constructible over $A$ 
\when \, there is an $\bold F^{\text{nsp}}_\kappa$-construction
${\cA}$ such that $A = A^{\cA} = A^{\cA}_0$ and $C = 
A^{\cA}_{\ell g({\cA})}$. 
\end{definition}

\begin{definition}
\label{3m.1}   
1) Let $A \le_\kappa C$ mean that $C$ is 
$\bold F^{\text{nsp}}_\kappa$-constructible over $A$.

\noindent
2)  We say that $(A^+,A)$ is $\kappa$-reduced \when \,:  
if $A \le_\kappa C$ and $\bar c \in {}^{\kappa >}(A^+)$ then 
tp$(\bar c,A)$ has a unique extension to a complete type over $C$.

\noindent
3) We say the $(N,M)$ is $\kappa$-nice \when \,:
\mn
\begin{enumerate}
\item[$(a)$]   $(N,M)$ is $\kappa$-reduced and $M \prec N$,
\sn
\item[$(b)$]   $M$ is $\kappa$-saturated,
\sn
\item[$(c)$]   $N$ is $\kappa$-saturated,
\sn
\item[$(d)$]   if $M \le_\kappa M^+$ then $M_{[N]} \prec M^+_{[N]}$,
see below.
\end{enumerate}
\mn
3A) Recall $M_{[B]}$ is $M$ expanded by $R_{\varphi(\bar x,\bar a)} =
\{\bar b \in {}^{\ell g(\bar x)}M:{\frak C} \models \varphi[\bar
b,\bar a]\}$ for $\varphi(\bar x,\bar y) \in \bbL(\tau_T)$ (with $\bar
x$ finite of course), $\bar a \in {}^{\ell g(\bar y)} N$ 
and recall Th$(M_{[B]})$ is dependent by \cite[\S1]{Sh:783}.

\noindent
4) We say that $(M,A)$ is pseudo $\kappa$-reduced when: 
if $\bar c \in {}^{\omega >}A,\|M_1\| < \kappa,M_1 
\subseteq M,q(\bar x) \in \bold S^{< \omega}(M)$ is 
finitely satisfiable in $M_1$ \then \, $q(\bar x)$, 
tp$(\bar c,M)$ are weakly orthogonal.
\end{definition}

\begin{observation}
\label{3m.2}  
For $\kappa$ regular:

\noindent
1) $\le_\kappa$ is a partial order.

\noindent
2) If $\langle A_i:i \le \alpha\rangle$ is increasing continuous and
$i < \alpha \Rightarrow A_i \le_\kappa A_{i+1}$ \then \, $A_0 \le_\kappa
A_\alpha$.

\noindent
3) In Definition \ref{3m.1}(2) it is enough to consider $\bar c \in
 {}^{\omega >}(C \backslash A)$. 

\noindent
4) If $A \le_\kappa B$ and $\bar c \in {}^{\kappa >}B$ 
\then\footnote{we could have chosen this as the definition.  This changes
the places we really need ``$\kappa$ regular".} \,
$\tp(\bar c,A)$ does not split over some $A' \subseteq A$ of cardinality $<
\kappa$.

\noindent
5) If the pair $(M,C)$ is $\kappa$-reduced \then \, $(M,C)$ 
is pseudo $\kappa$-reduced.

\noindent
6) If $\tp(\bar a,A)$ does not split over $B$ and $B \subseteq A$ has
   cardinality $< \kappa$ \then \, $A \le_\kappa A + \bar a$.
\end{observation}

\begin{PROOF}{\ref{3m.2}}  
Easy; e.g. part (6) by \ref{3k.4}(2) and part (4) by
\ref{3k.2}(3),(4)) and part (5) by \ref{3k.0.7}.
\end{PROOF}

\begin{claim}
\label{3m.3}   [$T$ is dependent and $\kappa =
 \,\text{\rm cf}(\kappa) > |T|$].

\noindent
1) For every $A$ there is a  $\kappa$-saturated $C$ 
such that $A \le_\kappa C$ and $|C| \le (|A| + |T|)^{< \kappa}$.

\noindent
2) If in addition $\mu \le (|A| + |T|)^{< \kappa}$ \then \, we can add
``$C$ is ``$\mu$-full $\kappa$-saturated"; 
clearer if $|C| \le (|A| +  |T|)^{< \kappa} + 2^{2^\kappa}$).
\end{claim}

\begin{PROOF}{\ref{3m.3}}
1) By \ref{3k.2} + \ref{3m.2}(2) and \ref{3m.2}(6). 

\noindent
2) Similarly (by \ref{3k.0.7}(4)).
\end{PROOF}

Now we arrive to the first result giving a decomposition.  The type
tp$(A^+,A)$ is decomposed in Theorem \ref{3m.4} by finding $M$ such that
$A \le_\kappa M$, (so the complete types over $A$ realized in $M$ are
somewhat definable) and $(A^+,M)$ is $\kappa$-reduced, so the type
tp$(A^+,M)$ is weakly orthogonal to types in 
$\bold S^{< \omega}(M)$ not splitting over subsets of $M$ of cardinality $<
\kappa$. 

\begin{theorem}
\label{3m.4}   
\underline{The Density of Reduced Pairs Theorem}
[$T$ dependent].

For any $A \subseteq A^+$ and $\kappa = \text{\rm cf}(\kappa) > |T|$
and $\lambda$ satisfying\footnote{no real harm if we replace ``$\theta
< \lambda \Rightarrow \lambda = \lambda^\theta \ge |A^+|"$ by $\theta
< \kappa \Rightarrow \lambda = \lambda^\theta + 2^{2^{\theta+|T|}} \ge
|A^+|$ and then we can use only the first version of \ref{3k.0.7}(4).}
$\theta < \kappa \Rightarrow \lambda = \lambda^\theta \ge |A^+|$ 
\mn
\begin{enumerate}
\item[$(A)$]   we can find $M$ such that $M$ is a model of 
cardinality $\lambda$ such that $A \le_\kappa M$ and $(A^+,M)$ 
is $\kappa$-reduced
\sn
\item[$(B)$]  for some $M$ as in clause (A) and $N$ the pair
$(N,M)$ is $\kappa$-reduced and even $\kappa$-nice  and $A^+ \subseteq N$.
\end{enumerate}
\end{theorem}

\begin{PROOF}{\ref{3m.4}}
\underline{Proof of (A)}:  

Our intention is to try to do a construction as described in
$\boxtimes$ below.  Having carried the induction the proof is divided
to two cases.   In the first we get the desired conclusion.  In the second, we
get a contradiction to $T$ being dependent; formally to the
maximality of the $k$ chosen in $(g)(\beta)$ of $\boxtimes$.

We choose $M_i,B_i,j_i,\bar c_i$
by induction on $i \le \lambda^+$ such that
\mn
\begin{enumerate}
\item[$\boxtimes$]   $(a) \quad M_i$ is $\prec$-increasing
continuous, $M_i$ of cardinality $\le \lambda + |i|$,
\sn
\item[${{}}$]  $(b) \quad j_i \le i,B_i \prec M_i,\|B_i\| < \kappa$,
\sn
\item[${{}}$]  $(c) \quad M_0$ is $\bold
F^{\text{nsp}}_\kappa$-constructible over $A$,
\sn
\item[${{}}$]  $(d) \quad M_{i+1}$ is 
$\bold F^{\text{nsp}}_\kappa$-constructible over $M_i$ and
$M_{i+1}$ is $\kappa$-saturated,
\sn
\item[${{}}$]   $(e) \quad \bar c_i \in {}^{\omega >}(M_{i+1})$ and
$B_i \subseteq M_{j_i}$ has cardinality $< \kappa$,
\sn
\item[${{}}$]  $(f) \quad$tp$(\bar c_i,M_i)$ does not split over $B_i$,
\sn
\item[${{}}$]  $(g) \quad$ if cf$(i) \ge \kappa$ and subclause
$(\alpha)$ below holds then subclause $(\beta)$ holds 

\hskip25pt where:

\underline{Subclause $(\alpha)$}: There are 
$j < i,m < \omega,B \subseteq M_j$ of cardinality $< \kappa$ and 

$p(\bar x) \in {\bold S}^m(M_i)$ 
which does not split over $B$ and $p(\bar x)$ has 

$\ge 2$ extensions in $\bold S^m(M_i \cup A^+)$.

\underline{Subclause $(\beta)$}:  There are $m = 
m_i < \omega,j = j_i < i,k = k_i < \omega$ and 

$\varphi(\bar x,\bar y) = \varphi_i(\bar x_i,\bar y_i) \in 
\bbL(\tau_T)$ with $\bar x,\bar y$ finite 
and $\bar b =$

$\bar b_i \in {}^{(\ell g(\bar y))}(M_j \cup A^+)$ 
and $\varepsilon_0 < \ldots < \varepsilon_{k-1}$ from the
interval 

$[j_i,i)$ such that:
\mn
\begin{enumerate}
\item[$\bullet$]  $B_i \subseteq M_j$,
\sn
\item[$\bullet$]  $\tp(\bar c_i,M_i) \cup \{\varphi_i(\bar x,\bar
b_i)\}$ and $\tp(\bar c_i,M_i) \cup \{\neg 
\varphi_i(\bar x,\bar b_i)\}$ are consistent,
\sn
\item[$\bullet$]  $\tp(\bar c_{\varepsilon_\ell},M_{\varepsilon_\ell})
= \text{ tp}(\bar c_i,M_{\varepsilon_\ell})$,
\sn
\item[$\bullet$]   ${\frak C} \models \varphi[\bar c_{\varepsilon_\ell},
\bar b]^{\text{if}(\ell\text{ even})}$,
\sn
\item[$\bullet$]   $k$ is maximal for the given 
$\varphi(\bar x,\bar y),\bar b,j_i$ (see $\circledast$ below; 
$k$ is well defined as $T$ is dependent, see $\circledast$ below),
\sn
\item[$\bullet$]   ${\frak C} \models 
\varphi[\bar c_i,\bar b]^{\text{if}(k\text{ is even})}$. 
\end{enumerate}
\end{enumerate}
\mn
So in stage $i$ we first choose $M_i$: if $i=0$ by clause (c), such
model $M_i$ exists by
\ref{3m.3}(1), if $i$ is a limit ordinal we choose $M_i$ as 
$\cup\{M_j:j < i\}$ and if $i=j+1$ (so $\bar c_j$ has already been defined)
 then choose
$M_i$ such that $M_j \cup \bar c_j \le_\kappa M_{j+1}$ and $M_i =
M_{j+1}$ is $\kappa$-saturated of cardinality $\lambda$ 
(and $\lambda$-full if you like), possible by Claim \ref{3m.3}(1).
Note that $M_j \le_\kappa M_j \cup \bar c_j$ by clause (f) and \ref{3m.2}(6) 
hence $M_j \le_\kappa M_{j+1}$ recalling \ref{3m.2}(2).

Note
\mn
\begin{enumerate}
\item[$\circledast$]   there is $n=n_{\varphi(\bar x,\bar y)}$
depending on $\varphi(\bar x,\bar y)$ and $T$ only such that in subclause
$(\beta)$ we have $\varphi_i(\bar x_i,\bar y_i) 
= \varphi(\bar x,\bar y) \Rightarrow k_i \le n$.
\end{enumerate}
\mn
[Why?  As by clause (f) in subclause (g)$(\beta)$ of $\boxplus(f)$, 
by \ref{3k.0.7}(2) the sequence 
$\langle \bar c_{\varepsilon_\ell}:\ell < k\rangle$ is an 
indiscernible sequence, so by $T$ being dependent we are done by
\ref{0n.17}.] 

Second, why can we choose $(m_i,j_i,B_i,\varphi_i,\text{tp}(\bar c_i,M_i))$ 
as required in clause (g)?  If cf$(i) < \kappa$ or the antecedent of clause 
(g), i.e. $(g)(\alpha)$ fails \then \, trivially yes (choose e.g. $\bar c_i$ as
the empty sequence).   Otherwise let $j<i,B \subseteq M_j$ be of
cardinality $< \kappa,m < \omega$ and $p(\bar x) 
\in \bold S^m(M_i)$ which does not split over $B$ and 
which has extensions $p_0(\bar x)
\ne p_1(\bar x)$ in $\bold S^m(M_i \cup A^+)$ with $p_0 \restriction
(M_i \cup A^+) \ne p_1 \restriction (M_i \cup A^+)$, so $p_0 \restriction
M_i = p = p_1 \restriction M_i$ does not split over $B$.

Hence for some $\bar b \in {}^{\omega >}(M_i \cup A^+)$ and $\varphi =
\varphi(\bar x,\bar y) \in \Bbb L(\tau_T)$ we have $\varphi(\bar
x,\bar b) \in p_1(\bar x),\neg \varphi(\bar x,\bar b) \in p_0(\bar
x)$; as $i$ is a limit ordinal \wilog \, $\bar b \in {}^{\omega >}(M_j
\cup A^+)$.  We now try to choose $\varepsilon_\ell$ by induction on
$\ell \le n_{\varphi(\bar x,\bar y)}$ such that:
\mn
\begin{enumerate}
\item[$\odot$]   $(a) \quad j \le \varepsilon_\ell < i$ and $k < \ell
\Rightarrow \varepsilon_k < \varepsilon_\ell$,
\sn
\item[${{}}$]  $(b) \quad \bar c_{\varepsilon_\ell}$ realizes
$p(\bar x) \restriction M_{\varepsilon_\ell}$,
\sn
\item[${{}}$]   $(c) \quad {\frak C} \models \varphi
[\bar c_{\varepsilon_\ell},\bar b]^{\text{if}(\ell\text{ is even)}}$,
\sn
\item[${{}}$]   $(d) \quad \varepsilon_\ell$ is minimal under (a)+(b)+(c).
\end{enumerate}
\mn
So by $\circledast$ for some $k \le n_{\varphi(\bar x,\bar y)}$ we
have: $\varepsilon_\ell$ is well defined iff $\ell < k$.  At last we
choose:
\mn
\begin{enumerate}
\item[$(*)$]    $(a) \quad B_i = B$
\sn
\item[${{}}$]    $(b) \quad \varphi_i = \varphi$
\sn
\item[${{}}$]    $(c) \quad j_i = j$
\sn
\item[${{}}$]    $(d) \quad k_i = k$
\sn
\item[${{}}$]   $(e) \quad \bar c_i$ realizes $p_1(\bar x)$ if $k$ is even
and realizes $p_0(\bar x)$ if $k$ is odd.
\end{enumerate}
\mn
So at last we are in a situation where the construction from
$\boxtimes(a)-(g)$ having been carried out.  So now comes the division
to cases.

Let $S_* = \{i < \lambda^+$: cf$(i) \ge \kappa$ and subclause
$(\alpha)$ of clause (g) holds for $i\}$.

Recall that $S$ is a stationary subset of an ordinal $\delta$ of
cofinality $> \aleph_0$ (e.g. a regular uncountable cardinal) when it
is not disjoint to any closed unbounded subset $E$ of $\delta$.
\bigskip

\noindent
\underline{Case 1}:  $S_*$ is a stationary subset of $\lambda^+$.  

Hence for $i \in S_*$, there are $j_i,B_i,
\varphi_i(\bar x_i,\bar y_i),\bar b_i$ and $k_i < \omega$ and
$\varepsilon_0(i) < \ldots < \varepsilon_{k_i-1}(i) < i$ as in
subclause $(\beta)$ of $\boxtimes(g)$ and by Fodor's
lemma (see e.g. \cite{J}) 
for some $m_* < \omega,j < \lambda^+,B,\varphi(\bar x,\bar
b),k_*,\langle \varepsilon_i:i < k^*\rangle$ and a stationary subset
$S$ of $S_* = \{\delta < \lambda^+$: cf$(\delta) \ge \kappa\}$ 
we have $\delta \in S \Rightarrow j_\delta = j \wedge B_\delta = B
\wedge \ell g(\bar c_\delta) = m_* \wedge \varphi_\delta(\bar x,\bar
b_\delta) = \varphi(\bar x,\bar b) \wedge k_\delta = k_* \wedge
\bigwedge\limits_{\ell < k_*} \varepsilon_\ell(\delta) = 
\varepsilon_\ell$.  Also \wilog \, by \ref{3k.0.7}(4) we have $\delta \in S
\Rightarrow \text{ tp}(\bar c_\delta,M_{\text{min}(S)}) = 
\text{ tp}(\bar c_{\text{min}(S)},M_{\text{min}(S)})$ recalling that
the number of such types is $\le 2^{|B|+|T|}$.
Choose $\delta(0) < \delta(1)$ from $S$ so both has cofinality $\ge
\kappa$ and $B_{\delta(0)} = B_{\delta(1)}$, tp$(\bar
c_{\delta(0)},M_{\delta(0)}) \subseteq \text{ tp}(\bar
c_{\delta(1)},M_{\delta(1)})$ by \ref{3k.0.7}(1) and $\bar b_{\delta(1)} =
\bar b_{\delta(0)}$ and $\varepsilon_i(\delta(0)) =
\varepsilon_i(\delta(1))$ for $i < k_*$.  But we could have chosen in
stage $\delta(1),\varepsilon_k$ for $k < k_*$ and
$k'_{\delta(1)} = k_* +1$ and $\varepsilon'_{k_*}(\delta(1)) 
= \delta(0)$, contradiction to the
maximality of $k$ in Subclause $(\beta)$ of $\boxtimes(g)$.
\bigskip

\noindent
\underline{Case 2}:  Not Case 1.  

Then for a club of 
$i < \lambda^+$ if cf$(i) \ge \kappa$ \then \, subclause $(\alpha)$ of
Clause (g) fails for $i$ hence $M_i$ exemplifies 
that we have gotten the desired conclusion in (A) of \ref{3m.4}.
\bigskip

\noindent
\underline{Proof of (B)}:  By induction on $i < \lambda^+$ we choose 
$M_i,M^+_i,B_i,j_i,\bar c_i$ such that
\mn
\begin{enumerate}
\item[$\boxtimes'$]   Clauses (a),(c),(d) of $\boxtimes$ and
\begin{enumerate}
\item[${{}}$]   $(h) \quad \langle M^+_j:j \le i\rangle$ is
$\prec$-increasing continuous and $A^+ \subseteq M^+_0$, 
\sn
\item[${{}}$]   $(i) \quad M_i \prec M^+_i$ and $M^+_i$ has cardinality
$\lambda$ and if $i$ is non-limit then 

\hskip25pt $M_i,M^+_i$ are $\kappa$-saturated,
\sn
\item[${{}}$]  $(j) \quad$ if cf$(i) \ge \kappa$ and there are $j_i
< i,m < \omega,B \prec M_{j_i}$ of cardinality 

\hskip25pt  $< \kappa$ and 
$p \in \bold S^m(M_i)$ which does not split over $B$ and has 

\hskip25pt $\ge 2$ extensions in $\bold S^m(M^+_j)$ 
\then \, subclause $(\beta)$ of clause $\boxtimes(g)$ 

\hskip25pt above holds (with $\bar b \in {}^{\ell g(\bar y)}(M^+_j))$.
\end{enumerate}
\end{enumerate}
\mn
The rest of the proof is similar to that of (A).
Alternatively, we choose $M_i,M^+_i$ by induction on $i \le \kappa$
such that Clauses (a),(c),(d) of $\boxtimes$ and (h),(i) of
$\boxtimes'$ and
\mn
\begin{enumerate}
\item[${{}}$]  $(k) \quad$ if $i = j+1$ then $(M^+_j,M_{i+1})$ is
$\kappa$-reduced.
\end{enumerate}
\mn
Then $(M^+_\kappa,M_\kappa)$ are as required on $(N,M)$.
\end{PROOF}
\bigskip

\noindent
\centerline {$* \qquad * \qquad *$}
\bigskip

\noindent
For the rest of the section we shall assume (as we use it all the time).
\begin{hypothesis}
\label{3n.0} 
$T$ is dependent.
\end{hypothesis}

\begin{definition}
\label{3n.1}  
1) $\bold S^{\text{nsp}}_{< \kappa}(A) =\{p \in \bold S(A):p$ does 
not split over some $B \subseteq A$ of cardinality $< \kappa\}$.

\noindent
2) $\bold S^{\text{nsp},\alpha}_{< \kappa}(A) \subseteq \bold S^\alpha(A)$ and 
$\bold S^{\text{nsp},< \alpha}_{< \kappa}(A) \subseteq 
\bigcup\limits_{\beta < \alpha} \bold S^\beta(A)$ are defined similarly.

\noindent
3) $\bold S^{\text{nsp}}_{\ge \kappa}(A) = \{p \in \bold S(A):p$ is
weakly orthogonal to $r$ for every $r \in \bold S^{\text{nsp}}_{<
\kappa}(A)\}$.

\noindent
4) $\bold S^{\text{nsp},\alpha}_{\ge \kappa}(A),
\bold S^{\text{nsp},< \alpha}_{\ge \kappa}(A)$ are defined similarly.
\end{definition}

\noindent
We may note
\begin{observation}
\label{3n.1d}
1) If tp$(\bar c_1,A)$ belongs to $\bold S^{\text{nsp},\alpha}_{\ge
   \kappa}(A)$ and $\bar c_2 \in {}^\beta{\gC}$ and Rang$(\bar c_2)
   \subseteq \text{ Rang}(\bar c_1)$ \then \, $\tp(\bar c_2,A)$ belongs to
   $\bold S^{\text{nsp},\beta}_{\ge \kappa}$.

\noindent
2) tp$(\bar c,A) \in \bold S^{\text{nsp},\alpha}_{\ge \kappa}(A)$ \Iff
   \, tp$(\bar c \rest u,A) \in \bold S^{\text{nsp},m}_{\ge
   \kappa}(A)$ for every finite $u \subseteq \ell g(\bar c)$.

\noindent
3) If tp$(\bar a,A) \in \bold S^m(A)$ is weakly orthogonal to tp$(\bar
   c,A)$ and does not split over $B \subseteq A$ and every $q \in
   \bold S^{< \omega}(B)$ is realized in $M$ \then \, $\tp(\bar a,A + \bar
   c)$ does not split over $B$.
\end{observation}

\begin{PROOF}{\ref{3n.1d}}
Straight.
\end{PROOF}

\begin{observation}
\label{3n.2}   
If $\kappa = \text{ cf}(\kappa) > |T|$, the model $M$ is
$\kappa$-saturated and $p \in {\bold S}^m(M)$, \then \, we can find $N,q$
such that: 
\mn
\begin{enumerate}
\item[$\circledast^1_N$]  $(a) \quad \|N\| = \|M\|^{< \kappa}$,
\sn 
\item[${{}}$]   $(b) \quad q \in \bold S^m(N)$ extends $p$,
\sn
\item[${{}}$]  $(c) \quad N$ is $\bold F^{\text{nsp}}_{<
\kappa}$-constructible over $M$,
\sn
\item[$\circledast^2_{N,q,\kappa}$]  $(a) \quad N$ is
$\kappa$-saturated and $q \in \bold S^m(N)$.
\sn
\item[${{}}$]  $(b) \quad$ if $r \in \bold S^{\text{nsp},< \omega}_{<
\kappa}(N)$ then $r,q$ are weakly orthogonal, i.e. 
$q \in \bold S^{\text{nsp},m}_{\ge \kappa}(N)$.
\end{enumerate}
\end{observation}

\begin{PROOF}{\ref{3n.2}}
Let $\bar c$ realize $p(\bar x)$ and let $C = \text{ Rang}
(\bar c)$, now we apply clause (A) of Theorem \ref{3m.4} with $M,M
\cup C,N$ here standing for $A,A^+,M$ there.
\end{PROOF}

\begin{theorem}
\label{3n.3}
\underline{The Tree-like Type Theorem}
Assume $q(\bar x) \in \bold S^{\text{\rm nsp},\alpha}_{\ge \kappa}(N)$
and $N$ is $\kappa$-saturated and $\kappa > \theta = |T| +|\alpha|$
and let $\bar z = \langle z_\alpha:\alpha < \theta\rangle$.
\underline{Then} we can find a sequence $\bar\psi = \langle
\psi_{\varphi(\bar x,\bar y)}(\bar x,\bar z):\varphi(\bar x,\bar y)
\in \bbL(\tau_T)\rangle$ of formulas such that for every $A \subseteq
M$ of cardinality $< \kappa$ there is $\bar c \in {}^\theta M$ such
that:
\mn
\begin{enumerate}
\item[$(a)$]   $\{\psi_{\varphi(\bar x,\bar y)}(\bar x,\bar c):\varphi(\bar
x,\bar y) \in \bbL(\tau_T)\} \subseteq q \rest \text{\rm Rang}(\bar c)
\subseteq q$,
\sn
\item[$(b)$]   for each $\varphi(\bar x,\bar y) \in \bbL(\tau_T)$ we
have $\psi_{\varphi(\bar x,\bar y)}(\bar x,\bar c) \vdash
\{\varphi(\bar x,\bar b):\bar b \in {}^{\ell g(\bar y)}A$ and
$\varphi(\bar x,\bar b) \in q\}$.
\end{enumerate}
\end{theorem}

\begin{PROOF}{\ref{3n.3}}
This follows from Claims \ref{3n.4}, \ref{3n.5} below.
\end{PROOF}

\begin{claim}
\label{3n.4}  
1) Assume that $\circledast^2_{N,q,\kappa}$
from the Claim \ref{3n.2} holds, which means $N$ is
$\kappa$-saturated and $q \in \bold S^{\text{\rm nsp}}_{\ge \kappa}(N)$.
\underline{Then}
\mn
\begin{enumerate}
\item[$\circledast^3_{N,q,\kappa}$]  if $M \prec N$ has
cardinality $< \kappa$ and $\varphi(x,y)$ is a formula with
parameters from $N$, \then \, for some $\psi(x,\bar d) =
\psi_{\varphi(x,y),M}(x,\bar d_{\varphi(x,y),M})
\in q$ and $\eta \in {}^M 2$ we have $\psi(x,\bar d) \vdash
p^{M,\eta}_{\varphi(x,y)}$ where

$p^{M,\eta}_{\varphi(x,y)} = \{\varphi(x,b)^{\eta(b)}:b \in M\}$; so
it is included in $q$. 
\end{enumerate}
\mn
2) Part (1) works also for $q \in \bold S^m(N)$, i.e. $q \in \bold
S^{\text{\rm nsp},m}_{\ge \kappa}(N)$ and 
$\varphi = \varphi(\bar x,\bar y)$ where
$\ell g(\bar x) = m,\ell g(\bar y) < \omega$.
\end{claim}

\begin{PROOF}{\ref{3n.4}}
Fix $M$ such that $M \prec N$ of cardinality $< \kappa$.

\noindent
1) First note that
\mn
\begin{enumerate}
\item[$(*)_1$]   if $D$ is an ultrafilter on $M$ then $q(x)$ is
weakly orthogonal to $r_D = \text{ Av}(D,N)$.
\end{enumerate}
\mn
[Why?  As $r_D$ does not split over $M$, by \ref{3k.4}(1).]

Second, note that
\mn
\begin{enumerate}
\item[$(*)_2$]   the following type cannot be finitely satisfiable in $M$:
\end{enumerate}

\begin{equation*}
\begin{array}{clcr}
r^*(y) = &\{(\exists x_1)(\psi(x_1,\bar d) \wedge
\varphi(x_1,y)):\psi(x,\bar d) \in q\} \cup \\
  &\{(\exists x_2)(\psi(x_2,\bar d) \wedge \neg \varphi(x_2,y)):
\psi(x,\bar d) \in q\}.
\end{array}
\end{equation*}
\mn
[Why?  Otherwise for some ultrafilter $D$ on $M$ we have
$\vartheta(y,\bar d) \in r^*(y) \Rightarrow \vartheta(M,\bar d) \in D$.  Let $b
\in {\frak C}$ realize Av$(D,N)$ so as $q(\bar x)$ is closed under
conjunctions, $q(x) \cup \{\varphi(x,b)\}$ and $q(x) \cup \{\neg
\varphi(x,b)\}$ are finitely satisfiable in ${\frak C}$,
and we get a contradiction to $(*)_1$.]
\mn
\begin{enumerate}
\item[$(*)_3$]  there is $\psi(x,\bar d) \in q$ such that
$\{(\exists x_1)(\psi(x_1,\bar d) \wedge \varphi(x_1,y)),(\exists
x_2)(\psi(x_2,\bar d) \wedge \neg \varphi(x_2,y))\}$ is satisfied
by no $b \in M$.
\end{enumerate}
\mn
[Why?  By the monotonicity in $\psi(x,\bar d)$ and $q$ being closed
under conjunctions this follows from $(*)_2$.]
\mn
\begin{enumerate}
\item[$(*)_4$]   let $\psi_{\varphi(x,y),M}(x,\bar d_{\varphi(x,y),M}) =
\psi(x,\bar d)$, from $(*)_3$,
\sn
\item[$(*)_5$]   for every $b \in M$ we have $N \models ``(\forall x)
(\psi(x,\bar d) \rightarrow \varphi(x,b))"$ or $N \models ``(\forall x)
(\psi(x,\bar d) \rightarrow \neg \varphi(x,\bar b))"$.
\end{enumerate}
\mn
[Why?  By logic this follows by $(*)_3$.]
\mn
\begin{enumerate}
\item[$(*)_6$]   there is $\eta \in {}^M 2$ such that for
every $b \in M$ we have $M \models ``(\exists x)(\psi(x,\bar d) \wedge 
\varphi(x,b))"$ iff $M \models ``\neg(\exists x)(\psi(x,y)) \wedge 
\neg \varphi(x,b)"$ iff $\eta(b)=1$.
\end{enumerate}
\mn
[Why?  By $(*)_5 + (*)_3$.]  

So we are done. 

\noindent
2) Similarly.
\end{PROOF}

\begin{claim}
\label{3n.5}  
1) In the previous claim \ref{3n.4}, fixing $p$, if
{\rm cf}$(\kappa) > |T|$ \then \, $\psi$ depends on $\varphi$ but 
does not depend on $M$ though $\bar d$ in general does, i.e. given
 $\varphi(\bar x,\bar y)$ we may assume without loss of generality that
$\psi_\varphi(\bar x,\bar d) = \psi_{\varphi(\bar x,\bar y)}
(x,\bar d_{\varphi(\bar x,\bar y),M})$.

\noindent
2) Assume $\circledast^2_{N,q,\kappa}$ from claim \ref{3n.2}, i.e. $N$
is $\kappa$-saturated and $q \in \bold S^{\text{\rm nsp},m}_{\ge \kappa}(N)$.
\underline{Then} the following partial order is $\kappa$-directed
\mn
\begin{enumerate}
\item[$(a)$]   elements: 
$q \restriction B$ for $B \subseteq N$ of 
cardinality\footnote{from some form of strongly dependent we should be able to
get ``essentially finite"} $\le|T|$ 
\sn
\item[$(b)$]   order:  $p_1 \le p_2$ if $p_2 \vdash p_1$.
\end{enumerate}
\end{claim}

\begin{PROOF}{\ref{3n.5}}
1) As if $N_1 \prec N_2 \prec N$ and $\|N_2\| < \kappa$, then
$\psi_{\varphi(\bar x,\bar y),N_2}(x,\bar d_{\varphi(\bar x,\bar y),N_2})$
can serve as $\psi_{\varphi(\bar x,\bar y),N_1}(\bar x,\bar
d_{\varphi(\bar x,q),N_2})$.

\noindent
2) Easy. 
\end{PROOF}

As a conclusion we can now show that a key fact in \cite{Sh:877} for
the theory $T = \text{ Th}(\bbQ,<)$ has a parallel for every dependent $T$.
\begin{conclusion}
\label{3n.6}
\underline{The Saturated Expansion Conclusion}
Assume
\mn
\begin{enumerate}
\item[$(a)$]   $N$ is $\kappa$-saturated,
\sn
\item[$(b)$]  $|A| < \kappa$,
\sn
\item[$(c)$]   if $\bar a \in A$ then $q_{\bar a} = 
\,\text{\rm tp}(\bar a,N) \in \bold S^{\nsp,<\omega}_{\ge \kappa}(N)$, 
see Definition \ref{3n.1}(3),
\sn
\item[$(d)$]  $N_{[A]}$ has elimination of quantifiers.
\end{enumerate}
\mn
\underline{Then} $N_{[A]}$ is $\kappa$-saturated.
\end{conclusion}

\begin{remark}
1) Recall $N_{[A]}$ is $N$ expanded by $R_{\varphi(\bar
x,\bar b)} = \{\bar a \in {}^{\ell g(\bar x)} N:{\frak C} \models
\varphi[\bar a,\bar b]\}$ for $\varphi(\bar x,\bar b)$ a formula with
parameters from $A$, see \cite[\S1]{Sh:783}.

\noindent
2) We can omit assumption (d) in \ref{3n.6}, but then get
$\kappa$-saturated only for quantifier free types.
\end{remark}

\begin{PROOF}{\ref{3n.6}}
Without loss of generality $\kappa$ is regular, this as it is
enough to prove $\lambda^+$-saturation for every $\lambda \in [|A|,\kappa)$.
Let $M \prec N$ be such that $\|M\| < \kappa$ and assume $p = p(\bar y) 
\in \bold S^m(|M|,N_{[A]})$ and we shall prove that some $\bar c \in {}^m
N$ realizes $p(\bar y)$.  Actually \wilog \, $M_{[A]} \prec N_{[A]}$ and
by assumption (d), equivalently $p(\bar y) \in \bold S^m(M \cup A) =
\bold S^m(M \cup A,{\frak C})$ is finitely satisfiable in $M$.  
Let $\bar{\bold c} = \langle c_\alpha:\alpha < \alpha^*\rangle$ list 
$A$ so $\alpha^* < \kappa$ and for $u \subseteq \alpha^*$ let $\bar
c_u = \langle c_\alpha:\alpha \in u\rangle,\bar x_u = \langle 
x_\alpha:\alpha \in u\rangle$.
  
Next note that by Claims \ref{3n.5}(1) and
\ref{3n.4}(2) (here clause (c) of the
assumption is used) applied to tp$(\bar c_u,N)$ noting $\bar y$ of
length $m$ is fixed and letting $\bar x_u = \langle x_\alpha:\alpha
\in u\rangle$, we have:
\mn
\begin{enumerate}
\item[$(*)_1$]   for every finite $u \subseteq \alpha^*$  
and formula $\varphi = \varphi(\bar x_u,\bar y,\bar z)
\in \bbL(\tau_T)$ there are 
$\psi_{\varphi(\bar x_u,\bar y,\bar z)}(\bar x_u,\bar d_{\varphi(\bar x_u,
\bar y,\bar z),M}) \in \text{ tp}(\langle c_\alpha:\alpha \in u\rangle,N)$ so 
$\bar d_{\varphi(\bar x_u,\bar y,\bar z),M} \in {}^{\omega >} N$ and 
$\eta$ a function from ${}^{\ell g(\bar y)+\ell g(\bar z)}M$ to $\{0,1\}$
such that:
\end{enumerate}

\[
\psi_{\varphi(\bar x_u,\bar y,\bar z)}(\bar x_u,\bar d_{\varphi(\bar x_u,\bar
y,\bar z),M}) \vdash 
\{\varphi(\bar x_u,\bar b,\bar c)^{\eta(\bar b \char 94 \bar
c)}:\bar b \in {}^{\ell g(\bar y)}M \text{ and } \bar c 
\in {}^{\ell g(\bar z)}M\}.
\]

\mn
Clearly $|p(\bar y)| < \kappa$ so there are $\zeta^* < \kappa$ and a sequence
$\langle(\varphi_\zeta(\bar x_{u_\zeta},\bar y,\bar z_\zeta),u_\zeta):
\zeta < \zeta^*\rangle$ listing the pairs $(\varphi(\bar x_u,\bar
y,\bar z),u)$ as above so we have

\[
p(\bar y) = \{\varphi_\zeta(\bar{\bold c} \restriction u_\zeta,\bar y,\bar e):
\zeta < \zeta^* < \kappa \text{ and } \bar e \in {}^{\ell g(\bar z_\zeta)}M\},
\]

\mn
so $u_\zeta \subseteq \alpha^*$ is finite.

For each $\zeta < \zeta^*$ we choose 
$\psi_\zeta(\bar x_{u_\zeta},\bar d_\zeta)$ as guaranteed by $(*)_1$
above (for $\varphi_\zeta(\bar x_{u_\zeta},\bar y,\bar z))$.

Let

\[
p'(\bar y) := \{(\forall \bar x_{u_\zeta})[\psi_\zeta(\bar x_{u_\zeta},
\bar d_\zeta) \rightarrow \varphi_\zeta(\bar x_{u_\zeta},\bar y,\bar
e)]:\zeta < \zeta^* \text{ and } \bar e \in {}^{\ell g(\bar z_\zeta)}M\}.
\]

\mn
Now
\mn
\begin{enumerate}
\item[$(*)_2$]   $p'(\bar y)$ is finitely satisfiable in $M$.
\end{enumerate}
\mn
[Why?  As $p(\bar y)$ finitely satisfiable in $M$, 
using the translation and the choice of
$\psi_\zeta(\bar x_{u_\zeta},\bar d_\zeta)$.  That is, let $p''(\bar y)$
be a finite subset of $p'(\bar y)$ so for some $k < \omega,\zeta_\ell
= \zeta(\ell) < \zeta^*,\bar c_\ell \in {}^{\ell g(\bar z_\zeta)} M$  
for $\ell <k$ we have $p''(\bar y) = 
\{(\forall \bar x_{u_{\zeta(\ell)}}[\psi_{\zeta(\ell)}
(\bar x_{u_{\zeta(\ell)}},\bar d_{\zeta(\ell)})\rightarrow 
\varphi_\zeta(\bar x_{u_{\zeta(\ell)}},\bar y,\bar e_\ell)]:
\ell < k\}$.  Now $\{\varphi_{\zeta_\ell}(\bar{\bold c}
\restriction u_{\zeta_\ell},\bar y,\bar e_\ell):\ell < k\}$ is a
finite subset of $p(\bar y)$ hence is realized by some $\bar b \in {}^m
M$, hence by $(*)_1$ the sequence $\bar b$ realizes $p''(\bar y)$.]
\mn
\begin{enumerate}
\item[$(*)_3$]   the type $p'(\bar y)$ is over $\cup\{\bar
d_\zeta:\zeta < \zeta^*\} \cup M \subseteq N$.
\end{enumerate}
\mn
[Why?  Check.]
\mn
\begin{enumerate}
\item[$(*)_4$]   $p'(\bar y)$ has cardinality $\le |A| + |T| + \|M\| <
\kappa$.
\end{enumerate}
\mn
[Why?  Obvious.]
\mn
\begin{enumerate}
\item[$(*)_5$]   there are $M^+$ and $\bar b$ such that:
\begin{enumerate}
\item[$(a)$]   $M \prec M^+ \prec N$,
\sn
\item[$(b)$]   tp$(M^+,M \cup \bigcup\{\bar d_\zeta:\zeta <
\zeta^*\})$ is finitely satisfiable in $M$,
\sn
\item[$(c)$]   $\bar b \in M^+$ realizing $p'(\bar y)$.
\end{enumerate}
\end{enumerate}
\mn
[Why?  Easy, e.g. using ultrapower, ``$N$ is $\kappa$-saturated" and
$(*)_2 + (*)_3$.]
\mn
\begin{enumerate}
\item[$(*)_6$]   $\bar b$ realizes $p(\bar y) \in \bold S^m(N \cup
A)$ and $\bar b \in {}^{\ell g(\bar y)}(M^+) \subseteq
{}^{\ell g(\bar y)} N$.
\end{enumerate}
\mn
[Why?  Follow the translations.]

\noindent
So we are done.  
\end{PROOF}

\begin{question}
\label{3n.7}  
1) Can we waive assumption (d) in \ref{3n.6}?

\noindent
2) Is the family of $(N,A)$ as in \ref{3n.6} ``dense under $\le_\kappa$"?
\end{question}

\begin{discussion}
\label{3n.8.17} 
1) Assume $\lambda = \lambda^{< \lambda} > \kappa =
\text{ cf}(\kappa) > |T|$ and we try to prove that there is a
$(\lambda,\kappa)$-limit model.

So let $M \in \text{ EC}_\lambda(T)$ be saturated and we try to analyze
the class of $N,M \prec N \in \text{ EC}_\lambda$, which are ``close
enough", in the sense of $(\lambda,\kappa)$-limit model.

So if $p \in \bold S^{\text{nsp}}_{< \lambda}(M)$, say $p$ does not
split over $B$, for some $B \subseteq M$ of cardinality $< \lambda$, 
then we can assume that in $N$ there are ``enough
elements" realizing ``types not-splitting over $B$" extensions of $p$.  So
hopefully we can analyze such $N$ by ${\cP} \subseteq \bold
S^{\text{nsp}}_{\ge \lambda}(M)$ pairwise perpendicular or ${\cP}
\subseteq \bold S^{\text{nsp},|T|}_{\ge \lambda}(M)$ such that for each $p \in
{\cP}$ the model $M_{[p]}$ from \ref{0n.4}(1A) 
has elimination of quantifiers and is saturated, it
is reasonable that this holds if we can expand $M$ by definition of $<
\lambda$ types $p \in {\cP}$.

What we need, i.e., what is necessary for this line of attack (but not
yet clear if sufficient to carry it), is:
\mn
\begin{enumerate}
\item[$(*)$]   if ${\cP}_\ell$ above has cardinality
$\lambda$ and is quite dense (e.g. using $\bold F$'s for $\ell=1,2$)
\then \, there is an automorphism of $M$ which maps ${\cP}_1$ onto
${\cP}_2$.
\end{enumerate}
\mn
This leads to the generic pair conjecture.

\noindent
About $\bold S^{\text{nsp}}_{< \lambda}(M)$ recall Definition
\ref{3n.1}(1).
\end{discussion}
\newpage

\section {The type decomposition theorem} \label{s:the}

\begin{context}
\label{tp.45.8}  
1) $T$ is a complete first order theory; dependent if not said otherwise.

\noindent
2) ${\frak C} = {\frak C}_T$ is a monster for $T$, etc. as in \ref{0n.1}.
\end{context}

Here we try to analyze a type $p \in \bold S^{\le \theta}(M)$ for
$\kappa$-saturated $M \prec {\frak C}$ where $\kappa > \theta \ge
|T|$, the characteristic case being $\kappa >> \theta$ ($\theta$ may be
$\aleph_0$, if $T$ is countable).
In the case of $\theta < |T|$, or even better $\theta < \aleph_0$ 
we know less but mention it.  We look at ``$T$ being stable" as
our dream, our paradise.  The hard reality is $T$ being just dependent.  In
some sense $T$ dependent should be like stable but we allow order,
e.g. Th$(\Bbb Q,<)$ or trees.  What we actually do is investigate the
$K_\ell$ (see Definition \ref{tp14.21}).

How helpful is this analysis?  We present two consequences.  The
first to some extent accomplished the professed aim:
 the Type Decomposition Theorem \ref{tp25.43}.

What is its meaning?  If 
$M$ is $\kappa$-saturated, $\bar d \in {}^{\theta^+>}{\frak C}$ and
$\kappa > \theta \ge |T|$ then we try to analyze the type tp$(\bar
d,M)$ in two steps:  for some $\bar c,B$:
\mn
\begin{enumerate}
\item[$(a)$]  $B \subseteq M$ has cardinality $< \kappa$, say $B=|N|$,
\sn
\item[$(b)$]   the type in the first step is similar to the types 
of stable theories,
i.e. tp$(\bar c,M)$ does not split over $B \subseteq M$; (we
can even demand tp$(\bar c,M)$ is finitely satisfiable in $B$),
\sn
\item[$(c)$]   the type in the second step, 
tp$(\bar d,M + \bar c)$ behaves as in trees;
e.g. letting $\bold x = (M,B,\bar{\bc},\bd)$ we have: on ${}^\theta M$
the partial orders $\le_{\bold x}$ is $\kappa$-directed (see
\ref{0n.22}) where we let
$\bar a_1 \le_{\bold x} \bar a_2$ iff tp$(\bar d,\bar c \char 94 \bar
a_2) \vdash \text{ tp}(\bar d,\bar c \char 94 \bar a_1)$.
\end{enumerate}
\mn
The reader may say that Clause $(b)$ is not a true parallel to a stable case,
as $|B|$ is not bounded by $\theta +|T|$ (but this is impossible even
for the theory of dense linear order).  Still a type not splitting
over a set is a weak form of definability.  Also we may wonder, what
is the meaning, when
$T$ is Th$(\bbQ,<)$?  If $M$ is $\kappa$-saturated each $p \in \bold
S(M)$ actually stands for a cut of $M$.  Now the cuts are divided to
those which have cofinality
$\ge \kappa$ from both sides (falling under (c)), and those which do
not (hence fall under (b)).

The second consequence deals with singular $\mu$ of cofinality $>
|T|$.  We ask: is there $M \prec {\frak C}$ which is exactly 
$\mu$-saturated, i.e. is $\mu$-saturated but not
$\mu^+$-saturated.  Now if $T = \Th(\bbQ,<)$ this is
impossible, that is, there is no such $M$.  If $T$ is stable there 
is no problem to find such $M$, the main case being $\cf(\mu) > |T|$
(or just $\cf(\mu) \ge \kappa(T)$, see \cite[Ch.III]{Sh:c}) 
and let $M$ be $\mu$-prime over an
indiscernible set of cardinality $\mu$.  The result says that for
 dependent $T$ there is something like that, this is \ref{tp16.14}.

\begin{lemma}
\label{tp16.14}
\underline{Singular Exact Saturation Lemma}  
Assume $\kappa$ is singular of cofinality $> |T|$ and $M \prec {\frak C}$ is
an exactly 
$\kappa$-saturated model, i.e. is $\kappa$-saturated but not
$\kappa^+$-saturated.  \Then \, there are $N$ and $A$ such that:
\mn
\begin{enumerate}
\item[$(a)$]    $N \prec M$ of cardinality $<
\kappa$ and $A \subseteq M$ of cardinality $\kappa$ and $M$ 
omits some $p \in \bold S(A)$ which does not split over $N$; in fact
\sn
\item[$(b)$]   there is $q \in \bold S(M)$ which does not split over $N$
such that $p = q \restriction A$,
\sn
\item[$(c)$]   there is an indiscernible sequence $\langle
\bar a_\alpha:\alpha < \kappa\rangle$ over $N$ of $\theta$-tuples from
$M$ such that 
{\rm Av}$(\langle \bar a_\alpha:\alpha < \kappa\rangle,
N \cup \{\bar a_\alpha:\alpha < \kappa\})$ is
omitted by $M$ (equivalently, we cannot choose $\bar a_\kappa \in M$) and
does not split over $N$.
\end{enumerate}
\end{lemma}

\begin{remark}
We can add in (c): 
\mn
\begin{enumerate}
\item[$(c)^+$]  moreover we can demand that there is an 
ultrafilter ${\cD}$ on $N$ such  that $M$ omits $p = 
\text{\rm Av}({\cD},A)$ where $A := \cup\{\bar a_\alpha:
\alpha < \kappa\} \cup N$ and $p(\bar x)$ is as in Clause (c).
\end{enumerate}
\end{remark}

\begin{theorem}
\label{tp25.43}
\underline{The Type Decomposition Theorem}  
Assume {\rm cf}$(\kappa) > \theta \ge |T|,M$ is
$\kappa$-saturated and $\bar{\bd} \in {}^{\theta \ge}{\frak C}$.
\underline{Then} for some $\bar{\bold c} \in {}^\theta{\frak C}$, 
recalling Definition \ref{3n.1} we have {\rm tp}$(\bar{\bc},M) 
\in \bold S^{\text{\rm nsp},\theta}_{< \kappa}(M)$ and 
$(\bold P,\le_{\bold P})$ is a $\kappa$-directed partial order 
\underline{where}:
\mn
\begin{enumerate}
\item[$(a)$]   $\bold P = \{\text{\rm tp}(\bar{\bd},A \cup
\bar{\bc}):A \subseteq M$ has cardinality $\le \theta\}$,
\sn 
\item[$(b)$]   $p_1 \le_{\bold P} p_2$ iff $p_2 \vdash p_1$.
\end{enumerate}
\end{theorem}

\begin{remark}  
1) In fact $(\bold P,\le_{\bold P})$ is 
$(\bold P_{\bold x,\theta},\le_{\bold x,\theta})$ from the Definition
\ref{tp14.21}(8) below.

\noindent
2) Note that being $\theta^+$-directed is obvious.

\noindent
3) Would it be more transparent to use the notation $p_2 \vdash p_1$ instead of
   $p_1 \le_{\bold P} p_2$?  A matter of taste, the author feels that not.
\end{remark}

\begin{definition}
\label{tp14.21}  
1) Let $K = K_1$ be the family of $\bold x$ satisfying
\mn
\begin{enumerate}
\item[$(a)$]   $\bold x = (A,B,\bar{\bc},\bar{\bd})$ but if
$A=|M|$, as usual, we may write $M$ instead of $A$ and 
if $B = \emptyset$ we may omit it,
\sn
\item[$(b)$]   $B \subseteq A$,
\sn
\item[$(c)$]   $I$ a linear order,
\sn
\item[$(d)$]   $\bar{\bc} = \langle \bar c_{t,n}:n < n_t,t
\in I\rangle$ where $n_t \le \omega$, each $\bar c_{t,n}$ a 
finite\footnote{we remark when it matters.} sequence and 
let\footnote{but abusing our notation, $\Rang(\bar{\bold c})$ is the
  set of elements of $\gC$ appearing in it; similarly in other cases}
 $\bar c_t = \bar c_{t,0} \char 94 \bar c_{t,1} \char 94 \ldots
\char 94 \bar c_{t,n_t-1}$, 
\sn
\item[$(e)$]    $\langle \bar c_{t,n}:n < n_t\rangle$ is an
indiscernible sequence over $A \cup \{\bar c_s:s \in I \backslash
\{t\}\}$, so if $n_t=1$ this is an empty statement,
\sn
\item[$(f)$]    if $t \in I$ then\footnote{of course, this implies that
clause (e) follows from a weak version, see \ref{tp14.28}(5) but see part (2).}
tp$(\bar c_t,\{\bar c_s:s <_I t\} \cup A)$ does not split over $B$
\sn
\item[$(g)$]    $\bar{\bold d}$ is a sequence of elements or finite
  sequences from $\gC$.
\end{enumerate}
\mn
2) Let $K_0$ be defined similarly omitting clause $(f)$.

\noindent
3) For $\lambda \ge \kappa$, cf$(\lambda) \ge \theta$ (or just
   $\lambda \ge \theta$), cf$(\kappa) \ge \theta$ and $\ell \in \{0,1\}$ let 
$K^\ell_{\lambda,\kappa,<\theta} = \{(M,B,\bar{\bold c},\bar{\bold d}) 
\in K_\ell:M$ is $\lambda$-saturated, $|B| < \kappa$ and $|\ell g(\bar d)|
+ |\ell g(\bar{\bold c})|| < \theta\}$; omitting $\ell$ means 1.
   If $\theta = \sigma^+$ instead of ``$< \theta$" we may write $\sigma$.

\noindent
4)
\begin{enumerate}
\item[$(a)$]   let $\bold x = (A_{\bold x},B_{\bold x},
\bar{\bold c}_{\bold x},\bar{\bold d}_{\bold x})$ for 
$\bold x \in K_0$ (or $M_{\bold x}$ instead of $A_{\bold x})$,
$I = I_{\bold x},\bar c_{\bold x,t} = \bar c_t,n_t = 
n_{\bold x,t}$ and $C_{\bold x} = \cup\{\text{Rang}
(\bar c_{t,n}):t \in I$ and $n < n_{\bold x,t}\}$,
\sn
\item[$(b)$]   we may\footnote{we may write $C_{\bold x} = \text{
Rang}(\bar{\bold c}),\bar c_{\bold x,t,n} = \langle c^{\bold
x}_{t,n,m}:m < \ell g(c_{\bold x,t,n})\rangle$ so in fact
$\bar c_{\bold x} = \langle c^{\bold x}_{t,n,m}:(t,m,m) \in J\rangle$
for the natural $J = J_{\bold x}$.}
 replace $\bar{\bold d}$ by $D_{\bold x} = \text{ Rang}(\bar{\bold d})$,
\sn
\item[$(c)$]  we may omit $\kappa$ if $\kappa = \lambda$.
\end{enumerate}
\mn
For $\lambda = \aleph_0$, let $``< \lambda"$ mean $A_{\bold x}$ is the
universe of $M_{\bold x} \prec {\frak C}$ (no saturation demand).

\noindent
5) We define a two-place relation $\le_0$ on $K_0:\bold x \le_0 \bold y$ iff 
$A_{\bold x} \subseteq A_{\bold y},
B_{\bold x} \subseteq B_{\bold y},I_{\bold x} \subseteq I_{\bold y},
\bar{\bold c}_{\bold x} = \bar{\bold c}_{\bold y} \restriction
I_{\bold x}$, i.e. $t \in I_{\bold x} \Rightarrow \bar c_{\bold y,t} =
\bar c_{\bold x,t}$ moreover $t \in I_{\bold x} \Rightarrow n_{\bold
y,t} = n_{\bold x,t}$ and $t \in I_{\bold x} \wedge n < n_{\bold x,t}
\Rightarrow \bar c_{\bold y,t,n} = \bar c_{\bold x,t,n},
\bar{\bold d}_{\bold x} \trianglelefteq \bar{\bold d}_{\bold y}$
and tp$(\bar{\bold c}_{\bold x},A_{\bold y})$ does not split over
$B_{\bold x}$ hence ``tp$(\bar{\bold c}_{\bold x},A_{\bold x})$
does not split over $B_{\bold x}$" follows.

\noindent
5A) $\bold x \le_1 \bold y$ mean $\bold x \le_0 \bold y$ and
$A_{\bold y} = A_{\bold x}$.

\noindent
6) We define a two-place relation $\le_2$ on $K_0:\bold x \le_2 \bold
y$ iff $\bold x \le_1 \bold y \wedge \bar{\bold d}_{\bold x} = \bar{\bold
d}_{\bold y}$.

\noindent
7) $\bold x \in K^0_{\lambda,\kappa,\theta}$ is called normal \when \,
Rang$(\bar{\bold c}_{\bold x}) \subseteq \text{ Rang}(\bar{\bold
d}_{\bold x})$.  

\noindent
8) For $\bold x \in K$, let $\bold P_{\bold x,\theta} = 
(\bold P_{\bold x,\theta} \le_{\bold x,\theta})$ be defined by:
\mn
\begin{enumerate}
\item[$(a)$]   $\bold P_{\bold x} = \{\tp(\bar\bd_{\bold x},A +
\bar\bc_{\bold x}):A \subseteq A_{\bold x}$ has cardinality $\le \theta\}$,
\sn
\item[$(b)$]   $\le_{\bold x,\theta}$ is the following two-place
relation on $\bold P_{\bold x,\theta}:p_1(\bar x_{\bar\bd_{\bold x}})
\le_{\bold x,\theta} p_2(\bar x_{\bar\bd_{\bold x}})$ iff $p_2 \vdash p_1$.
\end{enumerate}
\mn
9) If $\theta = |T| + |\ell g(\bar{\bd}_{\bold x})| + |\ell
g(\bar{\bc}_{\bold x})|$, i.e. we use ``$< \theta^+$", we may omit it.
\end{definition}

\begin{claim}
\label{tp14.28}  
1) $K_1 \subseteq K_0$.

\noindent
2) $\le_i$ is a partial order on $K_0$ for $i=0,1,2$.

\noindent
3) If $i \in \{0,1,2\},\langle \bold x_\alpha:\alpha < \delta\rangle$ is
 $\le_i$-increasing in $K^0_{\lambda,\kappa,\theta}$ where $\delta$ is
 a limit ordinal, $[\cf(\delta)
 \ge \theta^+ \Rightarrow \bigwedge\limits_{j < \delta} 
\bar\bc_{\bold x_j} =  \bar\bc_{\bold x_0}]$ and\footnote{This follows
 from ``$\langle \bold x_\alpha:\alpha < \delta\rangle$ is
 $\le_i$-increasing" when $i=2$.}$[i \le 1 \wedge \text{\rm cf}(\delta)
 \ge \theta^+ \Rightarrow (\bigwedge\limits_{\alpha < \delta} 
\bar{\bd}_{\bold x_\alpha} = \bar{\bd}_{\bold x_0})],[i =
 0 \Rightarrow \lambda \le \text{\rm cf}(\delta)]$ and 
$\delta < \text{\rm cf}(\kappa) \vee (\bigwedge\limits_{\alpha <
 \delta} B_{\bold x_\alpha} = B_{\bold x_i})$, \then \, it
 has a $\le_i$-{\rm lub} $\bold x_\delta := \cup\{\bold x_\alpha:\alpha <
\delta\} \in K^0_{\lambda,\kappa,\theta}$ defined by
$A_{\bold x} = \cup\{A_{\bold x_\alpha}:\alpha < \delta\},
B_{\bold x} = \cup\{B_{\bold x_\alpha}:\alpha < \delta\},
I_{\bold x} = \cup\{I_{\bold x_\alpha}:\alpha < \delta\},
\bar{\bc}_{\bold x} = \cup\{\bar{\bc}_{\bold x_\alpha}:
\alpha < \delta\}$, i.e. $\bar c_{\bold x,t} = \bar c_{\bold
 x_\alpha,t}$ when $t \in I_{\bold x_\alpha}$ and $\bar{\bd}_{\bold x} = 
\cup\{\bar{\bd}_{{\bold x}_\alpha}:\alpha < \delta\}$.

\noindent
3A) In part (3), if $\alpha < \delta \Rightarrow \bold x_\alpha \in
K^1_{\lambda,\kappa,\theta}$ \then \, $\bold x_\delta \in
K^1_{\lambda,\kappa,\theta}$. 

\noindent
4) If $\bar{\bd} \in{}^{\theta^+ >}{\frak C}$ and $M$ is
$\kappa$-saturated and $\kappa > \theta$ \then \, $\bold x =
(M,\emptyset,<>,\bar{\bd})) \in K^\ell_{\kappa,\theta}$ for $\ell =0,1$. 

\noindent
5) In the definition of $\bold x \in K_1$: in clause (e) it 
suffices to demand that: if $n_t > 1$ \then \,
$\langle \bar c_{t,n}:n < n_t\rangle$ 
is indiscernible over $A \cup \{\bar c_{s,m}:s <_I t,m<n_s\}$.

\noindent
6) For every $\bold x \in K_{\lambda,\kappa,\theta}$ there is a normal 
$\bold y \in K_{\lambda,\kappa,\theta}$ satisfying $\bold x 
\le_1 \bold y,\bar{\bc}_{\bold x} = \bc_{\bold y}$ 
and {\rm Rang}$(\bar{\bd}_{\bold y}) = 
\text{\rm Rang}(\bar{\bd}_{\bold x}) \cup \text{\rm Rang}
(\bar{\bc}_{\bold x})$.  Hence $\bold y \in 
\text{\rm mxK}_{\lambda,\kappa,\theta}
\Leftrightarrow \bold x \in \text{\rm mxK}_{\lambda,\kappa,\theta}$, see
Definition \ref{tp25.32} below.

\noindent
7) If $i=0$ and $\langle \bold x_\alpha:\alpha < \delta\rangle$ is
$\le_i$-increasing in $K^0_{\lambda,\kappa,\theta}$ and {\rm cf}$(\delta) <
\theta^+,\delta < \text{\rm cf}(\kappa)$, \then \, the sequence has a
$\le_i$-upper bound $\bold x_\delta \in K^0_{\lambda,\kappa,\theta}$,
note that we have not said ``lub".

\noindent
7A) In part (7), if $\alpha < \delta \Rightarrow \bold x_\alpha \in
K^1_{\lambda,\kappa,\theta}$ \then \, we can add $\bold x_\delta \in
K^1_{\lambda,\kappa,\theta}$.
\end{claim}

\begin{PROOF}{\ref{tp14.28}}
Easy e.g.

\noindent
7), 7A)  The problem is when part (3) does not cover it, so $\lambda >
\aleph_0$.  It is clear how to choose $\bar{\bc}_{\bold x_\delta},
\bar{\bd}_{\bold x_\delta}$ 
and $B_{\bold x_\delta}$, but we should choose a
$\lambda$-saturated $M_{\bold x_\delta}$.

Let $B = \cup\{B_{\bold x_\alpha}:\alpha < \delta\},I = \cup\{I_{\bold
x_\delta}:\alpha < \delta\}$ and $\bar{\bc} = \langle \bar c_t:t \in
I\rangle$ with $\bar c_t = \bar c_{\bold x_\alpha,t}$ when $\alpha \in
I_{\bold x_\alpha}$; similarly $\bar{\bold d}$. 

First, choose a $\lambda$-saturated $M$ extending $\cup\{M_{\bold
x_\alpha}:\alpha < \delta\}$ but what about ``{\rm tp}$(\bar c_{\bold
x_\alpha},M)$  does not split over $B_{\bold x_\alpha}$ for each
$\alpha < \delta$"?

Now for each $\alpha < \gamma < \delta$, tp$(\bar{\bc}_{\bold
x_\alpha},M_\gamma)$ does not split over $B_{\bold x_\alpha}$ which means
$p_{\alpha,\gamma}(\bar x) = \text{ tp}(\bar{\bc}_{\bold
x_\alpha},M_\gamma)$ does not split over $B_{\bold x_\alpha}$ hence
$p_\alpha(\bar x) := \cup \{p_{\alpha,\beta}(\bar x):
\beta \in (\alpha,\delta)\} = \tp(\bar{\bc}_{\bold x_\alpha},
\bigcup\limits_{\beta < \delta} M_\beta)$ does not 
split over $B_{\bold x_\alpha}$.  Also for
$\alpha < \delta$ by \ref{3k.0.7}(1) there is 
$\bar{\bc}'_\alpha$ such that tp$(\bar{\bc}'_\alpha,M)$ does not
split over $B_{\bold x_\alpha}$ and extends $p_\alpha$.  As $M_{\bold
x_\alpha} \supseteq B_{\bold x_\alpha}$ is $\lambda$-saturated and
$\lambda \ge \kappa$ by Definition \ref{tp14.21}(3) clearly the model
$M_{\bold x_\beta}$ is $\kappa$-saturated and 
$|B_{\bold x_\alpha}| < \kappa$ by the
definition of $K^i_{\lambda,\kappa,\theta}$.  
Recalling \ref{3k.0.7}(1), by the last
two sentences $\bar{\bc}'_\beta \restriction I_{\bold x_\alpha}$ 
realizes tp$(\bar{\bc}'_\alpha,M)$ for
$\alpha < \beta < \delta$ hence \wilog \, $\alpha < \beta < \delta
\Rightarrow \bar{\bc}'_\beta \restriction I_{\bold x_\alpha} =
\bar{\bc}'_\alpha$.

Hence there is an elementary mapping $f$ mapping with domain
$\cup\{\bar{\bc}'_\alpha:\alpha < \delta\} \cup M$, mapping 
$\bar{\bc}'_\alpha$ to $\bar{\bc}_{\bold x_\alpha}$ 
for $\alpha < \delta$, and
extending {\rm id}$_{\cup\{M_\alpha:\alpha < \delta\}}$.  Now $M_{\bold
x_\delta} := f(M)$ will do, i.e. let $M_{\bold y} = M,B_{\bold y} =
B,\bar{\bc}_{\bold y} = \bar{\bc},\bar{\bd}_{\bold y} = \bar{\bd}$.
\end{PROOF}

\begin{definition}
\label{tp25.32}  
1) For $\ell = 0,1$ let 
mxK$^\ell_{\lambda,\kappa,<\theta}$ be the family of $\bold x \in
K^\ell_{\lambda,\kappa,<\theta}$ which are $\le^+_2$-maximal in
$K^\ell_{\lambda,\kappa,< \theta}$, i.e. for no $\bold y$ do we have
$\bold x <^+_2 \bold y \in K^\ell_{\lambda,\kappa,< \theta}$, see
below; if $\ell=1$ we may omit it.

\noindent
2) For $i = 0,1,2$ let $\le^+_i$ be the following two-place relation on
$K_0:\bold x \le^+_i \bold y$ \Iff \, $\bold x \le_i \bold y$, see
 Definition \ref{tp14.21} and if $\bold x \ne \bold y$ \then \,
for some $t \in I_{\bold y} \backslash I_{\bold x}$ satisfying
$n_{\bold x,t} \ge 2$ we have: ${\frak C} \models 
\varphi[\bar{\bd}_{\bold x},\bar c_{t,1},\bar b] \wedge \neg
\varphi[\bar{\bd}_{\bold x},\bar c_{t,0},\bar b]$ for some $\varphi =
\varphi_t(\bar x,\bar y,\bar z) \in \Bbb L(\tau_T)$ and $\bar b
\subseteq A_{\bold x} \cup \bigcup\{\bar c_{\bold y,s}:s \in I_{\bold y} 
\backslash \{t\}\}$. 

\noindent
3) Again, if $\theta = \sigma^+$ instead of ``$< \theta$": we may write
$\sigma$, and if $\kappa = \lambda$ we may omit $\lambda$.

\noindent
4) Of course, $\bold x <^+_i \bold y$ means $\bold x \le^+_i
\bold y \wedge \bold x \ne \bold y$.
\end{definition}

\begin{observation}
\label{tp25.31}  
Let $i=0,1,2$.  

\noindent
1) For $\ell=0,1$ the two-place relation $\le^+_i$ is
   a partial order on $K^\ell$.

\noindent
2) If $\bold x_1 \le_i \bold x_2 <^+_i \bold x_3 \le_i \bold x_4$ then
   $\bold x_1 <^+_i \bold x_4$.

\noindent
3) If $\bold x <^+_i \bold z$ are from $K^\ell_{\lambda,\kappa,<
\theta}$ \then \, there is $\bold y \in K^\ell_{\lambda,\kappa,<
   \theta}$ such that $\bold x <^+_2 \bold y \le_i \bold z$ and 
$\bar{\bd}_{\bold y} = \bar{\bd}_{\bold x},
I_{\bold y} \backslash I_{\bold x}$ is finite.

\noindent
4) The parallel of \ref{tp14.28}(3),(3A) holds for $\mxK$.
\end{observation}

\begin{PROOF}{\ref{tp25.31}}
E.g.

\noindent
3) Let $t \in I_{\bold z} \backslash I_{\bold x}$ and $\varphi =
   \varphi(\bar x,\bar y,\bar z) \in \bbL(\tau_T)$ and $\bar b
   \subseteq A_{\bold x} \cup \bigcup\{\bar c_{\bold y,s}:s \in
   I_{\bold y} \backslash \{t\}\}$ be such that $\gC
\models \varphi[\bar{\bold d}_{\bold x},\bar c_{\bold y,t,1},\bar b]
\wedge \neg \varphi[\bar{\bold d}_{\bold x},\bar c_{\bold y,t,0},\bar b]$.
   We choose a finite $I \subseteq I_{\bold y} \backslash \{t\}$
   such that $\bar b \subseteq \cup\{\bar c_{\bold y,s}:s \in I\} \cup A_{\bold
   x}$.  Now define $\bold y$ by: $M_{\bold y} = M_{\bold x},I_{\bold
   y} = I_{\bold x} \cup I \cup \{t\},\bar{\bc}_{\bold y} = \bar c_{\bold z}
   \rest I_{\bold y},\bar d_{\bold y} = \bar{\bd}_{\bold x}$ and
   $B_{\bold y} = B_{\bold x}$.

Now check.
\end{PROOF}

The following claim may be good for digesting the meaning of
mxK$^\ell_{\lambda,\kappa,\theta}$. 
\begin{claim}
\label{tp25.30}   
\underline{The L.S.T. Claim for {\rm mxK}}

If $\bold x \in \text{\rm mxK}_{\kappa,\theta}$ and 
$M = M_{\bold x}$ \then \, for 
some function $F$ with domain $[M]^{< \kappa}$ satisfying 
$F(A) \in[M]^{\le 2^{|A|+|T|}}$ for $A \in
\text{\rm Dom}(F)$, we have: if $M_1 \prec M$ is closed under $F$ and
contains $B_{\bold x}$ \then \,
$(M_1,B_{\bold x},\bar{\bc}_{\bold x},\bar{\bd}_{\bold x})$
belongs to {\rm mxK}$^\ell_{\kappa,\theta}$.
\end{claim}

\begin{remark}
By \ref{tp25.31}(4), it suffices to consider $F$ with domain $[M]^{<
  \aleph_0}$. 
\end{remark}

\begin{PROOF}{\ref{tp25.30}}  
We can choose $F(\emptyset) \in [M]^1$ and for notational transparency we
fix a set $J$ of
cardinality $\aleph_0$ disjoint to $I_{\bold x}$.

Note that for every $N \prec M$ satisfying $B_{\bold x} \subseteq
N$ we have $\bold x_N := (N,B_{\bold x},\bar\bc_{\bold x},
\bar\bd_{\bold x}) \in K_{\kappa,\theta}$; call such $N$ a candidate.  
So to choose $F$ let us
analyze the cases $B_{\bold x} \subseteq N \prec M$ but $\bold x_N
\notin \text{ mxK}^\ell_{\kappa,\theta}$.  Considering Definition
\ref{tp25.32}, it suffices by \ref{tp25.31}(3)
to consider the case $\bold x_N <^+_2 \bold y,I_{\bold y} 
\backslash I_{\bold x}$ is finite and is $\subseteq J,t_*
\in I_{\bold y} \backslash I_{\bold x},\varphi_*(\bar x,\bar y,\bar
z),\bar b_*$ as there, we can ignore the possibility that also
some other $t \in I_{\bold y} \backslash I_{\bold x} \backslash
\{t_*\}$ works.

We let $\bar b_0$ list $B_{\bold y},I_{\bold y} \backslash I_{\bold x} =
\{t_0,\dotsc,t_{k-1}\},t_* = t_{\ell(*)},\ell(*) < k,m_t > \ell g(\bar
c_{y,t,0}),n_t = n_{\bold y,t}$ for $t \in I_{\bold y}$. 
Also \wilog \, $\varphi = \varphi(\bar x,\bar y;\bar z_1,\bar z_2),\ell
g(\bar y) = \ell g(\bar c_{\bold y,t_*,0}),\ell g(\bar x) = \ell
g(\bar{\bd}_{\bold x}),\bar b = \bar b_1 \char 94 \bar b_2,\bar b_1 \in
{}^{\omega >}(M_{\bold y}) = {}^{\omega >}N,\ell g(\bar z_1) = \ell
g(\bar b_1)$ and abusing our notation, $\ell g(\bar z_2) = \ell
g(\bar b_2)$ where $\bar b_2 = \langle 
\bar c_{\bold y,t,k}:k < n_{\bold y,t},t \in I_{\bold y}\rangle$ and
$\bar b_{2,<t} = \langle c_{\bold y,s,k}:k < n_t,s < t$ so $t \in
I_{\bold y} \rangle$ and $\bar b_{2,\ne t} = \langle \bar c_{\bold y,s,k}:
k < n_s,s \in I_{\bold y} \backslash \{t\}\rangle$, so 
$\gC \models \varphi[\bar{\bd}_{\bold x},\bar
c_{\bold y,t_*,1},\bar b_1,\bar b_{2,t_*}] \wedge \neg 
\varphi[\bar{\bd}_{\bold x},\bar c_{\bold y,t_*,0},\bar b_1,
\bar b_{2,t_*}]$ and let 
$\ell g(\bar z^*_t) = \ell g(\bar c_{\bold y,t}),\ell g(\bar
z_{t,\ell}) = \ell g(\bar c_{\bold y,t,\ell})$ so $\bar z^*_t 
= \bar z_{t,0} \char 94 \ldots \char 94 \bar z_{t,n_t-1}$ and
$\bar z_2 = (\ldots,\bar z^*_t,\ldots)_{t \in I_{\bold y}}$.

All this information will be called a witness against the candidate
$N$ and we denote it by $\bold w$.

Let $\bold s$ consist of the following pieces of information on the
witness $\bold w$ and in this case 
we shall say that $\bold w$ materializes $\bold s$ 
and $\bold s$ is a case for $N$.
\mn
\begin{enumerate}
\item[$\boxtimes$]  $(a) \quad I = I_{\bold y}$ and $\langle
t_\ell:\ell < k\rangle$ (so we will write $I_{\bold s},t_{\bold
  s,\ell}$) and $\ell(*)$,
\sn
\item[${{}}$]  $(b) \quad n_t,m_t$ for $t \in I$,
\sn
\item[${{}}$]  $(c) \quad \varphi = \varphi(\bar x,\bar y,\bar
z_1,\bar z_2)$ hence $\ell g(\bar b_1) = \ell g(\bar z_1)$,
\sn
\item[${{}}$]  $(d) \quad \zeta(0) = \ell g(\bar b_0)$ but not $\bar
  b_0$ itself,
\sn
\item[${{}}$]  $(e) \quad q_0 = \text{ tp}(\bar b_0,\emptyset)$ and
$q_1 = \tp(\bar b_0 \char 94 \bar b_1 \char 94 \bar b_2,\emptyset)$,
  so from $q_1$ we know when 

\hskip25pt $b_{2,\ell_1} = c_{t_\ell,n,\ell_2}$
\sn
\item[${{}}$]  $(f) \quad$ the scheme of non-splitting of tp$(\bar c_t,M)$
for $t \in I_{\bold y}$ from 

\hskip25pt clause $(f)$ of \ref{tp14.21}(1), that is 
$\Xi_t = \{(\psi(\bar z_{\bar c_{\bold y,t}},\bar y'),
q(\bar y',\bar y_{\bar b_0}))$: for some 

\hskip25pt $\bar b \in {}^{\ell g(\bar y')}N$ we have
$\gC \models \psi[\bar c_{\bold y,t},\bar b]$ and $q(\bar y',\bar y_{\bar
b_0}) = \text{ tp}(\bar b \char 94 \bar b_0,\emptyset)\}$.
\end{enumerate}
\mn
We shall write $I = I_{\bold s},t_* = t_*(\bold s),\varphi 
= \varphi_{\bold s}(\bar x_{\bold s},\bar{\bold y},
\bar z_{\bold s,1},\bar z_{\bold s,2}),q_{\bold s,0} = q_0$, etc. 
and let $r_{\bold s} = \tp(\bar\bc_{\bold y},N)$.  
We call $\bold s$ a case when it is
a case for some candidate $N$.  If $\bold s$ is a case and $\bar b_0
\in {}^{\zeta(0)}(M_{\bold x})$ realizes $q_{\bold s,0}$, 
\then \, we can choose
$\bar{\bold c}_{\bold s,\bar b_0} = \langle \bar c_{\bold x,\bar
b_0,t}:t \in I_{\bold s}\rangle$ such that tp$(\bar c_{\bold s,\bar
b_0,t},\cup\{\bar c_{\bold s,\bar b_0,s}:s <_{I_{\bold s}} t\} \cup
M_{\bold x})$ is defined by the scheme $\Xi_t$ with the parameter
$\bar b_0$, this type is determined by $\bold s,\bar b_0$ and $\bold x$ (though
not the $\bar c_{\bold s,\bar b_0,t}$'s themselves).  Without loss of
generality $t \in I_{\bold x} \Rightarrow \bar c_{\bold s,\bar b_0,t}
= \bar c_{\bold x,t}$. 

Now clearly
\mn
\begin{enumerate}
\item[$(*)_1$]  for a candidate $N$ we have $\bold x_N \in
K_{\kappa,\theta}$ \Iff \, for every case $\bold s$ and $\bar
b_0 \in {}^{\zeta(0)}N$ realizing $q_{\bold s,0}$ there is no $\bar b_1 
\in {}^{\ell g(\bar y)}N$ such that $\bold s_1,\bar b_0,\bar b_1,
\bar b_2 = \bar c_{\bold s,\bar b_0}$ are as above.
\end{enumerate}
\mn
So
\mn
\begin{enumerate}
\item[$(*)_2$]  for every case $\bold s$ and $\bar b_0 \in
{}^{\zeta(0)}(M_{\bold x})$ realizing $q_{\bold s,0}$ and $\bar b_1 
\in {}^{\ell g(\bar z_{\bold s,1})}M$, 
we cannot choose $\bar{\bc}',\bar{\bc}''$
realizing tp$(\bar{\bc}_{\bold s,\bar b_0},M_{\bold x})$ such that
$\bar{\bc}' \rest I_{\bold x} = \bar{\bc}_{\bold x} = \bar{\bc}''
\rest I_{\bold x}$ and $\gC \models \varphi_{\bold s}[\bd_{\bold x},
\bar c'_{t_*(\bold s)},\bar b_1,\bar{\bc}' \rest (I_{\bold s} \backslash
\{t_*(\bold s)\})] \wedge \neg \varphi_{\bold s}
[\bd_{\bold x},\bar c''_{t_*(\bold s)},\bar b_1,\bar{\bc}' \rest
(I_{\bold x} \backslash \{t_*(\bold s)\})]$.
\end{enumerate}
\mn
[Why?  As then $\bold x \notin \text{ mxK}_{\kappa,\theta}$.]

Hence 
\mn
\begin{enumerate}
\item[$(*)_3$]  for every case $\bold s$ and $\bar b_0 \in
{}^{\zeta(0)}(M_{\bold x})$ realizing $q_{\bold s,0}$ and $\bar b_1 
\in {}^{\ell g(\bar z_{\bold s,1})}(M_{\bold x})$ there is a finite set
$C = C_{\bold s}(\bar b_0,\bar b_1) \subseteq M_{\bold x}$ such that: if
$N$ is a candidate which includes $\bar b_0,\bar b_1,C$ \then \, there
is no witness $\bold w$ against $N$ with $\bold s = \bold s_{\bold
w},\bar b_0 = \bar b_{\bold w,0},\bar b_1 = \bar b_{\bold w,1}$.
\end{enumerate}
\mn
Also (\ref{3k.0.7}(4))
\mn
\begin{enumerate}
\item[$(*)_4$]  for $B \subseteq M$ of cardinality $< \kappa$ let $C(B)$
be a subset of $M$ of cardinality $\le 2^{|B|+|T|}$ in which every $p
\in \bold S^{< \omega}(B)$ is realized.
\end{enumerate}
\mn
Lastly, let $F$ be defined for $B \in [M]^{< \kappa},F(B) =
\cup\{C_{\bold s}(\bar b_0,\bar b_1):\bold s$ a case and $\bar
b_0,\bar b_1$ suitable sequences from $B\} \cup C(B)$.

Now the number of cases is $\le \aleph_0 + \aleph_0 + \aleph_0 + |T| +
\theta + 2^{|T|+|A|} = 2^{|T|+|A|}$, so $F(B) \in [M]^{\le \theta}$.
So we are done.  
\end{PROOF}

\bigskip

\begin{theorem}
\label{tp25.33}
\underline{The Existence Theorem}  If $\ell = 0,1$ and 
$\cf(\kappa) >
\theta \ge |T|$ and $\bold x \in K^\ell_{\lambda,\kappa,\theta}$
\then \, there is $\bold y$ such that $\bold x \le_2 \bold y \in$ 
{\rm mxK}$^\ell_{\lambda,\kappa,\theta}$.
\end{theorem}

\begin{remark}
\label{tp25.34}  
1) If we use $K^\ell_{\lambda,\kappa,<\theta}$ instead of
``$\theta \ge |T|$" we should demand ``cf$(\theta) > |T|$".

\noindent
2) We may get more.  E.g. demand $I_1 = I_{\bold x},I_2$ is well
ordered and $I_{\bold y} = I_1 \cup I_2,I_1 <
I_2$, i.e. $s_1 \in I_1 \wedge s_2 \in I_2 \Rightarrow s_1 <_I s_2$. 

\noindent
2A) Also this claim holds (by the same proof) 
when we replace clause $(f)$ in Definition
\ref{tp14.21}(1) by
\mn
\begin{enumerate}
\item[$(f)_2$]   tp$(\bar{\bold c}_{\bold x},A_{\bold x})$ is
finitely satisfiable in $B_{\bold x}$.  
\end{enumerate}
\mn
Then in part (2) of the remark we may add 
\mn
\begin{enumerate}
\item[$(*)$]   for $t \in I_2$,
tp$(\bar c_t,M_{\bold x} \cup \{\bar c_{\bold y,s}:s <_{I_{\bold y}}
t\})$ is finitely satisfiable in $B_{\bold y}$.  
\end{enumerate}
\mn
See more in \ref{tp.77}, \ref{tp.84} and \ref{tp.98}.

\noindent
2B) In this case we may say ``Clause (f)$_1$, 
of \ref{tp14.21}(1)" instead of Clause (f).

\noindent
3) We can be more relaxed in the demands on $\langle \bold
x_\alpha:\alpha < \theta^+\rangle$ in the proof e.g. 
it suffices to demand
\mn
\begin{enumerate}
\item[$\circledast'$]    $(a) \quad \bold x_\alpha \in 
K^0_{\lambda,\kappa,\theta}$,
\sn
\item[${{}}$]   $(b) \quad \bold x_\alpha$ is $\le_1$-increasing
continuous, natural to demand ``$\le_2$-increasing", 

\hskip25pt  that is, $\bar{\bd}_{{\bold x}_\alpha} 
= \bar{\bd}_{\bold x_0}$ but not necessary, 
\sn
\item[${{}}$]   $(c) \quad$ for each $\alpha < \theta^+$ 
(or just for stationarily many $\alpha < \theta^+$) we have 

\hskip25pt $\bold x_\alpha <^+_1 \bold x_{\alpha +1}$.
\end{enumerate}
\mn
[Why?  Let $\bar\bd_{\bold x_1} \rest u_\alpha$ ($u_\alpha$ finite)
$t_\alpha,\varphi_\alpha(\bar x_\alpha,\bar y_\alpha,\bar b_\alpha)$
witness $\bold x_\alpha <^+_1 \bold x_{\alpha +1}$ when $\alpha \in
S_0$ where $S_0 := \{\alpha < 
\theta^+:\bold x_\alpha <^+_1 \bold x_{\alpha +1}\}$.  For $\alpha \in
S_0$ let $h(\alpha) = \Min\{\gamma:\bar\bd_{\bold x_\alpha} \rest
u_\alpha = \bar\bd_{\bold x_\gamma} \rest u_\alpha$, equivalently
$u_\alpha \subseteq \Dom(\bar{\bold d}_{\bold x_\gamma})$ and Rang$(\bar
b_\alpha) \cap (A_{\bold x_\alpha} + \bar{\bc}_{\bold x_\alpha})
\subseteq A_{\bold x_\gamma} + \bar{\bc}_{\bold x_\gamma}\}$, 
clearly $h(\alpha) < \alpha$ for $\alpha$ limit $\in S_0$.

So by Fodor's Lemma for some $\beta < \alpha$ and $u$
the set $S = \{\delta \in S_0:\delta$ is a
limit ordinal as in clause (c) above $u_\delta = u$
and $h(\delta) = \beta\}$ is stationary.  As $\theta \ge |T|$, for 
some\footnote{if $\bar c_{t_\alpha,n}$ is infinite we let $u_\alpha
\subseteq \ell g(\bar c^{\bold x_\alpha}_{t_{\alpha,0}})$ be finite
such that ${\frak C} \models \varphi_\alpha[\bar{\bold d}_{\bold
x},\bar c^{\bold x_{\alpha +1}}_{t_{\alpha,0}} \restriction
u_\alpha,\bar b_\alpha] \wedge \neg \varphi_\alpha
[\bar{\bd}_{\bold x_\alpha},\bar c^{\bold x_{\alpha +1}}_{t_\alpha,1},
\bar b_\alpha]$ and the rest is the same.} 
formula $\varphi$ the set $S_2 := \{\delta \in S_1:\varphi_\delta =
\varphi\}$ is a stationary subset of $\theta^+$ and we continue as in
the proof.

\noindent
4) How does part (3) of the remark help?  E.g. if we like to get
$\bold y \in \text{ mxK}^\ell_{\lambda,\kappa,\theta}$ which is normal
and Rang$(\bar{\bold d}_{\bold y})$ is the universe of some $N \prec \gC$.
\end{remark}

\begin{PROOF}{\ref{tp25.33}}
Assume this fails.  We try to choose $\bold x_\alpha$ by
induction on $\alpha < \theta^+$ such that
\mn
\begin{enumerate}
\item[$\circledast$]   $(a) \quad \bold x_\alpha \in
K^\ell_{\lambda,\kappa,\theta}$,
\sn
\item[${{}}$]   $(b) \quad \bold x_\beta \le_2 \bold x_\alpha$ for
$\beta < \alpha$,
\sn
\item[${{}}$]   $(c) \quad$ if $\alpha = \beta +1$ then $\bold x_\beta
<^+_2 \bold x_\alpha$, i.e.
$\bold x_\alpha$ witness $\bold x_\beta \notin 
\text{\rm mxK}^\ell_{\lambda,\kappa,\theta}$

\hskip25pt  (i.e. is like $\bold y$ in Definition \ref{tp25.32});
\sn
\item[${{}}$]   $(d) \quad \bold x_0 = \bold x$.
\end{enumerate}
\mn
For $\alpha = 0$ use clause (d), for $\alpha = \beta +1$ we use our
assumption toward contradiction.  For $\alpha$ limit use
\ref{tp14.28}(3).  Note that $\bar{\bd}_{\bold x_\alpha} =
\bar{\bold d}_{\bold x}$ for $\alpha < \theta^+$ by clauses (d) + (b).

Having carried the induction, for each $\alpha < \theta^+$ there are $t_\alpha,
\varphi_\alpha(\bar x,\bar y_\alpha,\bar z_\alpha),\bar b_\alpha$ 
satisfying:
\mn
\begin{enumerate}
\item[$(*)^1_\alpha$]   $(a) \quad t_\alpha \in I_{\bold x_{\alpha +1}} 
\backslash I_{\bold x_\alpha}$
\sn
\item[${{}}$]  $(b) \quad \varphi_\alpha(\bar x,\bar y_\alpha,\bar z_\alpha) 
\in \bbL(\tau_T)$
\sn
\item[${{}}$]  $(c) \quad \bar b_\alpha \subseteq M_{\bold x_\alpha} \cup
\bigcup\{\bar c_{\bold x_{\alpha +1},s}:s \in I_{\bold x_{\alpha +1}}
\backslash \{t_\alpha\}\}$ 
\sn
\item[${{}}$]  $(d) \quad \ell g(\bar b_\alpha) = \ell g(\bar
z_\alpha),\ell g(\bar y_\alpha) = \ell g(\bar c_{\bold x_{\alpha
+1},t_\alpha,0})$
\sn
\item[${{}}$]  $(e) \quad {\frak C} \models ``\neg \varphi_\alpha
[\bar{\bd}_{\bold x},\bar c_{\bold x_{\alpha +1},
t_\alpha,0},\bar b_\alpha] \wedge \varphi_\alpha[\bar{\bd}_{\bold x},
\bar c_{\bold x_{\alpha +1},t_\alpha,1},\bar b_\alpha]"$ 
so 

\hskip25pt $\ell g(\bar x) = \ell g(\bar{\bd}_{\bold x})$ 
and $n_{t_\alpha} \ge 2$.
\end{enumerate}
\mn
Clearly the
sequence $\langle t_\alpha:\alpha < \theta^+\rangle$ is without repetitions.
Now let $J_\alpha \subseteq I_{\bold x_{\alpha +1}} \backslash
\{t_\alpha\}$ be finite such that $\bar b_\alpha \subseteq M_{\bold
x_\alpha} \cup \bigcup\{\bar c_{x_{\alpha +1},s}:s \in
J_\alpha\}$.  We can find a pair $(\varphi(\bar x,\bar
y,\bar z),n_*,m_*)$ such that the set $S_0 = \{\delta <
\theta^+:|J_\delta| \le n_*$ and $\varphi_\delta(\bar x,\bar y_0,\bar
z_0) = \varphi(\bar x,\bar y,\bar z)\}$ is infinite (even stationary).
By Ramsey's theorem (or Fodor's Lemma) we can find an infinite
(and even stationary) set $S_1 \subseteq S_0 \subseteq \theta^+$
such that $\delta \in S_1 \Rightarrow J_\delta \cap\{t_\beta:\beta 
\in S_1\} = \emptyset$.  Note that there are
$\le\theta$ possibilities for $\varphi$, not necessarily $\le |T|$
because though
$\varphi_\alpha(\bar x,\bar y_\alpha,\bar z_\alpha)$ depend only on $\bar x
\restriction u$ for some finite $u \subseteq \ell g(\bar x)$ there are
$\le \ell g(\bar{\bd}_{\bold x}) + \aleph_0$ possibilities for $u$.

Next we shall prove that in this case
$\varphi(\bar x;\bar y,\bar z)$ has the independence
property (for $T$), a contradiction.

For every $w,v \subseteq S_1$ and $\eta \in {}^w 2$ let,
noting that $w = \Dom(\eta)$:
\mn
\begin{enumerate}
\item[$(*)_2$]   $(a) \quad A_{\eta,v} = A_{\bold x} \cup 
\{\bar c_{\bold x_{\beta +1},t,n}:t \in J_\beta$ and 
$n<n_{\bold x_{\beta +1},t}$ for some $\beta \in v\} \cup$

\hskip25pt  $\{\bar c_{\bold x_{\alpha +1},
t_\alpha,\eta(\alpha)}:\alpha \in w\}$,
\sn
\item[${{}}$] $(b) \quad f_{\eta,v}$ is the function with:
\sn
\begin{enumerate}
\item[${{}}$]   $(\alpha) \quad$ domain $A_{\eta,v}$,
\sn 
\item[${{}}$]   $(\beta) \quad$ is the identity on $A_{\bold x} \cup \{\bar
c_{\bold x_{\beta +1},t,n}:t \in J_\beta$ and

\hskip35pt  $n < n_{\bold x_{\beta +1},t}$ for some $\beta \in v\}$,
\sn 
\item[${{}}$]  $(\gamma) \quad f_{\eta,v}(\bar c_{\bold x_{\alpha +1},t_\alpha,
\eta(\alpha)}) = \bar c_{\bold x_{\alpha +1},t_\alpha,0}$ 
for $\alpha \in w = \text{ Dom}(\eta)$.
\end{enumerate}
\end{enumerate}
\mn
Now 
\mn
\begin{enumerate}
\item[$(*)_3$]   $f_{\eta,v}$ is an elementary mapping.
\end{enumerate}
\mn
[Why?  Without loss of generality $w \subseteq v$ are finite so
$\subseteq \alpha(*)$ for some $\alpha(*) < \theta^+$ 
and prove this by induction on $|v|$.  We just use:  for $\alpha \in v$ the
sequence $\langle \bar c_{\bold x_{d(*)},t,n}:
n < n_{\bold x_{\alpha(*)},t}\rangle$ is indiscernible over $A_{\bold
x_{\alpha(*)}} \cup \{\bar c_{\bold x_{\alpha(*)},s}:s \in I_{\bold
x_{\alpha(*)}} \backslash \{t\}\}$, by Definition
\ref{tp14.21}(1), Clause (e).]

Now let $g_{\eta,v} \in \text{\rm aut}({\frak C})$ extends $f_{\eta,v}$.
So $\alpha \in v \Rightarrow g_{\eta,v}(\bar b_\alpha) = 
f_\eta(\bar b_\alpha) = \bar b_\alpha$ and $g_{\eta,\nu}
(\bar c_{\bold x_{\alpha +1},t_\alpha,\eta(\alpha)}) = 
\bar c_{\bold x_{\alpha +1},t_\alpha,0}$ 
for $\alpha \in w$; hence by the choice of $J_\alpha$ so for $\eta
\in{}^w 2$ we have
\mn
\begin{enumerate}
\item[$(*)_4$]   ${\gC} \models \varphi[g_{\eta,v}(\bar{\bd}_{\bold x}),
\bar c_{\bold x_{\alpha +1},t_\alpha,0},\bar b_\alpha]$ iff 
$\eta(\alpha)=1$.
\end{enumerate}
\mn
So $\langle \varphi(\bar x,\bar c_{\bold x_{\alpha +1},t_\alpha,0},
\bar b_\alpha):t \in u\rangle$ is an independent sequence of formulas,
see \ref{0n.17}(a);
as $w$ is any subset of $S_1$ we get a contradiction as promised.   
\end{PROOF}

\begin{claim}
\label{tp25.36}
\underline{The Weak Orthogonality Claim}  

Assume $\ell=0,1$ and $\bold x \in 
\text{\rm mxK}^\ell_{\lambda,\kappa,\theta}$.

\noindent
1) If $m < \omega,B' \subseteq M_{\bold x}$ 
and\footnote{this is just to ensure that $M$ realizes
every $q \in \bold S^{< \omega}(B')$.} $|B'| < \kappa$  and $q \in
\bold S^m(M_{\bold x} \cup C_{\bold x})$ does not split over $B'$
\then \, {\rm tp}$(\bar{\bd}_{\bold x},
M_{\bold x} \cup C_{\bold x})$ is weakly orthogonal to $q$.

\noindent
1A) It suffices\footnote{even more as we can increase the linear order
 $I_{\bold x}$.}  that $q = \text{\rm tp}(\bar c_1 \char 94 
\bar c_2,M_{\bold x} \cup C_{\bold x})$ and $\tp(\bar c_1,M_{\bold x})$
does not split over some $B' \in [M]^{< \kappa}$, {\rm tp}$(C_{\bold
x},M_{\bold x} + \bar c_1)$ does not split over $B_{\bold x}$ and
{\rm tp}$(\bar c_2,M_{\bold x} + \bar c_1 + C_{\bold x})$ does not split
over some $B''' \in [M_{\bold x}]^{< \kappa}$.

\noindent
2) If ${\bold x}_\alpha \in \text{\rm mxK}^\ell_{\lambda,\kappa,\theta}$ for
$\alpha < \delta$ is $\le_1$-increasing, $\delta < \theta^+$ and 
${\bold x}_\delta := \cup\{{\bold x}_\alpha:\alpha <\delta\}$ 
belongs to $K^\ell_{\lambda,\kappa,\theta}$ (recall
\ref{tp14.28}(3)) \then \, it also belongs to 
{\rm mxK}$^\ell_{\lambda,\kappa,\theta}$.
\end{claim}

\begin{PROOF}{\ref{tp25.36}}
1) Assume toward contradiction that those types are not
weakly orthogonal.  Let $q = q(\bar y),\bar y = \langle y_k:k<m\rangle$ and let
$\bar x = \langle x_\alpha:\alpha < \alpha_{\bold x}\rangle$ recalling
$\alpha_{\bold x} = \ell g(\bar{\bd}_{\bold x})$ and $p(\bar x) =
\text{ tp}(\bar{\bd}_{\bold x},M_{\bold x} \cup C_{\bold x})$.   
So for some formula $\varphi(\bar
x,\bar y,\bar z)$ and $\bar e \in {}^{\ell g(\bar z)}(M_{\bold x} \cup
C_{\bold x})$ the type $p(\bar x) \cup q(\bar y)$ does not decide
$\varphi(\bar x,\bar y,\bar e)$, i.e. $r_{\bold t}(\bar x,\bar y) =
p(\bar x) \cup q(\bar y) \cup \{\varphi(\bar x,\bar y,\bar e)^{\bold t}\}$
is consistent (= finitely satisfiable in ${\frak C}$) for $\bold t =
0,1$ and let $\bar c'_0,\bar c'_1$ be such that $\bar{\bd} \char 94
\bar c'_{\bold t}$ realizes $r_{\bold t}(\bar
x,\bar y)$ for $\bold t = 0,1$.  Now it cannot be that tp$(\bar
c'_{\bold t},M_{\bold x} \cup C_{\bold x} \cup \bar{\bd}_{\bold x})$ does
not split over $B'$ for both $\bold t = 0,1$ (by \ref{3k.0.7}(1), as
every $p \in \bold S^{< \omega}(B')$ is realized in $M$ recalling $M$
is $\kappa$-saturated and $\kappa > |B'|$).  So choose $\bar c_0 \in
\{\bar c'_0,\bar c'_1\}$ such that the type tp$(\bar c_0,M_{\bold x} \cup
C_{\bold x} \cup \bar{\bd})$ splits over $B'$.

Now by \ref{3k.0.7}(1) there is $\bar c_1 \in {}^m{\frak C}$ such that
tp$(\bar c_1,M_{\bold x} \cup C_{\bold x} \cup \bar{\bd}_{\bold x} \cup
\bar c_0)$ extends $q(\bar y)$ and does not split over $B'$.  Hence
also tp$(\bar c_1,M_{\bold x} \cup C_{\bold x} \cup \bar{\bd}_{\bold x})$
does not split over $B'$, hence it is different from 
$\text{tp}(\bar c_0,M_{\bold x} 
\cup C_{\bold x} \cup \bar{\bd}_{\bold x})$.  We can continue and
choose $\bar c_n \, (n=2,3\ldots)$ realizing the complete type 
over $M_{\bold x} \cup C_{\bold x} \cup \bar{\bd}_x \cup 
\bigcup\{\bar c_k:k<n\}$ which extends $q$ and does
not split over $B'$.  Hence
\mn
\begin{enumerate}
\item[$(*)_0$]    for every $n < \omega$, tp$(\bar c_n,M_{\bold x}
\cup C_{\bold x} \cup \bigcup\{\bar c_k:k < n\})$ extend $q(\bar y)$
and does not split over $B'$.  
\end{enumerate}
\mn
So by \ref{3k.0.7}(2)
\mn
\begin{enumerate}
\item[$(*)_1$]   $\langle \bar c_n:n < \omega\rangle$ is an
indiscernible sequence over $M_{\bold x} \cup C_{\bold x}$.
\end{enumerate}
\mn
Also (by induction on $\gamma \le \omega$) by \ref{3k.0.7}(3) we have:
\mn
\begin{enumerate}
\item[$(*)_2$]   if $t \in I_{\bold x}$ then $\langle \bar c_{\bold
x,t,n}:n < n_t\rangle$ is an indiscernible sequence
over $M_{\bold x} \cup \{\bar c_{\bold x,s,m}:
s \in I \backslash \{t\}\} \cup \{\bar c_n:n < \gamma\}$.
\end{enumerate}
\mn
Now we define $\bold y = (M_{\bold y},B_{\bold y},\bar{\bc}_{\bold
y},\bar{\bd}_{\bold y})$ by
\mn
\begin{enumerate}
\item[$\circledast$]  $(a) \quad M_{\bold y} = M_{\bold x}$,
\sn
\item[${{}}$]   $(b) \quad B_{\bold y} = B_{\bold x} \cup B'$,
\sn
\item[${{}}$]   $(c) \quad \bar{\bd}_{\bold y} = \bar{\bd}_{\bold x}$,
\sn
\item[${{}}$]   $(d) \quad I_{\bold y} = I_{\bold x} \cup \{s(*)\}$
were $[t \in I_{\bold x} \Rightarrow t <_{I_{\bold y}} s(*)]$,
\sn
\item[${{}}$]  $(e) \quad \bar{\bc}_{\bold y} \restriction I_{\bold x}
= \bar{\bc}_{\bold x}$,
\sn
\item[${{}}$]  $(f) \quad n_{\bold y,s(*)} = \omega$ (or any number $\in
[2,\omega]$) and $\bar c_{\bold y,s(*),n} = \bar c_n$.
\end{enumerate}
\mn
Now $\bold y$ contradicts the assumption $\bold x \in 
\text{ mxK}^\ell_{\lambda,\kappa,\theta}$. 

\noindent
1A) Similarly (recalling that we can use $\{s_1\} + I_{\bold x}$).

\noindent
2) Easy. 
\end{PROOF}

The following claim is a crucial step toward proving the Type
Decomposition Theorem \ref{tp25.43}.
\begin{claim}
\label{tp25.38}  
If $\bold x \in \text{\rm mxK}^\ell_{\lambda,\kappa,\theta}$ 
(or just $\bold x \in K$ and $\circledast$ below) \then \, for every
$A \subseteq M_{\bold x}$ of cardinality $< \kappa$ and $\varphi =
\varphi(\bar x,\bar y,\bar z) \in \bbL(\tau_T)$ satisfying $\ell
g(\bar x) = \ell g(\bar{\bd}_{\bold x}),
\ell g(\bar z) = \ell g(\bar c_{\bold x}),\ell g(\bar y) = m$, 
there is $\psi(\bar x,\bar e,\bar{\bold c}_{\bold x}) \in
\text{\rm tp}(\bar{\bd}_{\bold x},M_{\bold x} \cup C_{\bold x})$ satisfying
$\bar e \in {}^{\omega >}M_{\bold x}$ such that 
$\psi(\bar x,\bar e,\bar c_{\bold x}) \vdash
\{\varphi(\bar x,\bar b,\bar{\bc}_{\bold x})^{\bold t}:\bar b \in 
{}^{\ell g(\bar y)} A$ and $\bold t \in \{0,1\}$ are such that ${\frak C}
\models \varphi[\bar{\bd}_{\bold x},\bar b,\bar{\bc}_{\bold x}]^{\bold
t}\}$ where
\mn
\begin{enumerate}
\item[$\circledast$]    if $q(\bar y) \in \bold S^{\ell g(\bar
y)}(M_{\bold x} \cup \bar{\bc}_{\bold x})$ is finitely 
satisfiable in some $A \subseteq M_{\bold x}$
of cardinality $< \kappa$ \then \, $q(\bar y)$ is weakly orthogonal
to {\rm tp}$(\bar{\bd}_{\bold x},M_{\bold x} \cup C_{\bold x})$. 
\end{enumerate}
\end{claim}

\begin{PROOF}{\ref{tp25.38}}
By \ref{tp25.36}, if $\bold x \in 
\text{ mxK}_{\lambda,\kappa,\theta}$ then $\circledast$ holds; hence
we can in any case assume $\circledast$.

Let $p(\bar x) = \text{ tp}(\bar{\bd}_{\bold x},M \cup 
C_{\bold x})$, so $\bar x = \langle x_i:i < \ell g(\bar{\bd}_{\bold x})
\rangle$ and recalling $\bar y = \langle y_\ell:\ell < m\rangle$ we 
define a set $r=r(\bar y)$ as follows:

\begin{equation*}
\begin{array}{clcr}
r(\bar y) := \{(\exists \bar x)(\varphi(\bar x,\bar y,
\bar c_{\bold x})^{\bold t} \wedge \psi(\bar x,\bar a,
\bar{\bc}_{\bold x})):&\bold t \in \{0,1\} \text{ and} \\
  &\psi(\bar x,\bar a,\bar{\bc}_{\bold x}) \in p(\bar x) \text{ and }
\bar a \in {}^{\omega >}(M_{\bold x})\}.
\end{array}
\end{equation*}
\mn
Now
\mn
\begin{enumerate}
\item[$\odot_1$]   $r(\bar y)$ is not finitely satisfiable in ${}^m A$.
\end{enumerate}
\mn
[Why?  If $r(\bar y)$ is finitely satisfiable in
${}^m A$, then there is an ultrafilter ${\cD}$ on ${}^m A$ such that
for every $\vartheta(\bar y,\bar a,\bar c_{\bold x}) \in r(\bar y)$,
the set $\{\bar b:\bar b \in {}^m A$ and 
$\models \vartheta[\bar b,\bar a,\bar c_{\bold x}]\}$ belongs to 
${\cD}$.  Let $q(\bar y) = \text{ Av}({\cD},M_{\bold x} \cup
\bar c_{\bold x})$, clearly $q(\bar y) \in \bold S^m(M_{\bold x} \cup
\bar{\bold c}_{\bold x})$ is finitely satisfiable in ${}^m A,A
\subseteq M$ and $|A| < \kappa$.
Let $\bar b^*  \in {}^m \gC$ 
realize $q(\bar y)$, so $(\psi(\bar x,\bar a,\bar{\bc}_{\bold x}) 
\in p(\bar x)) \wedge \bold t  \in \{0,1\} \Rightarrow 
{\frak C} \models (\exists \bar x)(\varphi(\bar x,\bar b^*,
\bar{\bc}_{\bold x})^{\bold t} \wedge \psi(\bar x,\bar a,
\bar{\bc}_{\bold x}))$.
Why?  This holds by the choices of $r(\bar y),{\cD}$ and $\bar
b^*$.  As $p(\bar x)$ is closed under conjunctions it follows that
$p(\bar x) \cup\{\varphi(\bar x,\bar b^*,\bar{\bc}_{\bold
x})^{\bold t}\}$ is finitely satisfiable in ${\frak C}$ for $\bold t=0,1$.
But this contradicts the assumption $\circledast$.]

Hence for some $n < \omega$ and $\psi_\ell(\bar x,\bar
a_\ell,\bar{\bc}_{\bold x}) \in p(\bar x)$ for $\ell < n$ we have
\mn
\begin{enumerate}
\item[$\odot_2$]   for no $\bar b \in {}^m A$ do we have $(\exists
\bar x)(\varphi(\bar x,\bar b,\bar{\bc}_{\bold x})^{\bold t}
\wedge \psi_\ell(\bar x,\bar a_\ell,\bar c_{\bold x}))$ for $\ell <
n,\bold t \in \{0,1\}$.
\end{enumerate}
\mn
Let $\psi(\bar x,\bar a,\bar{\bc}_{\bold x}) =
\bigwedge\limits_{\ell < n} \psi_\ell(\bar x,\bar a_\ell,\bar{\bc}_{\bold
x})$, so clearly
\mn
\begin{enumerate}
\item[$\odot_3$]   $(a) \quad \bar a \in {}^{\omega >}(M_{\bold x})$,
\sn
\item[${{}}$]   $(b) \quad \psi(\bar x,\bar a,\bar{\bc}_{\bold x})
\in p(\bar x)$,
\sn
\item[${{}}$]    $(c) \quad$ for no $\bar b \in{}^m A$ do we have
$\bigwedge\limits_{\bold t = 0}^{1} (\exists \bar x)(\varphi(\bar x,\bar
b,\bar{\bc}_{\bold x})^{\bold t} \wedge \psi(\bar x,\bar
a,\bar{\bc}_{\bold x}))$.
\end{enumerate}
\mn
So for every $\bar b \in {}^m A$ for some $\bold t = \bold t(\bar b)
\in \{0,1\}$ we have ${\frak C} \models \neg(\exists \bar
x)(\varphi(\bar x,\bar b,\bar{\bc}_{\bold x})^{\bold t} \wedge
\psi(\bar x,\bar a,\bar{\bc}_{\bold x}))$ hence $\psi(\bar x,\bar
a,\bar{\bc}_{\bold x}) \vdash \neg \varphi(\bar x,\bar
b,\bar{\bc}_{\bold x})^{\bold t(\bar b)}$.  As $\psi(\bar x,\bar
a,\bar{\bc}_{\bold x}) \in p(\bar x) = \text{ tp}(\bar{\bd}_{\bold x},
M_{\bold x} \cup \bold c_{\bold x})$ it follows that $\neg \varphi(\bar x,\bar
b,\bar{\bc}_{\bold x})^{\bold t(\bar b)} \in p(\bar x)$.

So $\psi(\bar x,\bar a,\bar{\bc}_{\bold x})$ is as required.
\end{PROOF}

\begin{claim}
\label{tp14.42}  
Assume $\bold x \in \text{\rm mxK}^\ell_{\lambda,\kappa,\theta}$.  
If $B' \subseteq M_{\bold x},|B'| 
< \kappa,\bar d \in {}^m(D_{\bold x}),\bar c \in {}^{\omega >}
(C_{\bold x}),\varphi = \varphi(\bar x,\bar y,\bar z),\ell g(\bar x) 
= \ell g(\bar d),\ell g(\bar z) = \ell g(\bar c)$ \then \,
 for some $\psi(\bar x,\bar y',\bar z')$ and $\bar e \in {}^{\ell
 g(\bar y')}(M_{\bold x})$ we have $\models \psi[\bar d,\bar e,\bar c]$ and

\[
\psi(\bar x,\bar e,\bar c) \vdash \{\varphi(\bar x,\bar b,\bar c)^{\bold
t}:\bar b \in {}^{\ell g(\bar y)}(B') \text{ and } {\frak C} \models
\varphi[\bar d,\bar b,\bar c]^{\bold t} \text{ and } \bold t \in \{0,1\}\}.
\]
\end{claim}

\begin{PROOF}{\ref{tp14.42}}
This just reformulates \ref{tp25.38}; see more in
\ref{tp26.7}(0). 
\end{PROOF}

Now at last
\begin{PROOF}{\ref{tp25.43}}
\underline{Proof of \ref{tp25.43}, The Type Decomposition Theorem}  

By \ref{tp14.28}(4) there is 
$\bold x \in \text{\rm K}^1_{\kappa,\le \theta}$ such that 
$\bold d_{\bold x} = \bar{\bd}$ and $M_{\bold x} = M$.  By the
Existence Theorem \ref{tp25.33} \wilog \, $\bold x \in 
\text{ mxK}^1_{\kappa,\le \theta}$. 
Clearly $(\bold P_{\bold x,\theta},\le_{\bold x,\theta})$ is a partial
order.   Assume that $\alpha(*) <
\kappa$ and $p_\alpha \in \bold P_{\bold x,\theta}$ for $\alpha <
\alpha(*)$.  Let $B = \bigcup\{\text{Dom}(p_\alpha):\alpha <
\alpha(*)\} \cup B_{\bold x}$, so $B \subseteq
M_{\bold x}$ has cardinality $< \kappa$.  Hence by \ref{tp14.42} 
for every $v \subseteq \ell g(\bar{\bc}_{\bold x})$ 
and $\varphi = \varphi(\bar x,\bar y,\bar z)$ satisfying 
$\ell g(\bar z) = \ell g(\bar{\bc}_{\bold x}),\ell g(\bar
x) = \ell g(\bar{\bd}_{\bold x}),\ell g(\bar y) < \omega$ there is
$\psi = \psi_\varphi(\bar x,\bar e_\varphi,\bar{\bc}_{\bold x}) \in  
\text{ tp}(\bar{\bd}_{\bold x},M \cup \bar{\bc}_{\bold x})$
where $\bar e_\varphi \in {}^{\ell g(\bar y)}(M_{\bold x})$ 
such that $\psi_\varphi(\bar x,
\bar e_\varphi,\bar{\bc}_{\bold x}) \vdash \{\varphi(\bar x,\bar
b,\bar{\bc}_{\bold x})^{\bold t}:\bar b \in {}^{\ell g(\bar y)}B$
and ${\frak C} \models \varphi[\bar{\bd}_{\bold x},\bar b,
\bar{\bc}_{\bold x}]^{\bold t}$ and $\bold t \in \{0,1\}\}$.

Let $A = \cup\{\text{Rang}(\bar e_\varphi):\varphi$ as above$\}$,
 clearly $|A| \le \theta$ and
let $p = \text{ tp}(\bar{\bd}_{\bold x},A \cup \bar{\bc}_{\bold x})$;
it is an upper bound, as required.
\end{PROOF}

\begin{discussion}
\label{tp14.44} 
The type decomposition theorems say
that we can analyze a type $p \in \bold S(M)$ in two steps; first
tp$(\bar c,M)$ does not split over some ``small" $B \subseteq M$.
Second, tp$(\bar d,\bar c + M)$ is like a type in the theory of trees
and lastly \wilog \, by \ref{tp14.28}(6) 
some initial segment of $\bar d$ realizes $p$.  As an
example, see \cite{Sh:877}:
\end{discussion}

\begin{exercise}
\label{tp14.45}  
Let $T = T_{\text{ord}}$ be Th$(\bbQ,<)$ and if $M \prec {\frak C}_T,M$ 
is $\kappa$-saturated, in the main case and $p \in \bold S(M)$ \then \,
\mn
\begin{enumerate}
\item[$(a)$]  $p$ induces a cut $\bar C_p$ of $M$ where $\bar C_p 
= \langle C_{p,1},C_{p,2}\rangle,C_{p,1} = \{a \in M:(a < x) \in p\}$
and $C_{p,2} = M \backslash C_{p,1}$,
\sn
\item[$(b)$]   now $\bar C_p$ has a pair $(\kappa_1,\kappa_2)$ of
cofinalities, that is $\kappa_1 = \text{ cf}(C_{p,1},<_M),\kappa_2 =
\text{ cf}(C_{p,2},>_M)$ and max$\{\kappa_1,\kappa_2\} \ge \kappa$,
\sn
\item[$(c)$]   now $p$ does not split over some subset $B$ of $M$ of
cardinality $\le \lambda < \kappa$ \Iff \, min$\{\kappa_1,\kappa_2\} \le
\lambda$,
\sn
\item[$(d)$]   for every $B \subseteq M$ of cardinality $< \text{
min}\{\kappa_1,\kappa_2\}$ for some 
$\varphi(x,\bar a) \in p$ we have $\varphi(x,\bar a) \vdash p \restriction
B$ (i.e. $p$ under $\vdash$ is $\min\{\kappa_1,\kappa_2\}$-directed);
in fact we can add that for some $a_1 \in C_{p,1}$ and $a_2 \in
C_{p,2}$ we have $\varphi(x,\bar a) = a_1 < x < a_2$, 
\sn
\item[$(e)$]   so for $\kappa = \text{ min}\{\kappa_1,\kappa_2\}$ we
have a decomposition which is trivial in some sense: either we
have $\bar c = <>$ or we have $\bar d = \bar c$,
\sn
\item[$(f)$]   if e.g. $\kappa_1 < \kappa_2$ and $B_1$ is an
unbounded subset of $C_{p,1}$ of cardinality $\kappa_1$ and $c,d$ realize
in ${\frak C}$ (where $M \prec {\frak C}$) the type $p$ and ${\frak C}
\models c < d$, \then \, tp$(c,M)$ is finitely satisfiable in $B_1$ and
for every $A \subseteq M$ of cardinality $< \kappa_2$ for some formula
$\varphi(x,\bar a) \in p$ we have $\varphi(x,\bar a) \wedge (c < x)
\vdash \text{ tp}(d,A \cup \{c\})$.
\end{enumerate}
\end{exercise}

\begin{discussion}
\label{tp14.46}  
Note:  if $T$ is stable and $\bold x \in 
\text{\rm mxK}_{\kappa,\theta}$ is normal \then \, 
$\bar{\bd}_{\bold x} \subseteq \text{\rm dcl(Rang}
(\bar{\bc}_{\bold x}))$ recalling $\dcl(A) = \{b:b$ is definable over
$A$, equivalently $\tp(c,A) = \tp(b,A) \Rightarrow c=b\}$.
\end{discussion}

\begin{claim}
\label{tp26.7}  
0) In \ref{tp14.42} if {\rm cf}$(\kappa) > 
\theta +|T|$ \then \, we can choose $\psi = \psi(\bar x,
\bar y',\bar z')$ such that it depends on $\bold x,\bar d,
\varphi(\bar x,\bar y,\bar z)$ but not on $B'$.

\noindent
0A) In \ref{tp14.42} if $\cf(\kappa) > 2^{\theta +|T|}$ 
\then \, we can fix also $q = \tp(\bar e,C_{\bold x} 
\cup \bar{\bd}_{\bold x})$.  If
$\cf(\kappa) > 2^{\theta+|T|+|B_{\bold x}|}$ then moreover we can fix $q =
\tp(\bar e,B_{\bold x} \cup \bc_{\bold x} \cup \bar\bd_{\bold x})$.

\noindent 
1) Assume that $\bold x \in K_0$ and
$\varphi = \varphi(\bar x,\bar y,\bar z) \in \bbL(\tau_T)$ are such that
$\bar y,\bar z$ are finite and
$\ell g(\bar x) = \ell g(\bar{\bd}_{\bold x})$.

\Then \, the following set is finite: 
$J = J_{\bold x,\varphi} = \{t \in I_{\bold x}$: there are
$\bar\alpha = \langle \alpha_\ell:\ell < \ell g(\bar y)\rangle$
and a sequence $\bar b \in {}^{\ell g(\bar y)}(A_{\bold x} \cup 
\{\bar c_{\bold x,s}:s \in I_{\bold x} \backslash \{t\}\})$ such that
$\alpha_\ell \le \ell g(\bar c_{\bold x,t,0})$
and ${\frak C} \models \varphi[\bar{\bold d}_{\bold x},\langle(\bar c_{\bold
x,t,0})_{\alpha_\ell}:\ell < \ell g(\bar y)\rangle,\bar b] \wedge \neg
\varphi[\bar d_{\bold x},
\langle(\bar c_{\bold x,t,1})_{\alpha_\ell}:\ell < \ell
g(\bar y)\rangle,\bar b)]\}$.

\noindent
2) Moreover, the bound depends just on $\varphi$ and $T$.

\noindent
3) For any $\bold x \in K_{< \aleph_0,< \aleph_0}$ and $\varphi$ there
is $\bold y$ satisfying $\bold x \le_2 \bold y \in K_{< \aleph_0,< \aleph_0}$
such that $I_{\bold y} \backslash I_{\bold x}$ is finite and the local version
of maximality holds, i.e.

\[
\bold y \le_2 \bold z \in K_{< \aleph_0,< \aleph_0} \Rightarrow
J_{\bold z,\varphi} = J_{\bold y,\varphi}.
\]
\mn
4) For any $\bold x \in K_{< \aleph_0,< \aleph_0}$ and sequence
$\langle \varphi_n:n < \omega\rangle$  we can find
$\langle \bold x_n:n < \omega\rangle$ such that $\bold x_0 = \bold
x,\bold x_n \le_2 \bold x_{n+1} \in K_{< \aleph_0,<\aleph_0},I_{\bold
x_n} \backslash I_{\bold x}$ finite and $\bold x_{n+1},
\varphi_n$ satisfies the demands on $\bold y,\varphi$ above for every $n$.
\end{claim}

\begin{PROOF}{\ref{tp26.7}}
0) As there are $\le \theta + |T| < \text{ cf}(\kappa)$
possible choices of $\psi$ and the set of possible $B$'s is
$(\text{cf}(\kappa))$-directed.

\noindent
0A) Similarly.

\noindent
1),2)  Let $n_1$ be minimal such that: for no $\bar
b_\ell \in {}^{\ell g(\bar y)} {\frak C},\bar c_\ell \in {}^{\ell
g(\bar z)} {\frak C}$ for $\ell < n_1$ is $\langle \varphi(\bar x,\bar
b_\ell,\bar c_\ell):\ell < n_1\rangle$ an independent sequence of
formulas, exist as $T$ is a dependent theory.  
Let $n_2$ be minimal such that if $u_i \in [n_2 \backslash
\{i\}]^{\ell g(\bar z)}$ for $i < n_2$ then for some $v \in
[n_2]^{n_1}$ we have $i,j \in v \Rightarrow i \notin u_j$ (the
$\Delta$-system lemma for finite sets, see \ref{0n.19}).  
Now $n_2$ is a bound as required by the proof of \ref{tp25.33}.

\noindent
3),4) Follows.
\end{PROOF}

\begin{conclusion}
\label{tp26.9}
If $\bold x \in \mxK^\ell_{\lambda,\kappa,\theta}$ and
$\cf(\kappa) > \theta + |T|$ \then \, we can find $\bar e \in
{}^\theta{\gC}$ and $\langle \psi_{\varphi(\bar x,\bar y,\bar z)}(\bar
x,\bar e,\bar{\bold c}_{\bold x}):\varphi(\bar x,\bar y,\bar z) \in
\bbL(\tau_T)\rangle$ satisfying
$\ell g(\bar x) = \ell g(\bar{\bd}_{\bold x}),\ell g(\bar z)
= \bar{\bc}_{\bold x},\ell g(\bar y) = \theta$ 
such that: for each $\varphi(\bar x,\bar y,\bar z)$ with $\ell g(\bar
x) = \ell g(\bar\bd_{\bold x}),\ell g(\bar y) = \theta,\ell g(\bar z)
= \bar\bc_{\bold x}$ (but $\varphi$ depends just on finitely many variables)
we have: 
\mn
\begin{enumerate}
\item[$\bullet$]  $\psi_{\varphi(\bar x,\bar y,\bar z)}
(\bar x,\bar e,\bar{\bc}_{\bold x}) \in
\tp(\bar{\bold d}_{\bold x},\bar e \cup \bar{\bc}_{\bold x})$,
\sn
\item[$\bullet$]  $\psi_{\varphi(\bar x,\bar y,\bar z)}
(\bar x,\bar e,\bar{\bc}_{\bold x}) \vdash
\{\varphi(\bar x,\bar b,\bar{\bc}_{\bold x})^{\bold t}:\bar b \in {}^{\ell
g(\bar y)}M$ and $\bold t \in \{0,1\}$ such that $\gC \models 
\varphi[\bar d_{\bold x},\bar b,\bar{\bc}_{\bold x}]^{\bold t}\}$
\sn
\item[$\bullet$]  $\tp(\bar e,M_{\bold x} \cup \bar{\bold c}_{\bold x} 
\cup \bar{\bold d}_{\bold x})$ is finitely satisfiable in $M_{\bold x}$.
\end{enumerate}
\end{conclusion}

\begin{PROOF}{\ref{tp26.9}}
By \ref{tp25.38} and \ref{tp26.7} and compactness.
\end{PROOF}

\noindent
On strongly/strongly$^2$ dependent theories see \cite{Sh:863}.
\begin{remark}
\label{tp26.14}  
1) If $T$ is strongly dependent \then \, in
the previous claim \ref{tp26.7}(1), if $\bold x$ satisfies
$t \in I_{\bold x} \Rightarrow
n_t = \omega$ \then \, for each $n < \omega$ even the set 
$\bold J_{\bold x} = \cup\{\bold J_{\bold x,\varphi}:\varphi(\bar
x,\bar y,\bar z) \in \bbL(\tau_T)$ as there and 
$\ell g(\bar z) < n\}$ is finite.

\noindent
2) If $Y$ is strongly dependent$^2$ \then \, above we can allow $n =
   \omega$.

\noindent
3) The proofs are similar.
\end{remark}
\bigskip

\centerline {$* \qquad * \qquad *$}
\bigskip

\noindent
We now turn to exact saturation.  We first prove Theorem
\ref{tp16.14}, second we prove in \ref{tp16.15} that some independent
$T$'s satisfies it, even give a sufficient criterion.  Third, we give
a sufficient condition for $T$ to satisfy the theorem - in
\ref{tp16.17} - the existence of a stable indiscernible set (\ref{tp16.16}).
\begin{PROOF}{\ref{tp16.14}}
\underline{Proof of \ref{tp16.14}, The Singular Exact Saturation Lemma}  

Let $\theta = |T|$.  As $M$ is not
$\kappa^+$-saturated, there are $A \subseteq M$ of cardinality $\le
\kappa$ and $p \in \bold S^1(A)$ omitted by $M$.  Let $d \in {\frak
C}$ realize $p$.  By Theorem \ref{tp25.43}, there is $\bar{\bc}
\in {}^{\theta \ge}{\gC}$ as there for (d), so in particular such that 
tp$(\bar{\bold c},M)$ does not split over some $N \prec M$ of
cardinality $< \kappa$.  Let $\langle B_i:i < 
\text{ cf}(\kappa)\rangle$ be a $\subseteq$-increasing sequence of
sets with union $N \cup A$ such that $i < \text{ cf}(\kappa) \Rightarrow |B_i|
< \kappa$ and $N \subseteq B_0$.  
Now we choose $A_i  \subseteq M$ by induction on $i <
\text{ cf}(\kappa)$ such that $A_i$ is of cardinality $\le \theta$ and
tp$(d,A_i \cup \bar{\bold c}) \vdash \text{ tp}(d,B^+_i)$
where $B^+_i := B_i \cup \bigcup\{A_j:j<i\}$.

The choice is possible as $|B^+_i| < \kappa$ by \ref{tp25.43},
i.e. by the choice of $\bar{\bold c}$.  Next
we can find $A_\kappa \subseteq M$ of cardinality $\le \theta$ such that
tp$(d,A_\kappa \cup \bar{\bold c}) \vdash \text{ tp}
(d,\bigcup\limits_{i < \kappa} A_i)$, possible as 
$|\bigcup\limits_{i < \text{ cf}(\kappa)} A_i| \le \theta + 
\text{ cf}(\kappa) < \kappa$.  Let $B^+ = \cup\{B^+_i:i < 
\text{ cf}(\kappa)\} \cup A_\kappa$ so $|B^+| = \kappa$ and 
we ask the question:
\mn
\begin{enumerate}
\item[$\odot$]   is there an elementary mapping 
$f$ (or automorphism of ${\frak C}$)
such that $f \restriction B^+$ is the identity and $f(\bar{\bc}) 
\in {}^\theta M$?
\end{enumerate}
\mn
If yes, then let $d' \in M$ realize $f(\text{tp}(d,A_\kappa 
\cup \bar{\bc}) \in \bold S(A_\kappa \cup f(\bar{\bc}))$ 
hence there is an elementary mapping $g$ satisfying 
$g \rest A_\kappa = \id_{A_\kappa} = f \rest A_\kappa,
g(\bar{\bc}) = f(\bar{\bc})$ and $g(d) = d'$.  Hence easily for each 
$i < \text{ cf}(\kappa)$ the sequence $\langle d' \rangle \char 94
f(\bar{\bc})$ realizes $f(\text{tp}(\langle d \rangle \char 94
\bar{\bc},A_i))$ hence it realizes also
$f(\text{tp}(\langle d \rangle \char 94 \bar{\bc},B^+_i)))$; so
$d'$ realizes $f(\tp(d,B^+_i))$. 
But as $B^+_i$ increases with $i$ it realizes
$f(\text{tp}(d,\cup\{B^+_i:i < \text{ cf}(\kappa)\})$, but $A
\subseteq \cup\{B_i:i < \cf(\kappa)\} \subseteq \cup\{B^+_i:i < \kappa\}$
 hence $d'$ realizes $\tp(d,A)$, but $d' \in M$ 
contradicting the choice of $p,A,d$.

Let $B^*_i = B^+_i \cup A_\kappa$ hence $\langle B^*_i:i <
\cf(\kappa)\rangle$ is $\subseteq$-increasing with union $B^+$ and
$|B_i| < \kappa$ for $i < \cf(\kappa)$.
So the answer to the question is no, which gives clauses (a),(b) of the desired
conclusion.  As for clause (c), we choose $\bar{\bold c}_\varepsilon$
by induction on $\varepsilon < \kappa$ such that:
\mn
\begin{enumerate}
\item[$(*)$]  $(a) \quad \bar{\bc}_\varepsilon \in {}^\theta M$,
\sn
\item[${{}}$]  $(b) \quad \bar{\bc}_\varepsilon$ realizes
tp$(\bar{\bc},\cup\{\bar{\bc}_\zeta:\zeta < \varepsilon\} \cup N)$,
\sn
\item[${{}}$]  $(c) \quad$ for even $\varepsilon$,
if possible $\bar{\bc}_\varepsilon$ does not
realize tp$(\bar{\bc},N \cup B^+)$ hence for some 

\hskip25pt $\alpha = \alpha_\varepsilon < \cf(\kappa),
\bar\bc_\varepsilon$ does not realize $\tp(\bar{\bc},B^*_\alpha)$,
\sn
\item[${{}}$]  $(d) \quad$ for even $\varepsilon$, if 
$\alpha_\varepsilon$ is well defined, it is minimal,
\sn
\item[${{}}$]  $(e) \quad$ for odd $\varepsilon,\alpha_\varepsilon =
  \alpha_{\varepsilon -1}$ and $\bar\bc_\varepsilon$ realizes
 $\tp(\bar\bc,B^*_\alpha)$ or even

\hskip25pt  $\tp(\bar\bc,\cup\{\bar\bc_\zeta:\zeta < \varepsilon\} \cup N \cup
  B^*_{\alpha_\varepsilon})$. 
\end{enumerate}
\mn
There is no problem in carrying the induction, by \ref{3k.0.7}(2) the
sequence $\langle \bar{\bc}_\varepsilon:\varepsilon < \kappa\rangle$ is
indiscernible over $N$.  Also obviously the sequence $\langle
\alpha_\varepsilon:\varepsilon <\kappa\rangle$ is non-decreasing; if
$\alpha_* = \cup\{\alpha_\varepsilon:\varepsilon <\kappa\}$ is equal
to $\cf(\kappa)$ we are done; that is letting $\bar{\bold c}_\kappa$
from ${}^\theta M$ be such that $\langle \bar{\bold
  c}_\varepsilon:\varepsilon \le \kappa\rangle$ is an indiscernible
sequence over $N$, necessarily $\bar{\bold c}_\kappa$ realizes
$\tp(\bar{\bold c}_1,N \cup B^*_i)$ for each $i < \kappa$ hence realizes
$\tp(\bar{\bold c},B^+)$ contradiction to ``$\odot$ fails".  But the
non-existence of $\bar{\bold c}_\kappa$ and the properties of $\langle
\bar{\bold c}_i:i < \kappa\rangle$ are as promised in clause (c) of
the Lemma, so we are indeed done.

Otherwise, for each $\varepsilon <
\kappa$ there is a formula $\varphi_\varepsilon(\bar x,\bar
e_\varepsilon) \in \tp(\bar\bc,B^*_{2 \varepsilon}) \subseteq
\tp(\bar\bc,B^*_\beta)$ with $\ell g(\bar x) = \ell g(\bar\bc)$ such
that $\gC \models \neg \varphi_\varepsilon[\bar\bc,\bar e_{2\varepsilon}]$,
  but by clause (e) of $(*)$
necessarily $\gC \models \varphi_\varepsilon
[\bar\bc_{2 \varepsilon +1},\bar e_{2 \varepsilon}]$.  
As $\kappa > |B^*_\beta| + |T|$ for some
  formula $\varphi(\bar x,\bar e)$ we have $\{\varepsilon <
  \kappa:\varphi_\varepsilon(\bar x,\bar e_\varepsilon) =\varphi(\bar
  x,\bar e_\varepsilon)\}$ is infinite.  But this contradicts $T$ being
  dependent, so we are done proving Clause (c).

Clause $(c)^+$, that is ``moreover, there is an ultrafilter ${\cD}$ on
$N$" follows when we use the version of $\mxK$ from \ref{tp25.34}(2A)
or \ref{tp.77} - \ref{tp.98} below.
\end{PROOF}

We may have hoped that \ref{tp16.14} characterize being dependent,
but this is not so.
Clarification when this property (characterization of exactly
$\kappa$-saturated $\kappa > \cf(\kappa) > |T|$, as in 
\ref{tp16.14}) holds is given by:

\begin{example}
\label{tp16.15} 
1) There is an independent $T$ such that: if
$T$ has an exactly $\kappa$-saturated model \then \, $\kappa$ is
regular.  In fact, this is a sufficient condition.

\noindent
2) The same holds for exactly $\kappa$-compact, $\kappa > \aleph_0$.
\end{example}

\begin{PROOF}{\ref{tp16.15}}
 We use $T$ which satisfies ``The Chang Trick" from his proof
of his two cardinal theorem $(\aleph_1,\aleph_0) \rightarrow
(\lambda^+,\lambda)$ when $\lambda = \lambda^{< \lambda}$; the use is
not an incident, he uses such $T$ to overcome a related problem in his proof.

The condition is:
\mn
\begin{enumerate}
\item[$\circledast$]   for some predicate $R(x,y) \in \tau_T$
written $xRy$ (or just a formula $\varphi_*(x,y) \in \bbL(\tau_T))$
we have\footnote{instead (a)+(b) we can have

$(a)' \quad (\forall x_0,\dotsc,x_{n-1})(\forall z)(z Ry =
\bigvee\limits_{\ell < n} z = x_\ell)$ for every $n$.}:
\begin{enumerate}
\item[$(a)$]  the empty set can be coded, that is $\exists y \forall
x(\neg xRy)$,
\sn
\item[$(b)$]   we can add to a coded set one element, that is $(\forall
x,y)(\exists y_1)(\forall x_1)[x_1 R y_1 \equiv (x_1 Ry \vee x_1 =x)]$.
\end{enumerate}
\end{enumerate}
\mn
Note: for any model $M$, if $R \notin \tau_M,M$ an 
infinite model, let $\langle u_b:b
\in M\rangle$ list the finite subsets of $M$, and we expand $M$ to $M^+$ by
choosing $R^{\mu^+} = \{(a,b):a \in u_b$ and $b \in M\}$, \then \,
Th$(M^+)$ is as required.

So assume $\kappa = \Sigma\{\kappa_i:i < \text{cf}(\kappa)\},\kappa_i
< \kappa_j < \kappa$ for $i < j < \text{ cf}(\kappa)$, 
cf$(\kappa) < \kappa$ and $M$ is $\kappa$-saturated.  Let
$A \subseteq M,|A| = \kappa$ and $p \in {\bold S}(A,M)$.  Let $A = 
\cup\{A_i:i < \text{ cf}(\kappa)\},|A_i| \le \kappa_i,A_i$ increasing
with $i$.  Let $c_i$ realize $p \restriction A_i$.

By induction on $i$, we choose $b_i \in M$ which realizes the type

\begin{equation*}
\begin{array}{clcr}
p_i(y) = \{c_j Ry:j < \text{ cf}(\kappa) \text{ and } j \ge i\} &\cup 
\{(\forall x)(xRy \rightarrow xR b_j):j<i\} \\
  &\cup\{(\forall x)(xRy \rightarrow \varphi(x,\bar a)):\varphi(\bar
x,\bar a) \in p \restriction A_i\}.
\end{array}
\end{equation*}
\mn
Arriving to the $i$-th stage by 
$\circledast$ and the induction hypothesis on $i,p_i(y)$ is
finitely satisfiable in $M$.

\noindent
[Why?  Let $p'(y) \subseteq p_i(y)$ be finite so it has the form $\{c_j
Ry:j \in u\} \cup \{(\forall x)(xRy \rightarrow xRb_j):j \in v\} \cup
\{(\forall x)(xRy \rightarrow \varphi_\ell(x,\bar a_\ell)):\ell < n\}$
where $u \subseteq [i,\text{cf}(\kappa))$ is finite, $v \subseteq i$ is
finite and $\varphi_\ell(\bar x,\bar a_\ell) \in p \rest A_i$ for 
$\ell < n$.  By
$\circledast$ we can find $c \in M$ such that $M \models (\forall x)
(xRc \equiv (\bigvee\limits_{j \in v} x = b_j)$, thus $c$
realizes $p'(y))]$.

But $|p_i| < \kappa$ so we can choose $b_i$.

Now $\{xRb_i:i < \text{ cf}(\kappa)\}$ is a set of formulas 
finitely satisfiable in $M$ of cardinality $< \kappa$ and any
element realizing it realizes $p$. 
\end{PROOF}

A weak complement to \ref{tp25.43} is \ref{tp16.17} but first recall:
\begin{definition}
\label{tp16.16}  [$T$ not necessarily dependent]. 

\noindent
1) $\bold I \subseteq {}^\alpha{\frak C}$ is a stable
indiscernible set \when \,: 
$\bold I$ is an infinite indiscernible set and
{\rm Av}$(\bold I,{\frak C})$ is well defined, i.e. for any $\varphi(\bar
x,\bar y),\ell g(\bar x) = \alpha$ and $\bar b \in 
{}^{\ell g(\bar y)}{\frak C}$ either $\varphi(\bold I,\bar b)$ or $\neg
\varphi(\bold I,\bar b)$ is finite.

\noindent
2) $\bold I$ is a dependent indiscernible sequence \when \,: 
$I$ is a linear order and $\bold I = \langle \bar a_t:t \in I\rangle$
is an indiscernible sequence and for every
formula $\varphi(\bar x,\bar b),\ell g(\bar x) = \ell g(\bar a_t)$,
there is a convex equivalent relation $E$ on $I$ with finitely many
equivalence classes such that $s Et \Rightarrow \varphi[\bar a_t,\bar
b] \equiv \varphi[\bar a_s,\bar b]$.
\end{definition}

\begin{fact}
\label{tp16.16n}
If $T$ is dependent, $\bold I \subseteq 
{}^\alpha {\frak C}$ is a stable indiscernible set
\Iff \,  $\bold I$ is an infinite indiscernible set.
\end{fact}

\begin{PROOF}{\ref{tp16.16n}}
By \cite[1.28]{Sh:715}.
\end{PROOF}

\begin{claim}
\label{tp16.17} 
1) Assume ($T$ is dependent and) there is
an infinite indiscernible set $\bold I \subseteq {\frak C}$.  If $\kappa^+ =
2^\kappa$ and $\kappa > |T|$ \then \, $T$ has an exactly $\kappa$-saturated
model.

\noindent
2) Assume ($T$ not necessarily dependent), $\bold I \subseteq {\frak
C}$ is a stable indiscernible set.  \Then \, the conclusion of part (1) holds.

\noindent
3) In parts (1),(2) the conclusion holds for $T$ if the assumption
   holds for $T^{\eq}$.
\end{claim}

\begin{remark}
\label{tp16.18}  
1) Of course, trivially if for some non-zero ordinal
$\alpha$ there is an infinite indiscernible set $\bold I \subseteq {}^\alpha
{\frak C}$ then for some $i < \alpha,\{(\bar a)_i:\bar a \in \bold I\}$
is an infinite indiscernible set.

\noindent
2) But we could use below indiscernible set $\bold I \subseteq
{}^\alpha {\frak C}$.

\noindent
3) On indiscernible sets for $T$ dependent see \cite[\S1]{Sh:715}, 
we use it freely.

\noindent
4) Of course: if $\{\bar a_t:t \in I\} \subseteq {}^\alpha{\frak C}$
is an infinite indiscernible set/a stable indiscernible set, $u
\subseteq \alpha$ and $\langle \bar a_t \restriction u:t \in I\rangle$
is not constant \then \, $\{\bar a_t \restriction u:t \in I\}
\subseteq {}^u{\frak C}$ is an infinite indiscernible set/a stable
indiscernible set.  So using singletons in \ref{tp16.17} is not a
loss.

\noindent
5)  Recall that if 
$I$ is a linear order and $I_1,I_2 \subseteq I$ are infinite,
$\bar a_t \in {}^{\beta}{\frak C}$ for $t \in I$ and $\bold I =
\langle \bar a_\gamma:t \in I\rangle$ is an indiscernible sequence
then $\bold I \restriction I_1$ is stable iff $\bold I \restriction
I_2$ is stable.

\noindent
This is easy.
\end{remark}


\begin{PROOF}{\ref{tp16.17}} 

\noindent
1) By part (2) (and \ref{tp16.16n}).

\noindent
2) Let $\bold I = \{a_\alpha:\alpha < \kappa\} \subseteq {\frak C}$
be an infinite stable indiscernible set.  Now easily
\mn
\begin{enumerate}
\item[$\odot_1$]   for any $\bar b \in {}^{\omega >} {\frak C}$ for
some $\bold J \in [\bold I]^{\le|T|}$, the set $\bold I \backslash
\bold J$ is an indiscernible set over $\bar b \cup \bold J$.
\end{enumerate}
\mn
[Why?  By the definition but we elaborate.  
First prove that for any $n$ and $\varphi =
  \varphi(\bar x_0,\dotsc,\bar x_{n-1},\bar b')$ for some finite
$\bold J_\varphi$ we have $\gC \models ``\varphi[a_0,\dotsc,a_{n-1},\bar b'] 
\equiv \varphi[a'_0,\dotsc,a'_{n-1},\bar b']"$ when
$a_0,\dotsc,a_{n-1},a'_0,\dotsc,a'_{n-1} \in \bold I \backslash 
\bold J_\varphi$ with no repetitions.  Second, use transitivity of
equivalence to show it suffices that $a_0,\dotsc,a_{n-1} \in \bold I
\backslash \bold J_\varphi$ with no repetition and
$a'_0,\dotsc,a'_{n-1} \in \bold I \backslash \bold J_\varphi$ with no
repetitions.  Lastly, we choose $\bold J_k \subseteq \bold I$ by induction 
on $k$ such that $\|\bold J_k\|
\le |T|,m < k \Rightarrow \bold J_m \subseteq \bold J_k \subseteq
\bold I$ and if
$k=m+1$ and $\varphi = \varphi(\bar x_0,\dotsc,\bar x_{n-1},\bar
b'),\bar b' \subseteq (\cup \bold J_m) \cup \bar b$ then we can above
choose $\bold J_\varphi \subseteq \bold J_k$.  Now $\bigcup\limits_{k}
\bold J_k$ is as required.]
\mn
\begin{enumerate}
\item[$\odot_2$]   the following conditions on $\bar b \in 
{}^{\omega >} {\frak C}$ are equivalent:
\begin{enumerate}
\item[$(a)$]  tp$(\bar b,\bold I)$, Av$(\bold I,\bold I)$ are
weakly orthogonal,
\sn
\item[$(b)$]  for some $\bold J \in [\bold I]^{\le|T|}$ we have
tp$(\bar b,\bold J) \vdash \text{ tp}(\bar b,\bold I)$.
\end{enumerate}
\end{enumerate}
\mn
[Why?  Easy using $\odot_1$.]
\mn
\begin{enumerate}
\item[$\odot_3$]  let $\bold D = \bold D_{\bold I} = \{p \in \bold
S^{< \omega}(\bold I):p$ weakly orthogonal to Av$(\bold I,\bold I)\}$.
\end{enumerate}
\mn
We define
\mn
\begin{enumerate}
\item[$\odot_4$]   we say $A$ is a $\bold D$-set if $\bar
a \in {}^{\omega >} A \Rightarrow \text{ tp}(\bar a,\bold I) \in 
\bold D_{\bold I}$
\sn
\item[$\odot_5$]   if $\bold I \subseteq A$ we let $\bold S^m_{\bold
D}(A) = \{\text{tp}(\bar b,A):A \cup \bar b$ is a $\bold D$-set and
$\ell g(\bar b) = m\}$.
\end{enumerate}
\mn
We note
\mn
\begin{enumerate}
\item[$\odot_6$]   $(a) \quad$ if $A$ is a $\bold D$-set and $\bold I \subseteq
A$ then Av$(\bold I,\bold I) \vdash \text{ Av}(\bold I,A)$; hence
if $A=|M|$, 

\hskip25pt then $M$ is not $\kappa^+$-saturated,
\sn
\item[${{}}$]  $(b) \quad A$ is a $\bold D$-set iff $A \cup \bold I$ is a
$\bold D$-set.
\sn
\item[$\odot_7$]   If $A$ is a $\bold D$-set of cardinality $<
\kappa$ and $p \in \bold S^m_{\bold D}(A \cup \bold I)$ \then \, for some
$\bold J \subseteq \bold I$ of cardinality $\le |A|+|T|$ we have $p
\restriction (A \cup \bold J) \vdash p$.
\end{enumerate}
\mn
[Why?  By $\odot_2$ and $\odot_6$.]
\mn
\begin{enumerate}
\item[$\odot_8$]   if $\langle A_\alpha:\alpha < \delta\rangle$ is an
$\subseteq$-increasing sequence of $\bold D$-sets, \then \, 
$A_\delta := \cup\{A_\alpha:\alpha < \delta\}$ is a $\bold D$-set.
\end{enumerate}
\mn
[Why?  By the definition of a $\bold D$-set.]
\mn
\begin{enumerate}
\item[$\odot_9$]   if $\langle A_\alpha:\alpha \le \delta\rangle$ is
$\subseteq$-increasing continuous sequence of $\bold D$-sets, 
$\bold I \subseteq A_0$ and $p \in
{\bold S}^m(A_\delta)$ then $p \in \bold S^m_{\bold D}(A_\delta)
\Leftrightarrow \bigwedge\limits_{\alpha < \delta} p \restriction
A_\alpha \in \bold S^m_{\bold D}(A_\alpha)$.
\end{enumerate}
\mn
[Why?  By the definition of $\bold S^m_{\bold D}(-)$ and $\odot_8$.]

Now comes a major point
\mn
\begin{enumerate}
\item[$\odot_{10}$]  if $A \subseteq {\frak C}$ and $|A| < \kappa$
\then \, we can find $\bold I_1$ such that: ${\bold I}_1$ is an
indiscernible set, $\bold I \subseteq \bold I_1,|\bold I_1 \backslash
\bold I| \le |T| + |A|$ and $A$ is a $\bold D_{\bold I_1}$-set.
\end{enumerate}
\mn
[Why?  Let $\theta = |A| + |T|$ and we try by induction on $\alpha <
\theta^+$ to choose an element $a_\alpha$ of ${\frak C}$ which
realizes Av$(\bold I,\{a_\beta:\beta < \alpha\} \cup \bold I)$ but
$\alpha$ is even iff $a_\alpha$ realizes 
$\Av(\bold I,\cup\{a_\beta:\beta < \alpha\} \cup A \cup 
\bold I)$.  But if we succeed to carry the induction clearly $\bold I^+ :=
\bold I \cup \{a_\alpha:\alpha < \theta^+\}$ is an indiscernible set,
and a stable one (recalling \ref{tp16.18}(5)) hence for some 
$\bold J \subseteq \bold I^+$ of
cardinality $\le |A| + |T|$, also $\bold I^+ \backslash \bold J$ is an
indiscernible set over $A \cup \bold J$, but necessarily $\bold J
\subseteq \bold I \cup \{a_\beta:\beta < \alpha\}$ for some $\alpha <
\theta^+$, easy contradiction to the choice of the $a_\alpha$'s.]
\mn
\begin{enumerate}
\item[$\odot_{11}$]   if $A_1 \subseteq A_2$ are 
$\bold D_{\bold I}$-sets, $|A_2| \le \kappa,|A_1| < \kappa$ 
and $p \in \bold S^m_{\bold D}
(A_1 \cup \bold I)$ then there is $q \in \bold S^m_{\bold D}
(A_2 \cup \bold I)$ extending $p$.
\end{enumerate}
\mn
[Why $\odot_{11}$?  By $\odot_9$ \wilog \, $A_2 = A_1 \cup \{b\}$, so
$|A_2| < \kappa$.  We can find $\bar c$ realizing $p(\bar y)$ and let
$A = A_1 \cup \{b\} \cup \bar c = A_2 \cup \bar c$.  So by $\odot_{10}$
there is $\bold I^+$ such that: $\bold I^+$ is an indiscernible set,
$\bold I \subseteq \bold I^+,|\bold I^+ \backslash \bold I| \le \theta
:= |A| + |T| = |A_1| + |T| < \kappa$ and $A$ is a $\bold D_{\bold
I^+}$-set.  As $A_2$ is a $\bold D_{\bold I}$-set we can find $\bold
J_1 \subseteq \bold I$ of cardinality $\le |A_2| + |T| \le \theta <
\kappa$ satisfying tp$(A_2,\bold J_1) \vdash \text{ tp}(A_2,\bold I)$.  Also
$A_1 \cup \bar c$ is a $\bold D_{\bold I}$-set (as $A_1$ is a $\bold
D_{\bold I}$-set and $\bar c$ realizes $p(\bar y) \in 
{\bold S}^m_{\bold D}(A_1 \cup \bold I))$ hence there is $\bold J_2 \subseteq
\bold I$ of cardinality $\le |A_1 \cup \bar c| + |T| = \theta <
\kappa$ such that tp$(A_1 \cup \bar c,\bold J_1) \vdash \text{ tp}(A_1
\cup \bar c,\bold I)$.  Lastly, as $A$ is a $\bold D_{\bold I^+}$-set
there is $\bold J_3 \subseteq \bold I^+$ of cardinality $\le |A| + |T|
= \theta$ such that tp$(A,\bold J_3) \vdash \text{ tp}(A,\bold I^+)$.

As tp$(A_2,\bold J_1) \vdash \text{ tp}(A_2,\bold I)$ necessarily
tp$(A_2,\bold J_1) \vdash \text{ tp}(A_2,\bold I^+)$.  Similarly
tp$(A_1 \cup \bar c,\bold J_2) \vdash \text{ tp}(A_2 \cup \bar c,\bold
I^+)$.  By cardinality considerations there is a permutation $h$ of
$\bold I^+$ which is the identity on $\bold J_1,\bold J_2$ and
$\bold J_3 \cap\bold I$ and maps $\bold J_3 \backslash \bold I$ into
$\bold I$.  As $\bold I^+$ is an indiscernible set, $h$ is an
elementary mapping (of ${\frak C}$).  As $h \restriction \bold J_1$ is the
identity and tp$(A_2,\bold J_1) \vdash \text{ tp}(A_2,\bold I^+)$, see
above also $h \cup \text{ id}_{A_2}$ is an elementary mapping hence
there is an automorphism $g$ of ${\frak C}$ extending $h \cup \text{
id}_{A_2}$.  As tp$(A_1 \cup \bar c,\bold J_2) \vdash \text{ tp}(A_1
\cup \bar c,\bold I^+)$ and $h \restriction \bold J_2 = 
\text{ id}_{\bold J_2},h \rest A_1
= \text{ id}_{A_1}$ (recalling $A_1 \subseteq A_2$) and $h(\bold I^+)
= \bold I^+$ necessarily $g(\bar c)$ realizes tp$(\bar c,A_1 \cup
\bold I^+)$ hence it realizes tp$(\bar c,A_1 \cup \bold I)$ which is
equal to $p(\bar y)$.  Also tp$(A_2 \cup \bar c,\bold J_3) \vdash
\text{ tp}(A_2 \cup \bar c,\bold I^+)$ hence tp$(A_2 \cup g(\bar
c),h(\bold J_3)) \vdash \text{ tp}(A_2 \cup h(\bar c),\bold I^+)$ hence
tp$(A_2 \cup g(\bar c),h(\bold J_3)) \vdash \text{ tp}(A_2 \cup g(\bar
c),\bold I)$, but $h(\bold J_3) \subseteq \bold I$.
So $A_2 \cup g(\bar c)$ is
a $\bold D_{\bold I}$-set hence $q(\bar y) =: \text{ tp}(g(\bar c),A_2
\cup \bold I)$ belongs to $\bold S^m_{\bold D_{\bold I}}(A_2 \cup
\bold I)$ so is as required.]

As $2^\kappa = \kappa^+$, by $\odot_9 + \odot_{11}$ there is $M \supseteq \bold
I$ of cardinality $\kappa^+$ which is $\kappa$-saturated and is a $\bold
D$-set hence by $\odot_6$ is not $\kappa^+$-saturated.

\noindent
3) Should be clear. 
\end{PROOF}

\begin{observation}
\label{tp16.19}  (Any complete first order $T$)

In ${\frak C}$ there is no infinite indiscernible set \underline{iff} for
some $n$ and $\varphi = \varphi(x_0,\dotsc,x_{n-1}) \in 
\bbL(\tau_T),\varphi$ is connected and anti-symmetric i.e. if
$a_0,\dotsc,a_{n-1} \in {\frak C}$ with no repetitions then for some
permutations $\pi_1,\pi_2$ of $\{0,\dotsc,n-1\}$ we have

\[
{\frak C} \models \varphi[a_{\pi_1(0)},\dotsc,a_{\pi_1(n-1)}] \wedge
\neg \varphi[a_{\pi_2(a)},\dotsc,a_{\pi_2(a_{n-1})}].
\]
\end{observation}

\begin{remark}  
\label{tp16.20}
1) The second condition is related to the property (E) of complete
 first order theories of Ehrenfeucht
\cite{Eh57} which says that the condition holds for some infinite set.

\noindent
2) Note that ${\frak C}$ may have no infinite indiscernible set but
${\frak C}^{\text{eq}}$ has.
\end{remark}

\begin{PROOF}{\ref{tp16.19}}
The implication $\Leftarrow$ is obvious.

So assume the first statement.  For $\alpha \le \omega$ and $\Delta
\subseteq \Delta_* := \{\varphi(\bar x):\varphi \in \Bbb L(\tau_T),\bar
x = \langle x_\ell:\ell < n\rangle\}$ let\footnote{yes: we use
  singletons $y$'s.} $\Gamma^\alpha_\Delta = \{y_k
\ne y_\ell:k < \ell < \alpha\} \cup
\{\varphi(y_{k_0},\dotsc,y_{k_{n-1}}) \equiv \varphi(y_{\ell_0},
\dotsc,y_{\ell_{n-1}}):n < \alpha,
\varphi(x_0,\dotsc,x_{n-1}) \in \bbL(\tau_T)$ and
$k_0,\dotsc,k_{n-1} < \alpha$ without repetitions and
$\ell_0,\dotsc,\ell_{n-1} < \alpha$ without repetitions$\}$.  Easily
$\Gamma^\omega_{\Delta_*}$ is not realized in ${\gC}$ by the
present assumption and $\langle \Gamma^k_\Delta:\Delta$ is a finite
subset of $\Delta_*$ and $k < \omega\rangle$ is $\subseteq$-increasing
with $k$ and $\Delta$ with union $\Gamma^\omega_{\Delta_*}$.  
Hence for some finite $\Delta \subseteq \Delta_*$ and 
$k < \omega$, the set $\Gamma^k_\Delta$ is not
realized in ${\gC}$.

Let $\langle \varphi_i(x_0,\dotsc,x_{n_i-1}):i < i(*)\rangle$ list
$\Delta$ so $i(*) < \omega$, so without loss of generality 
$n_i < k$ for $i < i(*)$.  Lastly, we define $\varphi(y_0,\dotsc,y_{k-1})$, 
it says: if $\langle y_\ell:\ell < k\rangle$ is without 
repetitions and $i$ is minimal such
that $\langle y_\ell:\ell < k\rangle$ is not a
$\{\varphi_i(x_0,\dotsc,x_{n_i-1})\}$-indiscernible set then
$\varphi_i(y_0,\dotsc,y_{n_i-1})$.  Now check.
\end{PROOF}

\begin{question}
\label{tp16.21}  
1) Is there a dependent $T$ such that even in 
${\frak C}^{\text{eq}}$ there is no infinite indiscernible set but some
singular $\kappa$ of cofinality $>|T|$ there is an exactly
$\kappa$-saturated model of $T$?  

\noindent
2) For a dependent theory $T$ 
characterize $\{\kappa:\kappa$ singular and $T$ has
exactly $\kappa$-saturated models$\}$.

\noindent
3) In both parts we may at least initially restrict ourselves to $\kappa$
strong limit of large enough cofinality such that $2^\kappa =
\kappa^+$.

\noindent
4) Try to eliminate the assumption ``$2^\kappa = \kappa^+$" in
\ref{tp16.17} at least when $\kappa$ is strong limit of cofinality
$> |T|$.  A natural way is via $\cP^-(n)$-diagrams (as in
\cite[Ch.XII]{Sh:c} and even closer in \cite{Sh:234}.
\end{question}
\bigskip

\noindent
\centerline{$* \qquad * \qquad *$}
\bigskip

\noindent
The following in a sense gives a spectrum for $\bar{\bold d}/M$.
\begin{claim}
\label{tp35.46} 
For $\theta \ge |T|$, a model $M$ and sequence
$\bar{\bold d} \in {}^{\theta \ge}{\frak C}$, 
there is a set $\Theta$ such that:
\mn
\begin{enumerate}
\item[$(a)$]   $\Theta \subseteq \Theta_* := \{\kappa:\kappa >
\theta$ and $M$ is $\kappa$-saturated$\}$,
\sn
\item[$(b)$]  $|\Theta| \le \theta$,
\sn
\item[$(c)$]    if $\kappa \in \Theta$ and {\rm cf}$(\kappa) >
\theta$ then there is $\bold x \in \mxK_{\kappa,\le \theta}$ such 
that $M_{\bold x} = M,\bar{\bold d}_{\bold x} =\bar{\bold d}$ and
$|B_{\bold x}| \le \theta + \text{\rm sup}(\Theta \cap \kappa)$,
\sn
\item[$(d)$]   if $\theta < \kappa \in \Theta$ and {\rm cf}$(\kappa) \le
\theta$ then $\text{\rm sup}(\Theta \cap \kappa) = \kappa$.
\end{enumerate}
\end{claim}

\begin{PROOF}{\ref{tp35.46}}
Straight.  

For each $\kappa \in \Theta' := \{\kappa' \in \Theta_*:\kappa'$ has
cofinality $> \theta\}$ we can find $\bold x_\kappa \in \text{
mxK}_{\kappa,\le \theta}$ such that $M_{\bold x} = M,
\bar{\bold d}_{\bold x} = \bar{\bold d}$ and for $\kappa \in \Theta_*
\backslash \Theta'$ let $\Theta_\kappa$ be a cofinal subset of $\kappa$
of cardinality cf$(\kappa) \le \theta$.  Let $f:\Theta' \rightarrow$ Card be
defined by $f(\kappa) = |B_{\bold x}| + \theta$.  Note that $\Theta'$
has a maximal member or $\Theta'$ has a cofinal subset of cardinality $\le
\theta$.  Now we shall 
choose $\Theta_n$ by induction on $n$ such that
$\Theta_n \subseteq \Theta_*,|\Theta_n| \le \theta$ and $n=m+1
\Rightarrow \Theta_m \subseteq \Theta_n$.  Let $\Theta_0$ be a cofinal
subset of $\Theta_*$ of cardinality $\le \theta$, see above why
possible.  If $n=m+1$, let $\Theta_n = \{f(\kappa):
\kappa \in \Theta_m \cap \Theta'\} \cup\{\Theta_\kappa:\kappa \in
\Theta_m \backslash \Theta'\} \cup \Theta_m$.  Now
$\cup\{\Theta_n:n < \omega\}$ is as required.
\end{PROOF}

\begin{discussion}
\label{tp26.47}  
Note that $\bold P_{\bold x}$ in
\ref{tp25.43} is $\kappa$-directed, but in general it is not definable in
$M_{\bold x}$ and even not definable in $(M_{\bold x})_{[\bar{\bold
c}_{\bold x}]}$ (or $M_{[B_{\bold x} + {\bold c}_{\bold x}]})$ even
by disjunction of types as it depends on $\bar{\bold
d}_{\bold x}$.  So we may consider $\bold P'_{\bold x} = \bold
P_{M,\bar{\bold c}_{\bold x},\alpha_{\bold x}} = \{p:p \in \bold
S^{\alpha(\bold x)}(A \cup \bar{\bold c}_{\bold x})$ and $A \subseteq
M$ has cardinality $< \kappa\}$ ordered as before.  
Now $\bold P'_{\bold x}$ is partially ordered but it is
not clear that it is $\kappa$-directed.  Moreover 
$(M_{\bold x})_{[\bar{\bold c}_{\bold x}]}$ is 
not $\kappa$-saturated, but is $(\bold
D_{\bold x},\kappa)$- sequence homogeneous for suitable $\bold
D_{\bold x}$ and $\bold D_{\bold x}$ is a good diagram (see e.g.
\cite{Sh:3}; see more in \cite{Sh:950}).  
So we can consider the families of such $\bold D$'s,
fixing ($T$ and) $\theta$.

But we can define the order in the $\kappa$-saturated 
$(M_{\bold x})_{[B_{\bold x}]}$ which is 
$\Bbb L_{\infty,\kappa}(\tau_T)$-equivalent to ${\frak C}_{[B_{\bold x}]}$.
In this model we have $\psi(\bar x,\bar y) \in \Bbb L_{\infty,\kappa}$
which is a partial order on the $\theta$-tuples, $\ell g(\bar x) =
\alpha_{\bold x} = \ell g(\bar y)$.

However, in our case we know more.  Letting $\Gamma = \Bbb L(\tau_T)$,
if cf$(\kappa) > \theta \ge |T|$ we know that we can 
find $\bar \psi = \langle \psi_\varphi(\bar x_\varphi,\bar
y_\varphi,\bar z_\varphi):\varphi \in \Gamma\rangle$ and the order 
on the set of $\bar{\bold e} = (\ldots \char 94 \bar e_\varphi \char
94 \ldots)_\varphi$ such that in \ref{tp25.38} we can choose $\psi =
\psi_\varphi$ (easy, see \ref{tp26.7}(0)).  
If cf$(\kappa) > 2^\theta$, we can fix
there also the type of $\bar{\bold e}$ over $C_{\bold x} 
\cup \bar{\bold d}_{\bold x}$.

So
\mn
\begin{enumerate}
\item[$(*)_1$]    let $\bold I = \{\bar e:\bar e$ as above$\}$, so
$\bold I$ is type-definable in $(M_{\bold x})_{[B_{\bold x}]}$.
\sn
\item[$(*)_2$]    $p_{\bar{\bold e}} = \{\psi_\varphi(\bar x_\varphi,
\bar{\bold e},\bar{\bold c}_{\bold x}):\varphi \in \Gamma\}$ for
$\bar{\bold e} \in \bold I$.
\sn
\item[$(*)_3$]    $(a) \quad \le_1$ defined by $\bar{\bold e}_1 \le
\bar{\bold e}_2$ if $p_{\bar{\bold e}_2} \vdash p'$ for some $p'$ such
that

\hskip25pt  $p_{\bar{\bold e}_1} \subseteq p' \in \bold S(B_{\bold x} \cup
\bar{\bold e}_1)$,
\sn
\item[${{}}$]   $(b) \quad \bold R$ is defined by $\bold e_1 R
\bold e_2$ \Iff \,  $\psi_\varphi(\bar x,\bar e_\varphi,\bar c) \vdash p' \cap
\{\varphi(\bar x,\bar b,\bar c)$:

\hskip25pt $\bar b \subseteq \text{ Rang}(\bar{\bold e}_1 \cup
B_{\bold x})\}$ for each $\varphi \in \Gamma$ where $p'$ is as above.
\end{enumerate}
\mn
There are other variants, we intend to return to this.
\end{discussion}

We now consider some variants of the main Definition \ref{tp14.21}.
\begin{definition}
\label{tp.77}  
1) In Definition \ref{tp14.21}
we add and define $K_\ell,K^\ell_{\lambda,\kappa,\theta}$, etc., also
for $\ell=2,3$ by replacing Clause $(f)_1 = (f)$ by $(f)_\ell$ where:
\mn
\begin{enumerate}
\item[$(f)_2$]   if $\ell=2$ then tp$(\bar c_t,A \cup \{\bar
c_{s,m}:s <_I t,m<n_s\})$ is finitely satisfiable in $B$
\sn
\item[$(f)_3$]    if $\ell=3$ then for some endless indiscernible
sequence $\bar{\bold b}_t = \langle \bar b_{t,r}:r \in J_t\rangle$ of
sequences from $B$, the sequence $\bar c_{t,0}$ 
realizes\footnote{we may consider ``$\bar c_t$ realizes seeming this makes no
difference.} the type
Av$(A \cup \{\bar c_{s,m}:s <_I t,m<n_t\},\bar{\bold b}_t)$.
\end{enumerate}
\mn
2) We define mxK$^\ell_{\lambda,\kappa,\theta}$ similarly.
\end{definition}

\begin{claim}
\label{tp.84}
1) $K_3 \subseteq K_2 \subseteq K_1$.

\noindent
2) If $\ell \in \{2,3\}$ and $\langle \bold x_\alpha:\alpha <
\delta\rangle$ is $\le_1$-increasing in
$K^\ell_{\lambda,\kappa,<\theta}$ and $\delta < \text{\rm cf}(\theta),\delta
< \text{\rm cf}(\kappa)$ \then \, $\bold x_\delta = 
\bigcup\limits_{\alpha < \delta} \bold x_\alpha$ 
defined as in \ref{tp14.28} belongs to
$K^\ell_{\lambda,\kappa,< \theta}$ and is a $\le_1$-lub of the
sequence.

\noindent
3) If $\ell=0,1,2,3$ and $\bar{\bold d} \in {}^{\theta >}{\frak C}$
and $M$ is $\kappa$-saturated, {\rm cf}$(\kappa) \ge \theta$ \then \, 
$\bold x = (M,\emptyset,<>,\bar{\bold d}) \in K^\ell_{\kappa,<\theta}$.

\noindent
4) Like \ref{tp25.33} for $\ell=2,3$, i.e.: if {\rm cf}$(\theta) >
|T|$, {\rm cf}$(\kappa) > \theta,\ell \le 3$ and $\bold x \in
K^\ell_{\lambda,\kappa,< \theta}$ 
\then\footnote{so if $\kappa = \cf(\kappa) > \theta \ge |T|,M$ is
  $\kappa$-saturated and $\bar\bd \in {}^{\theta^+ >}{\frak C}$ then
  for some $\bar\bc \in {}^{\theta^+>}{\frak C}$ and $B \in [M]^{< \kappa}$ we
  have $(M,B,\bar\bc,\bar\bd) \in \mxK_2$.}
 for some $\bold y$ we have
$\bold x \le_1 \bold y \in \text{\rm mxK}^\ell_{\lambda,\kappa,<
\theta}$; so in {\rm tp}$(\bold c_{\bold x},M_{\bold y} + \bold
c_{\bold x})$ we can get {\rm tp}$(\bar c,M)$ is
finitely satisfiable in $M_{\bold y} + (C_{\bold y} \backslash C_{\bold x}))$.

\noindent
5) If $\bold x \in K_2$ and $\bar c \in C_{\bold x}$ 
\then \, {\rm tp}$(\bar c,A_{\bold x})$ is finitely
satisfiable in $B_{\bold x}$.
\end{claim}

\begin{PROOF}{\ref{tp.84}}
Similar to the proofs for $\ell=1$.
\end{PROOF}

\begin{claim}
\label{tp.98}  
1) In \ref{tp25.36} we can deal with
$K^\ell_{\lambda,\kappa,<\theta},\ell=2,3$, i.e.
if $\ell=2$ we should strengthen the assumption to ``$q$ is
finitely satisfiable in $B'$".

\noindent
2) If $\ell=3$ we should strengthen the assumption to $q = \,
\text{\rm Av}(M_{\bold x} \cup C_{\bold x},\bold I),\bold I$ an endless
indiscernible sequence of cardinality $< \kappa$.

\noindent
3) In \ref{tp25.38} we can deal with {\rm mxK}$^2_{\lambda,\kappa,\theta}$.
\end{claim}

\begin{PROOF}{\ref{tp.98}}
Similar to the proof for $\ell=1$.
\end{PROOF}
\newpage

\section {Existence of strict type decompositions} \label{s:existence}

We here get a better decomposition, but at a price: using normal
ultrafilters (so measurable or supercompact cardinals).  Why is the
decomposition from \S2 not enough for our purposes?  See \ref{d10} below.

\begin{hypothesis}
\label{pr.7}  
We assume $T$ is dependent, ${\frak C} = {\frak C}_T$ a monster; if
not said otherwise, we assume (1) or just (2) where:
\mn
\begin{enumerate}
\item[$(1)$]    $(a) \quad \lambda = \kappa$ is a measurable cardinal,
\sn
\item[${{}}$]  $(b) \quad {\cD}$ is a normal ultrafilter on $I = \lambda$, so
$I$ is a linear order,
\sn
\item[${{}}$]  $(c) \quad M_\alpha \prec {\frak C}$ is
$\prec$-increasing, $\|M_\alpha\| < \lambda$ for $\alpha < \lambda$,
\sn
\item[${{}}$]  $(d) \quad M_\lambda = \bigcup\limits_{\alpha < \lambda} 
M_\alpha$, by Clauses (a) + (c) + (e) necessarily $M$ is saturated,
\sn
\item[${{}}$]  $(e) \quad M_\alpha$ is $\|M_\beta\|^+$-saturated for
$\beta < \alpha$,
\sn
\item[$(2)$]    $(a) \quad I$ is the following partial order, which is
$(< \kappa)$-directed, and:
\begin{enumerate}
\item[${{}}$]  $(\alpha) \quad$ set of elements
$\{a \in [\lambda]^{<\kappa}:a \cap \kappa \in \kappa\}$ and
\sn
\item[${{}}$]  $(\beta) \quad s \le_I t$ \Iff \, $s \subseteq t
\wedge |s| < \text{ min}(\kappa \backslash t)$,
\end{enumerate}
\item[${{}}$]  $(b) \quad {\cD}$ is a fine normal 
ultrafilter on $I$ and it follows that $\kappa$ is a measurable 

\hskip35pt  cardinal $\le \lambda$,
\sn
\item[${{}}$]  $(c) \quad M_t \prec {\frak C},\|M_t\| < \kappa$ and
$s <_I t \Rightarrow M_s \prec M_t$, 
\sn
\item[${{}}$]   $(d) \quad M_\lambda = \cup\{M_t:t \in I\}$ (by
  (a)+(c)+(e), $M_\lambda$ is $\kappa$-saturated)
\sn
\item[${{}}$]   $(e) \quad$ if $s <_I t$ then
$M_t$ is $\|M_s\|^+$-saturated.
\end{enumerate}
\end{hypothesis}

\begin{remark}
\label{pr.8}
1) So in \ref{pr.7} we can define:
\mn
\begin{enumerate}
\item[$(A)$]   like (2) without the normality and
\sn
\item[$(B)$]   $(a)(b),(c),(d)$ of part (2).
\end{enumerate}
\mn
2) Note that we have $(1) \Rightarrow (2) \Rightarrow (A)$ 
and $(2) \Rightarrow (B)$.
\end{remark}

\begin{notation}
\label{pr.9}
1) In \ref{pr.7}(1) let
$\kappa_I(t) := t$ for $t \in I$, this notation is introduced only for having 
a uniform treatment of (1) and (2).

\noindent
2) In \ref{pr.7}(2) let $\kappa_I(t) = \Min(\kappa \backslash t)$ 
for $t \in I$.
\end{notation}

\begin{definition}
\label{pr.10}  [under \ref{pr.7}(1) or (2) or alternatively (B) from
\ref{pr.8} so these notions depend on $\langle M_t:t \in I\rangle$.]

\noindent
1) For ${\cU} \subseteq I$ (usually $\in {\cD}$), so is a 
partial order, we say 
$\langle \bar a_t:t \in {\cU}\rangle$ is indiscernible in 
$M_\lambda$ over $A$ \when \, ($A \subseteq M$ and):
\mn
\begin{enumerate}
\item[$(a)$]   $\ell g(\bar a_t)$ is constant, possibly infinite, and 
$\bar a_t \subseteq M$ for $t \in {\cU}$,
\sn
\item[$(b)$]   for each $n$ for some $p_n$ for every $t_0 <_I \ldots
<_I t_{n-1}$ from ${\cU}$ 
we have tp$(\bar a_{t_0} \char 94 \ldots \char 94 \bar
a_{t_{n-1}},A,M) = p_n$, the ninth paragraph of \S0.
\end{enumerate}
\mn
2) We say $\langle \bar a_t:t \in {\cU}\rangle$ is fully
indiscernible (in $M_\lambda$) over $A$ \when \, Clauses (a),(b) above
holds and
\mn
\begin{enumerate}
\item[$(c)$]    if $s <_I t$ are from ${\cU}$ then $\bar a_s
\subseteq M_t$, recalling $\bar M$ is from \ref{pr.7},
\sn
\item[$(d)$]   if $s \in {\cU}$ then \, recalling $\bar M$ is from 
\ref{pr.7}\footnote{By normality 
(i.e. if (1) or (2) or (B) holds) then this
follows.} the sequence $\langle a_t:t \in 
{\cU} \cap I_{\ge s}\rangle$ is indiscernible over $M_s \cup A$ 
 where, of course, $I_{\ge s} := \{t \in I:s \le_I t\}$.
\end{enumerate}
\mn
3) In parts (1),(2) of the definition we say 
$k$-indiscernible when in Clause (b) we demand $n \le k$.
\end{definition}

\noindent
We shall in Theorem \ref{pr.42} (and see \ref{pr.28}) 
below prove the existence of:
\begin{definition}
\label{pr.35}  
1) For an infinite linear order $J$ we say 
$\langle(\bar c_t,\bar d_t):t \in J\rangle$ is a 
strict$_1$ $(\kappa,< \theta)$-decomposition
over $(M,B)$ (and over 
$M$ means ``for some $B \subseteq M$ of cardinality $< \kappa$") \when \,:
\mn
\begin{enumerate}
\item[$(a)$]  $B \subseteq M$ is of cardinality $< \kappa$ and
$M \prec {\frak C}$ is $\kappa$-saturated, but if we
write $< \aleph_0$ instead of $\kappa$ we mean $M \prec {\frak C}$ and
if we write $0$ instead of $\kappa$ we replace $M$ by a set $\supseteq B$,
\sn
\item[$(b)$]  $\alpha = \ell g(\bar c_t),\beta = \ell g(\bar d_t)$ are
$< \theta$,
\sn
\item[$(c)$]  if $t_0 <_J \ldots <_J t_n$ then
tp$(\bar c_{t_n},M + \bar c_{t_0} \char 94 \bar d_{t_0} + \ldots
+ \bar c_{t_{n-1}} \char 94 d_{t_{n-1}})$ 
is increasing with $n \le k$ and does not split over $B$,
\sn
\item[$(d)$]   $\langle(\bar c_t \char 94 \bar d_t:t \in J\rangle$ 
is an indiscernible sequence over $M$,
\sn 
\item[$(e)$]   if $s <_J t$ then tp$(\bar d_t,\bar c_t +
\bar d_s) \vdash \text{ tp}(\bar d_t,\bigcup\{\bar c_r \char 94 \bar
d_r:r \le_J s\} \cup \bigcup\{\bar c_r:r \in J\} \cup M)$,
\sn
\item[$(f)$]   for every $A \subseteq M$ of cardinality $< \kappa$
for some $\bar c \char 94 \bar d \in {}^{\alpha + \beta}M$, the
sequence $\langle \bar c \char 94 \bar d \rangle 
\char 94 \langle \bar c_t \char 94 \bar d_t:t \in J\rangle$ is 
an indiscernible sequence over $A$, so if $\kappa =
0$ this is an empty demand.
\end{enumerate}
\mn
2) We say strict$_{-1} \,
(\kappa,\theta)$-decomposition if (in part (1)) we omit Clauses (e) and (f). 

\noindent
3) We say strict$_0 \, (\kappa,< \theta)$-decomposition 
\If \, we omit (f) and weaken (e) to (e)$^-$, where
\mn
\begin{enumerate}
\item[$(e)^-$]  if $s <_I t$ then tp$(\bar d_t,\bar c_t + \bar d_s)
\vdash \text{ tp}(\bar d_t,\cup\{\bar c_r \char 94 \bar d_r:r \le s\} \cup
\bar c_t \cup M)$.
\end{enumerate}
\mn
4) Notation:
\mn
\begin{enumerate}
\item[$\bullet$]  If $\theta = \sigma^+$ instead of ``$< \theta$" we may write
$\sigma$.
\sn
\item[$\bullet$]    If $\kappa = 0$ 
then $M$ is replaced by a set $B$, if we write $< \aleph_0$ instead of
$\kappa$ then $M$ is just a model.
\sn
\item[$\bullet$]   Strict$_1$ may be written strict.
\end{enumerate}
\end{definition}

\noindent
A natural question about those notions of indiscernibility is about
existence results.  Now \ref{pr.21} is a well known set-theoretic existence and
\ref{pr.28} is existence for dependent theories.
\begin{fact}
\label{pr.21} 
1) If $A \subseteq \gC,|A| < \kappa,\alpha < \kappa,{\cU}_1 \in 
{\cD}$ and $\bar a_t \in {}^\alpha \gC$ for $t \in {\cU}_1$ \then \, 
for some ${\cU}_2 \subseteq {\cU}_1$ from ${\cD}$ the sequence 
$\langle \bar a_t:t \in {\cU}_2\rangle$ is indiscernible over $A$.

\noindent
2) If in addition $\bar a_t \in {}^\alpha(M_\lambda)$ \then \, we can add
``fully indiscernible".

\noindent
3) If $\iota \in \{-1,0,1\}$ and $\langle (\bar c_t,\bar d_t):t \in
   J\rangle$ is a strict$_\iota \, (\kappa,<\theta)$-decomposition over
$(M,B)$ and $M \subseteq B_1 \subseteq B$ \then \, it is a
strict$_\iota \, (\kappa,< \theta)$-decomposition over $(M,B)$.
Similarly we can replace $M$ by $M'$ if $B_1 \subseteq M'
   \subseteq M$ and $M_1$ satisfies Clause (a) in Definition \ref{pr.35}.
\end{fact}

\begin{PROOF}{\ref{pr.21}}
1),2) By well known set theory (see Kanamori Magidor \cite{KnMg78}).
\end{PROOF}

\begin{observation}
\label{pr.23}
1) For some $\cU_* \in \cD$, for every $t \in I$ the model $M_t$ is 
$\kappa(t)$-saturated and $\kappa(t) > |T|$.

\noindent
2) If $M \prec \gC$ is $\kappa$-saturated, $J$ an infinite linear
   order $\theta \ge \aleph_0,B \in [M]^{< \kappa}$ and $\bar c_t =
   \langle \rangle = \bar d_t$ for $t \in J$ \then \, $\langle (\bar
   c_t:\bar d_t):t \in J \rangle$ is a strict$_\iota$ $(\kappa,<
   \theta)$-decomposition over $(M,B)$ for $\iota =1$ hence for $\iota
   \in \{-1,0\}$, too.
\end{observation}

\begin{PROOF}{\ref{pr.23}}
Obvious.
\end{PROOF}

\begin{theorem}
\label{pr.28}  
Assume $\theta$ satisfies $\kappa > \theta \ge |T|$ and recall
$\lambda \ge \kappa$.  
For every $\gamma(*) < \theta^+$ and
$\bar d \in {}^{\gamma(*)}{\gC}$ there are $B$ and
$\langle (\bar c^\omega_t,\bar d^\omega_t):t \in {\cU}
\cup \{\lambda\}\rangle$ such that:
\mn
\begin{enumerate}
\item[$\boxtimes$]  $(a) \quad {\cU} \in {\cD}$,
\sn 
\item[${{}}$]  $(b) \quad B \subseteq M_\lambda,|B| < \kappa$,
\sn
\item[${{}}$]   $(c) \quad \ell g(\bar d^\omega_t) = \gamma(*) 
+ \theta \cdot \omega$,
\sn
\item[${{}}$]  $(d) \quad \bar d \trianglelefteq \bar
d^\omega_\lambda$,
\sn
\item[${{}}$]   $(e) \quad \bold x = (M_\lambda,B,\bar c^\omega_\lambda,
\bar d^\omega_\lambda) \in \text{\rm mxK}_{\kappa,\theta}$,
\sn 
\item[${{}}$]   $(f) \quad \bar c^\omega_t \char 94 \bar d^\omega_t
\subseteq M_\lambda$ realizes {\rm tp}$(\bar c^\omega_\lambda \char 94 \bar
d^\omega_\lambda,M_t)$ for $t \in {\cU}$, 
\sn 
\item[${{}}$]  $(g) \quad \langle \bar c^\omega_t \char 94 \bar
d^\omega_t:t \in {\cU}\rangle$ is fully indiscernible (in $M_\lambda$)
over $B$ and even over

\hskip25pt  $B \cup \bar c^\omega_\lambda \cup 
\bar d^\omega_\lambda$, (see Definition \ref{pr.10}(2)),
\sn
\item[${{}}$]  $(h)_1 \quad$ if $t_0 <_I \ldots <_I t_m <_I \ldots
<_I t_n$ belongs to ${\cU}$, so $m<n$ and possibly

\hskip25pt  $m+1=n$; moreover possibly $0=m$  \then \,

\hskip25pt  $\tp(\bar d^\omega_{t_{m+1}},\bar c^\omega_{t_{m+1}} +$
$\bar c^\omega_{t_{m+2}} + \ldots + \bar c^\omega_{t_n} + \bar d^\omega_{t_n})
\vdash$

\hskip25pt $\tp(\bar d^\omega_{t_{m+1}},\bar c^\omega_{t_0} 
+ \ldots + \bar c^\omega_{t_{m+1}} + \ldots +
\bar c^\omega_{t_n} + \bar d^\omega_{t_0} + \ldots + \bar
d^\omega_{t_m} + M_{t_0})$,
\sn 
\item[${{}}$]   $(h)_2 \quad$ if $s <_I t$ are from ${\cU}$ 
\then \, $\bar c^\omega_s \char 94 \bar d^\omega_s$ is from $M_t$,
(actually follows 

\hskip25pt from clause (g)),
\sn
\item[${{}}$]   $(h)_3 \quad$ if $t_0 <_I \ldots <_I
t_n$ are from ${\cU}$ \then \, {\rm tp}$(\bar d^\omega_{t_1},
\bar d^\omega_{t_0} + \bar c^\omega_{t_1}) \vdash \,\text{\rm tp}
(\bar d^\omega_{t_1},M_{t_0} +$

\hskip25pt $\bar d^\omega_{t_0} + \bar c^\omega_{t_0} +
\bar c^\omega_{t_1} + \ldots + \bar c^\omega_{t_n})$ (actually this is
the case $m=0$ in $(h)_1$).
\end{enumerate}
\end{theorem}

\begin{remark}
\label{pr.29}
We easily can add:
\mn
\begin{enumerate}
\item[${{}}$]  $(i) \quad \bold x$ is normal, i.e. 
{\rm Rang}$(\bar c_\omega) \subseteq \, \text{\rm Rang}(\bar d_\omega)$.
\end{enumerate}
\end{remark}

\begin{PROOF}{\ref{pr.28}}

First by induction on $n$ we choose $\bar d_n,\bar c_n,
B_n,\langle (\bar c^n_t,\bar d^n_t):t \in {\cU}_n\rangle$ and 
if $n>0$ also $\bar e_n,\bar e^n_t$ (for $t \in {\cU}_n$) such that:
\mn
\begin{enumerate}
\item[$\circledast_n$]  $(a) \quad \bar d_n \in 
{}^{\gamma(*)+\theta \cdot n}{\frak C}$ and $\bar d_0 = \bar d$ and
$\bar e_n \in {}^\theta \gC$
\sn 
\item[${{}}$]  $(b) \quad \bold x_n = (M_\lambda,B_n,\bar c_n,\bar d_n) \in 
\text{\rm mxK}_{\kappa,\theta}$ is normal,
\sn
\item[${{}}$]  $(c) \quad \bar d_m \triangleleft \bar d_n$ if $m <n$,
\sn
\item[${{}}$]  $(d) \quad \bar c_m = \bar c_n \restriction I_{\bold x_m}$
 if $m<n$,
\sn
\item[${{}}$]  $(e) \quad \bar c^n_t \char 94 \bar d^n_t$ is from
$M_\lambda$ and realizes tp$(\bar c_n \char 94 \bar d_n,M_t)$ 
for $t \in {\cU}_n$,
\sn
\item[${{}}$]  $(f) \quad \bar d_n = \bar d_m \char 94 \bar e_n$
and $\ell g(\bar e_n) = \theta$ and $\bold x_m \le_1 \bold x_n$ if $n=m+1$,
\sn
\item[${{}}$]  $(g) \quad \langle \bar c^n_t \char 94 \bar d^n_t:t
\in {\cU}_n\rangle$ is fully indiscernible over $B_n + \bar c_n + \bar d_n$,
\sn
\item[${{}}$]  $(h) \quad {\cU}_n \in {\cD}$ decrease with $n$ and is
$\subseteq \cU_*$ from \ref{pr.23}(1),
\sn
\item[${{}}$]  $(i) \quad$ if $s<t$ are from ${\cU}_n$ then $M_s
+ \bar c^n_s + \bar d^n_s \subseteq M_t$ (follows by (g)),
\sn
\item[${{}}$]  $(j) \quad$ {\rm tp}$(\bar d_m,\bar e^n_t + \bar c_m) \vdash
\text{ tp}(\bar d_m,M_t + \bar c_m + \bar c^m_t + \bar d^m_t)$ if $n =
m+1,t \in {\cU}_n$,
\sn
\item[${{}}$]  $(k) \quad (M_t,B_n,\bar c_n,\bar d_n) \in
\text{\rm mxK}_{\kappa_I(t),\theta}$ for $t \in {\cU}_n$,
\sn
\item[${{}}$]   $(l) \quad$ if $n = m+1,k < \omega$ and $t_0 <_I
\ldots <_I t_k$ are from ${\cU}_n$ then 

\hskip25pt tp$(\bar d^m_{t_1},\bar c^m_{t_1} + \bar e^n_{t_0}) \vdash
\text{ tp}(\bar d^m_{t_1},\sum\limits_{\ell=0}^{k} \bar c^m_{t_\ell} +
\bar d^m_{t_0} + \bar e^n_{t_0} + M_{t_0})$.
\end{enumerate}
\bigskip

\noindent
\underline{Case 1}:  $n=0$.

First, let $\bar d_0 = \bar d$.  Recalling \ref{tp14.28}(4) 
clearly ${\bold y}_n =: (M_\lambda,\emptyset,<>,\bar d_0) 
\in K_{\kappa,\theta}$.

Second, by Claim \ref{tp25.33} we can find $B_n,\bar c_n$ 
such that $\bold y_n \le_1 {\bold x}_n := (M_\lambda,B_n,
\bar c_n,\bar d_n) \in \mxK_{\kappa,\theta}$.

Third, for $t \in I$ we can choose $\bar c^n_t \char 94 \bar d^n_t$ from
$M_\lambda$ which realizes tp$(\bar c_n \char 94 \bar d_n,M_t)$.

Fourth, by \ref{pr.21} we choose ${\cU}_0 \in {\cD}$ 
such that $\langle \bar c^n_t \char 94 
\bar d^n_t:t \in {\cU}_0\rangle$ is a fully indiscernible sequence
over $B_0 + \bar d_0 + \bar c_0$ 
and (by the normality of the filter ${\cD}$) in particular
$\bar c^n_s \char 94 \bar d^n_s \subseteq M_t$ when $s<t \in 
{\cU}_0$ are from $I$ and $\cU_0$ is $\subseteq \cU_*$ from \ref{pr.21}.
It is easy to check that all the demands hold, recalling $\bar e_n$
for $n=0$ is not required in $\circledast_n$.
\bigskip

\noindent
\underline{Case 2}:  $n = m+1$.

First, by clause (k) for $m$ and Conclusion \ref{tp14.22} for each 
$t \in {\cU}_m$ recalling $M_t + \bar c^m_t + \bar d^m_t$ is
$\subseteq M_\lambda$ and of cardinality $< \kappa$ 
and $M_\lambda$ is $\kappa$-saturated, 
there is $\bar e^{n,*}_t \in {}^\theta(M_\lambda)$ such that: 
\mn
\begin{enumerate}
\item[$(*)_1$]  $\tp(\bar d_m,\bar c_m + \bar e^{n,*}_t)
\vdash \tp(\bar d_m,M_t + \bar c_m + \bar c^m_t + \bar d^m_t)$.
\end{enumerate}
\mn
Second, by \ref{pr.21} choose ${\cU}'_n \subseteq {\cU}_m$ which belongs to
${\cD}$ such that $\langle \bar c^m_t \char 94 \bar d^m_t \char 94 \bar
e^{n,*}_t:t \in {\cU}'_n\rangle$ is fully indiscernible over $B_m 
+ \bar c_m + \bar d_m$. 

We shall now prove
\mn
\begin{enumerate}
\item[$\odot_1$]   if $t_0 < \ldots < t_k$ are from ${\cU}'_n$
\then \, tp$(\bar d^m_{t_1},\bar c^m_{t_1} + \bar e^{n,*}_{t_0}) \vdash
\text{ tp}(\bar d^m_{t_1},\bar c^m_{t_1} + \ldots + \bar c^m_{t_k} + M_{t_0})$.
\end{enumerate}
\mn
Toward this, by clause (e) of $\circledast_m$ we have
\mn
\begin{enumerate}
\item[$(*)_2$]  tp$(\bar c^m_{t_1} \char 94 \bar d^m_{t_1},M_{t_1})
= \text{ tp}(\bar c_m \char 94 \bar d_m,M_{t_1})$.
\end{enumerate}
\mn
By $(*)_2$ there is an elementary mapping $f$ mapping $\bar c_m\char
94 \bar d_m$ to $\bar c^m_{t_1} \char 94 \bar d^m_{t_1}$ which is the identity
on $M_{t_1}$.  But $M_{t_0} + \bar c^m_{t_0} + \bar d^m_{t_0} \subseteq
M_{t_1}$ by $\circledast_m(i)$ and $\bar e^{n,*}_{t_0} \subseteq
M_{t_1}$ by the full indiscernibility, i.e. by the choice of $\cU'_n$ 
above, hence by applying $f$ on $(*)_1$ for $t_1$ and monotonicity we get
\mn
\begin{enumerate}
\item[$(*)_3$]  tp$(\bar d^m_{t_1},\bar c^m_{t_1} + 
\bar e^{n,*}_{t_0}) \vdash \text{ tp}(\bar d^m_{t_1},M_{t_0} 
+ \bar c^m_{t_1} + \bar c^m_{t_0} + \bar d^m_{t_0})$.
\end{enumerate}
\mn
Now by clause (k) of $\circledast_m$
\mn
\begin{enumerate}
\item[$(*)_4$]   $(M_{t_0},B_m,\bar c_m,\bar d_m) \in
\text{\rm mxK}_{\kappa_I(t_0),\theta}$.
\end{enumerate}
\mn
But by clause (e) of $\circledast_m$
\mn
\begin{enumerate}
\item[$(*)_5$]   tp$(\bar c^m_{t_1} \char 94 \bar d^m_{t_1},M_{t_0})
= \text{ tp}(\bar c_m \char 94 \bar d_m,M_{t_0})$.
\end{enumerate}
\mn
By $(*)_4 + (*)_5$
\mn
\begin{enumerate}
\item[$(*)_6$]   $(M_{t_0},B_m,\bar c^m_{t_1},\bar d^m_{t_1}) \in
\text{\rm mxK}_{\kappa_I(t_0),\theta}$.
\end{enumerate}
\mn
Also easily by clauses (k) + (i) of $\circledast_m$ applied to $t =
t_2,t_3,\dotsc,t_k$ (recalling \ref{3k.4}(3))
\mn
\begin{enumerate}
\item[$(*)_7$]   tp$(\bar c^m_{t_2} \char 94 \ldots \char 94 
\bar c^m_{t_k},M_{t_0} + \bar c^m_{t_1} + \bar d^m_{t_1})$ does not split
over $B_m$.
\end{enumerate}
\mn
By $(*)_6 + (*)_7$ and the weak orthogonality claim \ref{tp25.36}(1) we have
\mn
\begin{enumerate}
\item[$(*)_8$]   tp$(\bar d^m_{t_1},M_{t_0} + \bar c^m_{t_1}) \vdash
\text{ tp}(\bar d^m_{t_1},M_{t_0} + \bar c^m_{t_1} + \ldots + \bar
c^m_{t_k})$.
\end{enumerate}
\mn
By $(*)_3 + (*)_8$
\mn
\begin{enumerate}
\item[$(*)_9$]   tp$(\bar d^m_{t_1},\bar c^m_{t_1} + \bar
e^{n,*}_{t_0}) \vdash \text{ tp}(d^m_{t_1},M_{t_0} + \bar c^m_{t_1} +
\ldots + \bar c^m_{t_k})$
\end{enumerate}
\mn
as promised in $\odot_1$.
\medskip

Now we continue to deal with Case 2, 
choose $F_n:{\cU}'_n \rightarrow {\cU}'_n$ such that $s \in
{\cU}'_n \Rightarrow s <_I F_n(s) \in {\cU}'_n$ and for $t \in
{\cU}'_n$ we let $\bar e^n_t := \bar e^{n,*}_{F_n(t)}$.  Let $\bar
e_n \in {}^\theta \gC$ 
be such that $\bar c_m \char 94 \bar d_m \char 94 \bar e_n$
realizes Av$(\langle \bar c^m_t \char 94 \bar d^m_t \char 94 \bar
e^n_t:t \in {\cU}'_n\rangle/{\cD},M_\lambda)$.
Let ${\cU}''_n \in {\cD}$ be $\subseteq 
{\cU}'_n$ and such that $s \in {\cU}''_n \wedge 
t \in {\cU}''_n \wedge s <_I t \Rightarrow F_n(s) <_I t$ and the sequence
$\langle \bar c^m_t \char 94 \bar d^m_t \char 94 \bar e^{n,*}_{F(t)}:t \in
{\cU}''_n\rangle$ is fully indiscernible over $\bar c_m \char 94
\bar d_m \char 94 \bar e_n$.

Let $\bar d_n = \bar d_m \char 94 \bar e_n$ and
$\bar d^n_t = \bar d^m_t \char 94 \bar e^n_t$ for $t
\in {\cU}''_n$.  

Let ${\bold y}_n := (M,B_m,\bar c_m,\bar d_n) \in
K_{\kappa,\theta}$ so clearly $\bold x_m \le_1 \bold y_n$ 
hence by \ref{tp25.33} there is ${\bold x}_n = (M_\lambda,B_n,\bar c_n,\bar
d_n) \in \text{ mxK}_{\kappa,\theta}$ such that ${\bold y}_n \le_2 
{\bold x}_n$ so $\bar c_n$ and $B_n$ are well 
defined\footnote{in fact, we can demand that $\tp(\bar c_n \restriction
  (I_{\bold x_n} \backslash I_{\bold x_m}),M + \bar c_m)$ does not
  split over $B_{\bold x_n}$.}
 and $\bar c_m = \bar c_n \restriction I_{\bold x_m}$.  
For $t \in \cU''_n$, let $\bar c^n_t$ be a sequence from $M_\lambda$ such
that $\bar c^m_t = \bar c^n_t \restriction \text{ Dom}(\bar c_m)$ and
tp$(\bar c^n_t \char 94 \bar d^n_t,M_t) = \text{ tp}(\bar c_n \char 94
\bar d_n,M_t)$, this is possible as $M$ is $\kappa$-saturated, $|M_t +
\bar c^m_t + \bar d^n_t| < \kappa,|\text{Rang}(\bar c_n)| < \kappa$ 
and $\bar c^m_t \char 94 \bar d^n_t$
realizes the type tp$(\bar c_m \char 94 \bar d_n,M_t)$.

Lastly, let ${\cU}_n$ be a subset of ${\cU}''_n$ which belongs to
${\cD}$ such that:
\mn
\begin{enumerate}
\item[$\odot_2$]   $\langle \bar c^n_t \char 94 \bar d^n_t:t \in 
{\cU}_n\rangle$ is fully indiscernible over $B_n + \bar c_n + \bar d_n$.
\sn
\item[$\odot_3$]    $(M_t,B_n,\bar c_n,\bar d_n) \in \text{\rm
mxK}_{\kappa,\theta}$ for every $t \in {\cU}_n$.
\end{enumerate}
\mn
[Why ${\cU}_n$ exists?  By \ref{pr.21} and \ref{tp25.30}.]

It is easy to check that $\bold x_n,\bar c_n,\bar d_n,\bar
d_n,\langle(c^n_t,\bar d^n_t,\bar e^n_t):t \in {\cU}_n\rangle$ are
as required.  E.g. clause (f) holds as
$\ell g(\bar d_0) = \ell g(\bar d) = \gamma(*)$ and $\ell g(\bar d_m) =
\gamma(*) + \theta \cdot m$ by clause (a) of $\circledast_m$
 and $\ell g(\bar e_n) = \theta$ by $\circledast_n(f)$ 
we can prove that $\ell g(\bar d_n) = 
\ell g(\bar d_m) + \theta = \gamma(*) + \theta \cdot m + \theta = 
\gamma(*) + \theta n$, so we clearly are done.  For clause(e) note
$\odot_1$ and the choices of $F$ and $\cU''_n,\cU_n$.

So we have carried the induction.  Second, let $\bar c_\omega = \bar
c^\omega_\lambda = \bigcup\limits_{n < \omega} \bar c_n,\bar d_\omega 
= \bar d^\omega_t = \bigcup\limits_{n < \omega} \bar d^n_t$ and 
$B = \cup\{B_n:n < \omega\}$ and ${\cU} = \cap\{{\cU}_n:n < \omega\}$.

\noindent
Let us check that $\boxtimes$ from the theorem holds indeed.
\bigskip

\noindent
\underline{Clause (a)}:  ${\cU} \in {\cD}$ as each ${\cU}_n \in {\cD}$ by
$\circledast_n(h)$ and ${\cD}$ is $\kappa$-complete and $\kappa >
\aleph_0$ recalling $\cU = \cap\{\cU_n:n <\omega\}$.
\bigskip

\noindent
\underline{Clause (b)}:  $B \in [M_\lambda]^{< \kappa}$ 
as $B_n \subseteq M_\lambda,|B_n| <
\kappa$ by $\circledast_n(b)$ 
for $n < \omega$ and $\kappa$ is regular uncountable recalling
and $B := \cup\{B_n:n < \omega\}$.
\bigskip

\noindent
\underline{Clause (c)}:  By $\circledast_n(a) + (c)$ for $n < \omega$.
\bigskip

\noindent
\underline{Clause (d)}:  $\bar d = \bar d_0 = \bar d_\omega \restriction
\gamma(*)$ is proved as in clause (c).
\bigskip

\noindent
\underline{Clause (e)}:  As ${\bold x}_n = (M,B_n,\bar c_n,\bar d_n) \in 
\text{\rm mxK}_{\kappa,\theta}$ by clause (b) of $\circledast_n$ 
and ${\bold x}_n \le_1 {\bold x}_{n+1}$ by $\circledast_{n+1}(f)$ 
and $(M,B,\bar c_\omega,
\bar d_\omega) = \cup\{{\bold x}_n:n < \omega\}$, clearly by claim
\ref{tp25.36}(2) we are done.
\bigskip

\noindent
\underline{Clause (f)}: By clause (e) of $\circledast_n$ (and the
choice of $\bar c_\omega,\bar d_\omega,\bar c^\omega_t,\bar d^\omega_t$, etc). 
\bigskip

\noindent
\underline{Clause (g)}:  Similarly by clause (g) of $\circledast_n$.
\bigskip

\noindent
\underline{Clause (h)$_1$}:  By clauses
$(h)_2 + (h)_3$ proved below. 
\bigskip

\noindent
\underline{Clause (h)$_2$}:   By clause (i) of $\circledast_n$.
\bigskip

\noindent
\underline{Clause (h)$_3$}:  Holds by clause $(\ell)$ of $\circledast_n$.
\end{PROOF}

\begin{theorem}
\label{pr.42}
1) If $M$ is $\kappa$-saturated of cardinality $\le \lambda,
\bar d \in {}^{\theta^+ >}{\frak C}$ 
\then \, we can find a strict $(\kappa,\theta)$-decomposition $\langle(\bar
c_n,\bar d_n):n <\omega\rangle$ over $M$ such that $\bar d \triangleleft
\bar d_0$.

\noindent
2) Instead Hypothesis \ref{pr.7} it is enough to demand: if $M$ is
 $\kappa$-saturated and $p = p(\langle x_i:i < \theta\rangle)$ 
is a type with parameters from $M \cup C,|C| \le \theta$ 
which is $(< \kappa)$-satisfiable in $M$, i.e.
every subset of $p$ of cardinality $< \kappa$ is realized in $M$ \then
 \, $p$ can be extended to $p^+ \in
 \bold S^\theta(M \cup C)$ which is $(<\kappa)$-satisfiable in $M$.
\end{theorem}

\begin{PROOF}{\ref{pr.42}}
1) We can choose $\bar M' = \langle M'_t:t \in I\rangle$ such that
$(M,\bar M')$ satisfies the demands on $(M_\lambda,\langle M_t:t \in
I\rangle)$ in Hypothesis \ref{pr.7}(1) or \ref{pr.7}(2) and apply
Theorem \ref{pr.28} (as assuming $\ell g(\bar d)= \theta$ or $\ell g(\theta) <
 \theta^+$ does not matter).

\noindent
2) The idea is to repeat the proof of \ref{pr.28}, but as of unclear
 value we leave it to the reader.
\end{PROOF}

\begin{corollary}
\label{pr.56}  
Assume $\kappa = \lambda > |T|$ is
weakly compact, $M_\alpha \in \text{\rm EC}_{< \lambda}(T)$ is
$\prec$-increasing continuous, $M = \cup\{M_\alpha:\alpha < \lambda\}$
is saturated.  \Then \, \ref{pr.28} and \ref{pr.49}(2) hold.
\end{corollary}

\begin{PROOF}{\ref{pr.56}}
Revise the proof of \ref{pr.28}, but in $\circledast_n$
weaken clauses (g),(h) to (g)$^-$,(h)$^-$ and use the proof of \ref{ps.21} in
the end where
\mn
\begin{enumerate}
\item[$(g)^-$]   ${\cU}_n \in [\kappa]^\kappa$ decreasing with
$n$ and $\langle \bar c^n_t \char 94 \bar d^n_t:t \in 
{\cU}_n\rangle$ is just fully $n$-indiscernible
\sn
\item[$(h)^-$]  $\cU_n$ does not belong to the weakly compact ideal.
\end{enumerate}
\mn
We leave the details to the reader.
\end{PROOF}

\begin{claim}
\label{pr.49}  
1) (even not assuming \ref{pr.7})

Assume $M$ is $\kappa$-saturated, $|T| \le \theta < \kappa$, {\rm
tp}$(\bar c_n,M + \bar c_0 + \ldots + \bar c_{n-1})$ does not split
over $B$ where $|B|<\kappa,B \subseteq M$ and $\ell g(\bar c_n) <
\theta^+$ and $\bar d_n = \bar c_n \restriction u$ for $n
<\omega$ (so $u \subseteq \text{\rm Dom}(\bar c_n))$ \then \,
$\langle(\bar c_n,\bar d_n):n <\omega\rangle$ is a strict
$(\kappa,\theta)$-decomposition over $(M,B)$.

\noindent
2) (Assuming \ref{pr.7}!)  Assume $\langle(\bar c_n,\bar d_n):n <
\omega\rangle$ is a strict $(\kappa,\theta)$-decomposition over
$(M_\lambda,B)$.  For any $\bar d \in {}^{\theta^+ >}{\frak C}$ we can find a
strict $(\kappa,\theta)$-decomposition $\langle(\bar c^+_n,\bar
d^+_n):n < \omega\rangle$ over $(M_\lambda,B)$ such that 
$\bar c_n \trianglelefteq
\bar c^+_n,\bar d_n \trianglelefteq \bar d^+_n$ for $n < \omega$ 
and $\bar d_0 \char 94 \bar d \trianglelefteq \bar d^+_0$.
\end{claim}

\begin{PROOF}{\ref{pr.49}}
1) Easy.

\noindent
2) Repeat the proof of \ref{pr.28} starting with $\bold x_0 =
(M_\lambda,B,\bar c_0,\bar d_0)$ and $\langle \bar c_t \char 94 \bar d_t:t \in
I\rangle$ such that $\bar c_t \char 94 d_t \in {}^{\ell g(\bar c_t) +
  \ell g(\bar d_t)}M$ realizing tp$(\bar c_0 \char 94 \bar d_0,M_t +
\sum\limits_{n < \omega} \bar c_{1+n} \char 94 \bar d_{1 +n})$.   
\end{PROOF}

\begin{claim}
\label{pr.84}   
The sequence $\langle (\bar c^2_t,\bar d^2_t):t
\in I_2\rangle$ is a strict$_0 \, (0,\theta)$-decomposition over
$(B_2,B_1)$ \when \,:
\mn
\begin{enumerate}
\item[$(a)$]   $\langle(\bar c^1_s,\bar d^1_s):s \in I_1)\rangle$ is
a strict$_0$\, $(0,\theta)$-decomposition (and $\ell g(\bar c^1_s) = \ell
g(\bar c^2_t),\ell g(\bar d^1_s) = \ell g(\bar d^1_s) = \ell g(\bar
d^2_t)$, of course),
\sn
\item[$(b)$]   for any $n$ if $I_\ell \models ``t_{\ell,0} < \ldots <
t_{\ell,n-1}"$ for $\ell =1,2$ then $\tp(\bar c^1_{t_{1,0}} \char 94
\bar d^1_{t_{1,0}} \char 94 \ldots \char 94$

\hskip15pt $\bar c^1_{t_{1,n-1}} \char 94 
\bar d_{t_{1,n-1}},B_2) = \tp(\bar c^2_{t_{2,0}} \char 94 \bar
d^2_{t_{2,0}} \char 94 \ldots \char 94 \bar c^2_{t_{2,n-1}} \char 94
\bar d^2_{t_{2,n-1}},B_2)$.
\end{enumerate}
\end{claim}

\begin{PROOF}{\ref{pr.84}}
Should be clear.
\end{PROOF}
\newpage

\section {Consequences of strict decomposition} \label{s:consequences}

Here we look again at the generic pair conjecture
(\cite[0.2]{Sh:877}).  
The non-structure side (in a strong version) is
proved there for $\lambda = \lambda^{< \lambda}$ non-strong
limit and in \cite{Sh:906} for $\lambda = \lambda^{<\lambda}$ strong
limit (i.e. strongly inaccessible).  

\noindent
The conjecture is (the instances of G.C.H. are used to make the
conjecture transparent):
\begin{conjecture}
\label{ps.1}
\underline{The generic pair conjecture}  

Assume\footnote{the 
``$2^\lambda = \lambda^+$" is just for making the formulation more
transparent, and by absoluteness is equivalent to the formulation not
assuming $2^\lambda = \lambda^+$.}
 $\lambda = \lambda^{<\lambda} > |T|,2^\lambda = \lambda^+,M_\alpha 
\in \text{ EC}_\lambda(T)$ is $\prec$-increasing continuous for $\alpha <
\lambda^+$ with $\cup\{M_\alpha:\alpha < \lambda^+\} \in 
\text{ EC}_{\lambda^+}(T)$ being saturated.  

\noindent
1) The $\lambda$-generic pair conjecture says that: $T$ is 
dependent \Iff \, for some club $E$ of $\lambda^+$ for all 
pairs $\alpha < \beta < \lambda^+$ from $E$
of cofinality $\lambda^+,(M_\beta,M_\alpha)$ has the same isomorphism type.

\noindent
2) For $\varepsilon < \lambda$ the $\lambda$-generic $\zeta$-tuple
conjecture says that: $T$ is dependent \Iff \, for some club $E$ of
   $\lambda^+$ for all increasing sequences $\langle
   \alpha_\varepsilon:\varepsilon \le \zeta\rangle$ of members of $E$
of cofinality $\lambda$, the structure $(M_{\alpha_\zeta},
M_{\alpha_\varepsilon})_{\varepsilon < \zeta}$ has the same isomorphism
type (equivalently, if $\langle
   \alpha_{\ell,\varepsilon}:\varepsilon \le \zeta\rangle$ is as above
   for $\ell=1,2$ then there is an isomorphism from
   $M_{\alpha_{1,\zeta}}$ onto $M_{\alpha_{2,\zeta}}$ mapping
   $M_{\alpha_{1,\varepsilon}}$ onto $M_{2,\alpha_\varepsilon}$ for
   $\varepsilon < \zeta$).
\smallskip

We concentrate on the pair.  Note that if $\kappa = \cf(\kappa) <
\lambda$, then the $\lambda$-generic $\kappa$-tuple conjecture implies
that for dependent $T$ there is a medium $(\lambda,\kappa)$-limit
model, see \cite{Sh:877}, but we do not succeed to deal with it here.

Here we prove the ``structure" side when $\lambda$ is measurable.
It seemed natural to assume that the first order theories of
such pair is complicated if $T$ is independent and 
``understandable" for dependent of $T$. 
\end{conjecture}

In fact, it may be better to ask
\begin{problem}
\label{ps.3}  
1) Assume $|T| < \theta \le \kappa \le \lambda =
\lambda^{< \kappa} < \kappa_2 \le \mu = \mu^{< \aleph_2}$ and $M_1
\prec M_2 \prec {\frak C},M_1$ is $\kappa_1$-saturated of cardinality
$\lambda,M_2$ is $\kappa_2$-saturated of cardinality $\mu$.  What can we say
on Th$(M_2,M_1)$?  On Th$_{\bbL_{\infty,\theta}(\tau(T))}(M_2,M_1)$?
\end{problem}

More generally
\begin{problem}
\label{gs.3}
1) Assume $n < \omega,|T| < \theta,\sigma < \theta 
< \kappa_0,\lambda_\ell = \lambda^{< \kappa_\ell}_\ell$ for
$\ell \le n,\lambda_\ell < \kappa_{\ell +1}$ for $\ell < n$.
Let $M_\ell$ be $\kappa_\ell$-saturated of cardinality $\lambda_\ell$
for $\ell \le n$ and $M_\ell \prec M_{\ell +1}$ for $\ell < n$.  What
can we say on $M^+ = \text{ Th}(M_n,\ldots,M_1,M_0)$,i.e. $M_n$ expanded by
unary predicates for $M_\ell$ for $\ell < n$?  When can we interpret
(with first order formulas with parameters) second order logic on
$\theta$?  i.e. classify $T$ by this.

\noindent
2) Similarly for $\bbL_{\sigma,\sigma}(\tau_{M^+})$. 

\noindent
3) Similarly allowing $n$ to be $< \theta$.

The proof here, if e.g. $\kappa = \lambda$ is measurable say that even
the $\bbL_{\infty,\kappa}$-theory of the pair is constant, 
but does not say much
even on the first order theory, (see \cite{KpSh:946}).  
It is known that for many ``complicated enough" theories $T$, for
$M_2,M_1$ as in \ref{ps.3}, in $\Th(M_2,M_1)$ we can interpret second 
order logic on $(\lambda,=)$.  This holds, e.g. for $T =$ Peano arithmetic.
\end{problem}

On $n$-independent theories see \cite[\S2]{Sh:886}.  Note that
\begin{claim}
\label{gs.4}  Assume $T$ is the model 
completion of $T_0$, defined below so seems 
``the simplest" $2$-independent theory.  If $M_0 \in
\text{\rm EC}_\lambda(T)$ and $M_1$ is a $\lambda^+$-saturated
$\prec$-extension of $M_0$ \then \, in $(M_1,M_0)$ we can interpret
second order logic on $M_0$ (i.e. quantification on two-place
relations, \when \,:
\mn
\begin{enumerate}
\item[$(*)$]   $\tau_{T_0}$ consists of $P_0,P_1,P_2$ (unary
  predicates) and $R$ (a ternary predicate) and a $\tau_{T_0}$-model
  $M$ is a model of $T_0$ iff $\langle P^M_0,P^M_1,P^M_2\rangle$ is a
  partition of $|M|$ and $R^M \subseteq P^M_0 \times P^M_1 \times P^M_2$.
\end{enumerate}
\end{claim}

\begin{PROOF}{\ref{gs.4}}
Obvious.  
\end{PROOF}
\bigskip

\noindent
Also for $T =$ theory of Boolean algebras (which is $n$-independent
for every $n$) the theory is complicated.
Of course, it would be better to eliminate the measurable assumption.
\begin{explanation}
\label{d10}
Why $K_{\lambda,\kappa,\theta}$ and $\mxK_{\lambda,\kappa,\theta}$
from \S2 does not suffice for us so that in \S3 we deal also with the
more complicated $\langle (\bar c_t,\bar d_t):t \in J\rangle$ from
\ref{pr.35}?  This is motivated by the proof of the generic pair
conjecture.

To understand it maybe better consider the class
\mn
\begin{enumerate}
\item[$(*)$]  $\bold N^2_\kappa = \{(N,M):M \prec N \prec \gC_T,M$ is
  $\kappa$-saturated and $N$ is $\|M\|^+$-saturated$\}$.
\end{enumerate}
\mn
Proving the generic pair conjecture for $\kappa$ we consider $\bar M =
\langle M_\alpha:\alpha < \kappa^+\rangle$, which is
$\prec$-increasing continuous, $M= \bigcup\limits_{\alpha} M_\alpha
\prec \gC_T$ is saturated of cardinality $\kappa^+$.

Assuming $T$ is dependent we should choose a thin enough club $E$ of
$\kappa^+$ such that $\{(M_\beta,M_\alpha):\alpha < \beta$ and
$\{\alpha,\beta\} \subseteq E \cap S^{\kappa^+}_\kappa\}$.  Now the club $E$
will be chosen such that all relevant pairs $(M_\beta,M_\alpha)$ are
similar enough to those of pairs from $\bold N^2_\kappa$.

So a sufficient condition for the conjecture is:
\mn
\begin{enumerate}
\item[$\boxplus$]  assume $(N_1,M_1),(N_2,M_2) \in \bold N^2_\kappa$,
  then we can find $\bar f = \langle f_s:s \in \cY \rangle$ such that:
\sn
\begin{enumerate}
\item[$(a)$]  $\cY$ is a partial order,
\sn
\item[$(b)$]  $\cY$ is $(< \kappa)$-complete, that is, any increasing
  chain (for $<_{\cY}$) of length $< \kappa$ has an upper bound,
\sn
\item[$(c)$]  $f_s$ is an $(N_1,N_2)$-elementary mapping,
\sn
\item[$(d)$]  $\Dom(f_s)$ has cardinality $< \kappa$,
\sn
\item[$(e)$]  $f_s$ maps $\Dom(f_s) \cap M_1$ onto $\Rang(f_s) \cap
 M_2$,
\sn
\item[$(f)_1$]  if $s \in \cY,A \in [N_1]^{< \kappa}$ \then \, for
  some $t \in \cY$ we have $s <_Y t \wedge A \subseteq \Dom(f_t)$,
\sn
\item[$(f)_2$]  if $s \in \cY$ and $A \in [N_2]^{< \kappa}$ \then \,
  for some $t \in \cY$ we have $s <_{\cY} t \wedge A \subseteq \Rang(f_t)$.
\end{enumerate}
\end{enumerate}
\mn
Now the approximation consists of $B_\ell \in [M_\ell]^{< \kappa}$ and
$\langle (\bar c^\ell_n,\bar d^\ell_n):n < \omega\rangle$, with $\bar
c^\ell_0,\bar d^\ell_0$ from $N_\ell$, which form a strict $(\kappa,<
\kappa)$-decomposition over $(M_\ell,B_\ell)$ for $\ell=1,2$ and $f$
an elementary mapping $h$ mapping $B_1$ onto $B_2$ and $(\bar c^1_n,\bar
d^1_n)$ to $(\bar c^2_n,\bar d^2_n)$.  So using $\bar c^\ell_n,\bar
d^\ell_n$ for $n > 0$ is to give us a condition with which we can
continue in a good induction hypothesis.

Now it should be clear what we like to have from strict decomposition;
however, the decomposition from \S2 are not enough.
\end{explanation}

\noindent
We first connect decomposition (i.e. the results of \ref{pr.28}) and
sufficient conditions for being an indiscernible sequence.
\begin{claim}
\label{pu.14}  
The sequence $\langle \bar a_\alpha:\alpha < \alpha^*\rangle$ is an
indiscernible sequence over $B$ \when \, for some $p,B$ we have:
\mn
\begin{enumerate}
\item[$(a)$]   $\bar a_\alpha$ is a sequence of members of $\gC$ and 
$\bar a_\alpha = \bar c_\alpha \char 94 \bar d_\alpha$ where $\ell
  g(\bar c_\alpha) =\ell g(\bar c_0),\ell g(\bar d_\alpha) = \ell
  g(\bar d_0)$ (not necessarily finite),
\sn
\item[$(b)$]  {\rm tp}$(\bar a_\alpha \char 94 \bar a_\beta,B) = p$
for $\alpha < \beta < \alpha^*$,
\sn
\item[$(c)$]   if $\alpha < \beta < \alpha^*$ then {\rm tp}$(\bar
d_\beta,\bar c_\beta + \bar d_\alpha + B) \vdash \text{\rm tp}(\bar
d_\beta,\bar c_\beta \cup \{\bar a_\gamma:\gamma \le \alpha\} \cup B)$,
\sn
\item[$(d)$]  {\rm tp}$(\bar c_\alpha,\cup\{\bar a_\beta:\beta <
\alpha\} \cup B)$ is increasing with $\alpha$, 
\sn
\item[$(e)$]   {\rm tp}$(\bar c_\alpha,\cup\{\bar a_\beta:\beta <
\alpha\} \cup B)$ does not split over $B$.
\end{enumerate}
\end{claim}

\begin{PROOF}{\ref{pu.14}}
For $u,v \subseteq \alpha^*$ let $B_{u,v} = 
\bigcup\{\bar a_\alpha:\alpha \in u\} \cup \bigcup\{\bar
c_\alpha:\alpha \in v\} \cup B$.  For $u,v \subseteq \alpha^*$ and
increasing functions $h_1$ from $u$ to $\alpha^*$ and $h_2$ to $v$ to
$\alpha_*$ such that $h_1 \rest (u \cap v) = h_2 \rest (u \cap v)$ (so
$h = h_1 \cup h_2$ is an increasing 
function $h$ from $u \cup v$ to $\alpha^*$) let 
$f = f_{h_1,h_2}$ be defined as follows:
\mn
\begin{enumerate}
\item[$\circledast$]   $(a) \quad \Dom(f) = B_{u,v}$,
\sn
\item[${{}}$]  $(b) \quad f \restriction B = \id_B$,
\sn
\item[${{}}$]   $(c) \quad f$ maps $\bar a_\alpha$ to $\bar
a_{h(\alpha)}$ for $\alpha \in u$,
\sn
\item[${{}}$]  $(d) \quad f$ maps $\bar c_\alpha$ to $\bar
c_{h(\alpha)}$ for $\alpha \in v$.
\end{enumerate}
\mn
Is $f_{h \rest u,h \rest v}$ a well defined function and even one to
one?  For this it suffices to check the following 
three demands, which follows by Clause (b) of the assumption
\mn
\begin{enumerate}
\item[$(*)_1$]   $(\alpha) \quad$ if $\alpha,\beta < \alpha^*,b \in B$
and $i < \ell g(\bar a_\alpha)$ \then \, $(\bar a_\alpha)_i = b
\Leftrightarrow (\bar a_\beta)_i = b$,
\sn
\item[${{}}$]  $(\beta) \quad$ if $\alpha,\beta < \alpha^*$ and
$i,j < \ell g(\bar a_\alpha)$ \then \, $(\bar a_\alpha)_i = (\bar
a_\alpha)_j \Leftrightarrow (\bar a_\beta)_i = (\bar a_\beta)_j$,
\sn
\item[${{}}$]   $(\gamma) \quad$ if $\alpha_1 < \alpha_2 < \alpha^*$ and
$\beta_1 < \beta_2 < \alpha^*$ and $i,j < \ell g(\bar a_\alpha)$ 
\then

\hskip25pt $(\bar a_{\alpha_1})_i = (\bar a_{\alpha_2})_j 
\Leftrightarrow (\bar a_{\beta_1})_i = (\bar a_{\beta_2})_j$.
\end{enumerate}
\mn
We prove by induction on $n$ that $\langle \bar a_\alpha:\alpha <
\alpha^*\rangle$ is an $n$-indiscernible sequence over $B$
(when $n < \alpha^*$).  For $n \le
2$ this is trivial by Clause (b) of the assumption.  So assume $n=m+1
> 2$ and we have proved it up to $m$.  So let $\alpha_0 < \ldots < \alpha_m
< \alpha^*,\beta_0 < \ldots < \beta_m < \alpha^*$ and we shall prove
that $\bar a_{\alpha_0} \char 94 \ldots \char 94 \bar a_{\alpha_m},\bar
a_{\beta_0} \char 94 \ldots \char 94 \bar a_{\beta_m}$ realize the
same type over $B$, this suffices.

Now by symmetry
\wilog \, $\alpha_m \le \beta_m$ let $h_0 =
\{(\alpha_\ell,\beta_\ell):\ell < m\},h_1 =
\{(\alpha_\ell,\alpha_\ell):\ell < m\},h_2 = h_1
\cup\{(\alpha_m,\beta_m)\}$ and $h_3 = h_0 \cup \{(\beta_m,\beta_m)\}$ 
and $h_4 = h_0 \cup \{(\alpha_m,\beta_m)\}$.

Let $f_0$ be the mapping $f_{h_0,h_0}$.  By the induction hypothesis
\mn
\begin{enumerate}
\item[$(*)_2$]   $f_0$ is an elementary mapping.
\end{enumerate}
\mn
Let $f_1$ be the mapping $f_{h_1,h_2}$, now by Clause (d) of the assumption
\mn
\begin{enumerate}
\item[$(*)_3$]  $f_1$ is an elementary mapping. 
\end{enumerate}
\mn
By $(*)_2$ we know that $\bar a_{\alpha_0} \char 94 \ldots \char 94
\bar a_{\alpha_{m-1}}$ and $f_0(\bar a_{\alpha_0} \char 94 \ldots
\char 94 \bar a_{\alpha_{m-1}})$ realize the same type over $B$ 
and they are included in $B_{\alpha_m,\alpha_m},B_{\beta_m,\beta_m}$
respectively.  But $\alpha_m \le \beta_m$ so both sequences are from
$B_{\beta_m}$ hence by Clause (e) of the assumption, i.e. as tp$(\bar
c_{\beta_m},B_{\beta_m,\beta_m})$ does not split over $B$, recalling
$(*)_2$ we have
\mn
\begin{enumerate}
\item[$(*)_4$]  $f_2 := f_{h_0,h_3}$ is an elementary mapping.
\end{enumerate}
\mn
By $(*)_3 + (*)_4$, comparing we have
\mn
\begin{enumerate}
\item[$(*)_5$]   $f_3 := f_{h_0,h_4}$ being $f_2 \circ f_1$ 
is an elementary mapping.
\end{enumerate}
\mn
 Note that
\mn
\begin{enumerate}
\item[$(*)_6$]   $f_3(\bar c_{\alpha_m} \char 94 \bar
d_{\alpha_{m-1}}) = \bar c_{\beta_m} \char 94 \bar d_{\beta_{m-1}}$.
\end{enumerate}
\mn
By Clause (b) of the assumption and $(*)_5 + (*)_6$
clearly 
\mn
\begin{enumerate}
\item[$(*)_7$]   $f_3(\text{tp}(\bar d_{\alpha_m},
\bar c_{\alpha_m} + \bar d_{\alpha_{m-1}} + B)) = \text{ tp}(\bar
d_{\beta_m},\bar c_{\beta_m} + \bar d_{\beta_{m-1}} + B)$.
\end{enumerate}
\mn
By Clause (c) of the assumption
\mn
\begin{enumerate}
\item[$(*)_8$]   $\tp(\bar d_{\alpha_m},
\bar c_{\alpha_m} + \bar d_{\alpha_{m-1}} + B) \vdash 
\tp(\bar d_{\alpha_m},\bar c_{\alpha_m} + \bar a_{\alpha_0} +
\ldots + \bar a_{\alpha_{m-1}} + B)$
\sn
\item[$(*)_9$]  $\tp(\bar d_{\beta_m},
\bar c_{\beta_m} + \bar d_{\beta_{m-1}} + B) \vdash 
\tp(\bar d_{\beta_m},\bar c_{\beta_m} + \bar a_{\beta_0} +
\ldots + \bar a_{\beta_{m-1}} + B)$.
\end{enumerate}
\mn
Together $f_4 = f_{h_4,h_4}$ is an elementary mapping and it maps $\bar
a_{\alpha_\ell}$ to $\bar a_{\beta_\ell}$ for $\ell \le m$ (and extend
id$_B$) so we are done. 
\end{PROOF}

\begin{observation}
\label{pu.16}   
The sequence $\langle (\bar c_t,\bar d_t):t \in I\rangle$ is a strict$_0 \,
(0,\theta)$-decomposition over $(B^0,B)$ \when \, for some $\langle
I_\ell,B,B^\ell:\ell < n\rangle$ we have:
\mn
\begin{enumerate}
\item[$(a)$]  $n < \omega$ and $n \ge 2$
\sn
\item[$(b)$]  the linear order $I$ is $I_0 + \ldots + I_n$ where
$I_\ell$ is infinite for $\ell=1,\dotsc,n-1$,
\sn
\item[$(c)$]  $\langle(\bar c_t,\bar d_t):t \in I_\ell\rangle$ is a 
strict$_0 \,(0,\theta)$-decomposition over $(B^\ell,B)$ for $\ell \le n$,
\sn
\item[$(d)$]  $\langle \bar c_t \char 94 \bar d_t:t \in I_\ell +
I_{\ell +1}\rangle$ is indiscernible over $B^\ell$ and $\ell g(\bar
c_t)$ for $t \in I_\ell \cup I_{\ell +1}$ is constant for $\ell < n$,
\sn
\item[$(e)$]  $B^{\ell +1} \supseteq \{\bar c_t \char 94
\bar d_t:t \in I_\ell\} \cup B^\ell$,
\sn
\item[$(f)$]   $B \subseteq B^0$.
\end{enumerate}
\end{observation}

\begin{PROOF}{\ref{pu.16}}
In Definition \ref{pr.35}, Clause (a) holds trivially as $B^0
\supseteq B$ by Clause (f) here (recalling that $\kappa$ there stands
for 0 here).  For Clause (b) of \ref{pr.35} the sequence $\langle \ell
g(\bar c_t):t \in I\rangle$ is constant as for each $\ell < n$ as the
sequence $\langle \ell g(\bar c_t):t \in I_\ell \cup I_{\ell
  +1}\rangle$ is constant (by Clause (c) here and (b) in \ref{pr.35})
and use transitivity of equality and $I_\ell$ for $\ell=1,\dotsc,n-1$
being non-empty by Clause (b) here.  Similarly $\langle \ell g(\bar
d_t):t \in I\rangle$ is constant, so \ref{pr.35}(b) indeed holds.

Similarly, Clause \ref{pr.35}(c) follows from \ref{pr.35}(d) proved
below and also Clause (c) here for $\ell=1$ (recalling $I_1$ is infinite).  
Clause \ref{pr.35}(e)$^-$
follows similarly using $B' \supseteq B^0$ by 
Clause (e) here and $(e)^-$ there 

So we are left with Clause \ref{pr.35}(d), that is $\langle \bar c_t
\char 94 \bar d_t:t \in I\rangle$ is an indiscernible sequence over
$B^0$.  For this we prove by induction on $k \le n$ that $\langle \bar
c_t \char 94 \bar d_t:t \in \cup\{I_\ell:\ell \in [n-k,n]\}\rangle$ is an
indiscernible sequence over $B^{n-k}$.  For $k=0,1$ this holds by
clause (d), for $k \ge 2$, let $s_0 <_I \ldots <_I s_{m-1},t_0 <_I
\ldots <_I t_{m-1}$ be from $\cup\{I_\ell:\ell \in [n-k,n]\}$ be
given.  Choose $s'_i,t'_i$ from $I_{n-k} \cup I_{n-k+1}$ for $i<m$
such that $s'_0 <_I \ldots <_I s'_{m-1},t'_0 <_I \ldots <_I t'_{m-1}$
and $s_\ell \in I_{n-k} \Rightarrow s'_\ell = s_\ell,s_\ell \notin
I_{n-k} \Rightarrow s'_\ell \in I_{n-k+1}$ and $t_\ell \in I_{n-k}
\Rightarrow t'_\ell = t_\ell,t_\ell \notin I_{n-k} \Rightarrow t'_\ell
\in I_{n-k+1}$.  This is possible as $I_{n-k+1}$ is infinite because 
$k \in [2,n]$.

Let $j \le m$ be such that $\ell < m \Rightarrow (s_\ell \in I_{n-k}
\Leftrightarrow \ell < j)$; so by the induction hypothesis we have
$\tp(\bar a_{s_j} \char 94 \ldots \char 94 \bar a_{s_{m-1}},B^{n-k+1})
= \tp(\bar a_{s'_j} \char 94 \ldots \char 94 \bar
a_{s'_{m-1}},B^{n-k+1})$.  As $\ell < j \Rightarrow \bar a_{s_\ell} =
\bar a_{s'_\ell} \subseteq B^{n-k+1}$ and $B^{n-k} \subseteq
B^{n-k+1}$ it follows that $\tp(\bar a_{s_0} \char 94 \ldots \char 94 \bar
a_{s_{m-1}},B^{n-k})$ is equal to tp$(\bar a_{s'_0} 
\char 94 \ldots \char 94 \bar a_{s'_{m-1}},B^{n-k})$.  This type by
clause (d) is equal to tp$(\bar a_{t'_0} \char 94 \ldots \char 94 \bar
a_{t'_{m-1}},B^{n-k})$.  Similarly to the proof in the beginning of
the paragraph, this type by the induction hypothesis is equal to
$\tp(\bar a_{t_0} \char 94 \ldots \char 94 \bar a_{t_{m-1}},B^{n-k})$,
so together we are done.
\end{PROOF}

\noindent
The following is a local version of \ref{pu.14} (see \ref{du.22})
\begin{claim}
\label{du.19}  
Assume ($n(*) < \omega$ and $\alpha(*) > n(*))$ for each $k < n(*)$)
\mn
\begin{enumerate}
\item[$(a)$]  $\bar c^k_\alpha \in {}^{\gamma(k,1)} {\frak C}$ for
$k \le n(*),\alpha < \alpha(*)$; this means that the length of $\bar
  c^k_\alpha$ may depend on $m$ but not on $\alpha$ and may be infinite,
\sn
\item[$(b)$]  $\bar d_\alpha = \bar d^0_\alpha 
\in {}^{\gamma(k,0)} {\frak C}$ for $\alpha < \alpha(*)$,
\sn
\item[$(c)$]  $\bar c^k_\alpha \triangleleft \bar c^{k+1}_\alpha$
for $\alpha < \alpha(*),k < n(*)$, 
\sn
\item[$(d)$]  $\bar e^k_\alpha \in {}^{\gamma(m,2)} {\frak C}$ 
for non-zero $k \le n(*),\alpha < \alpha(*)$,
\sn
\item[$(e)$]  for all $\alpha < \beta < \alpha(*)$ 
the type {\rm tp}$(\bar c^{k+1}_\beta 
\char 94 \bar e^{k+1}_\alpha \char 94 \bar d_\beta \char 94 
\bar e^k_\beta,B)$ is the same,
\sn 
\item[$(f)$]  $\tp(\bar d_\beta \char 94 \bar e^k_\beta,\bar
c^k_\beta + \bar e^{k+1}_\alpha + B) \vdash \text{\rm tp}(\bar
d^k_\beta \char 94 \bar e^k_\beta,\bar c^k_\beta + 
\sum\limits_{\gamma < \beta} \bar c^k_\gamma \char 94 \bar d_\gamma + B)$, 
\sn
\item[$(g)$]  {\rm tp}$(\bar c^k_\alpha,\cup\{\bar c^{k+1}_\beta \char
  94 \bar d_\beta \char 94 \bar e^{k+1}_\beta:\beta < \alpha\} \cup B)$
is increasing with $\alpha$,
\sn
\item[$(h)$]  {\rm tp}$(\bar c^k_\alpha,\cup\{\bar c^{k+1}_\beta 
\char 94 \bar d_\beta \char 94 \bar e^{k+1}_\beta:\beta < \alpha\} 
\cup B)$ does not split over $B$,
\sn
\item[$(i)$]  {\rm tp}$(\bar c^{n(*)}_\alpha \char 94 
\bar d_\alpha \char 94 \bar e^k_\alpha,B)$ is the 
same for all $\alpha < \alpha(*)$.
\end{enumerate}
\mn
\Then \, $\langle \bar c^0_\alpha \char 94
\bar d_\alpha:\alpha < \alpha^*\rangle$ is an 
$n(*)$-indiscernible sequence over $B$.
\end{claim}

\begin{remark}
\label{du.22}
1) In what sense is \ref{du.19} a local version of \ref{pu.14}?  In
the second we get only $n(*)$-indiscernibility.  Note that the role of $\bar
d_\alpha$ there is played by $\bar d_\alpha,\bar e^m_\alpha (m < n(*))$ here.

\noindent
2) The claim is not used in the rest of the section.
\end{remark}

\begin{PROOF}{\ref{du.19}}
We prove by induction on $n < n(*)$ that
\mn
\begin{enumerate}
\item[$\odot_1$]   if $k \le n(*)-n$ and $\alpha_0 < \ldots <
\alpha_n$ and $\beta_0 < \beta_1 < \ldots < \beta_n$ then

$\bar c^k_{\alpha_0} \char 94 \bar d_{\alpha_0} \char 94 \ldots \char 94
\bar c^k_{\alpha_{n-1}} \char 94 \bar d_{\alpha_{n-1}} \char 94 
\bar c^k_{\alpha_n} \char 94 \bar d_{\alpha_n} \char 94 
\bar e^k_{\alpha_n}$ and 

$\bar c^k_{\beta_0} \char 94 \bar d_{\beta_0} \char 94
\ldots \char 94 \bar c^k_{\alpha_{n-1}} \char 94 \bar d_{\alpha_{n-1}}
\char 94 \bar c^k_{\beta_n} \char 94 \bar d_{\beta_n} \char 94 \bar
e^k_{\beta_n}$ realize the same type over $B$.
\end{enumerate}
\mn
For $n = n(*)-1,k=0$ we get the desired conclusion.
\smallskip

\noindent
For $n=0$ this holds by clause (i) of the assumption.  So assume
$n=m+1$ and $k \le n(*)-n$ and we have proved this for $m$.  Note that
$k+1 \le n(*)-m$.  So let $\alpha_0 < \ldots < \alpha_{m+1} <
\alpha^*,\beta_0 < \ldots < \beta_{m+1} < \alpha^*$ be given and \wilog \,
$\alpha_n \le \beta_n$ and we shall proof the equality of types from
$\odot_1$ in this case, this suffice.  Now
\mn
\begin{enumerate}
\item[$(*)_1$]  $\bar c^{k+1}_{\alpha_0} \char 94 \bar d_{\alpha_0}
\char 94 \ldots \char 94 \bar c^{k+1}_{\alpha_{m-1}} \char 94 
\bar d_{\alpha_{m-1}} \char 94 \bar c^{k+1}_{\alpha_m} \char 94
\bar d_{\alpha_m} \char 94 \bar e^{k+1}_{\alpha_m}$ and

$\bar c^{k+1}_{\beta_0}
\char 94 \bar d_{\beta_0} \char 94 \ldots \char 94 \bar
c^{k+1}_{\beta_{m-1}} \char 94 \bar d^{k+1}_{\alpha_{m-1}} \char 94
\bar c^{k+1}_{\alpha_m} \char 94 \bar d_{\alpha_m} \char 94
e^{k+1}_{\alpha_m}$ realize

the same type over $B$.
\end{enumerate}
\mn
[Why?  By the induction hypothesis.]
\mn
\begin{enumerate}
\item[$(*)_2$]  tp$(\bar c^k_{\beta_n},\cup\{\bar c^{k+1}_\gamma
\char 94 \bar d_\gamma \char 94 \bar e^{k+1}_\gamma:\gamma < \beta_n\}
\cup B)$ extends tp$(\bar c^k_{\alpha_n},\cup\{\bar c^{k+1}_\gamma 
\char 94 \bar d_\gamma \char 94 \bar e^{k+1}_\gamma:\gamma <
\alpha_n\} \cup B)$.
\end{enumerate}
\mn
[Why?  By clause (g) of the assumption.]
\mn
\begin{enumerate} 
\item[$(*)_3$]  tp$(\bar c^k_{\beta_n},\cup\{\bar
c^{k+1}_\gamma \char 94 \bar d_\gamma \char 94 \bar
e^{k+1}_\gamma:\gamma < \beta_n\} \cup B)$ does not split over $B$.
\end{enumerate}
\mn
[Why?  By clause (h) of the assumption.]
\mn
\begin{enumerate}
\item[$(*)_4$]  $\bar c^{k+1}_{\alpha_0} \char 94 \bar d_{\alpha_0}
\char 94 \ldots \char 94 \bar c^{k+1}_{\alpha_{m-1}} \char 94 
\bar d_{\alpha_{m-1}} \char 94 \bar c^{k+1}_{\alpha_m} \char 94 \bar
d_{\alpha_m} \char 94 \bar e^{k+1}_{\alpha_m} \char 94 \bar
c^k_{\beta_{m+1}}$ and

$\bar d^{k+1}_{\beta_0} \char 94 
\bar d_{\beta_0} \char 94 \ldots \char 94 \bar c^{k+1}_{\beta_{m-1}}
\char 94 \bar d_{\beta_{m-1}} \char 94 \bar c^{k+1}_{\beta_m} \char
94 \bar d_{\beta_m} \char 94 \bar e^{k+1}_{\beta_m} \char 94 \bar
c^k_{\beta_{m+1}}$ realize 

the same type over $B$.
\end{enumerate}
\mn
[Why?  By $(*)_1 + (*)_3$.]
\mn
\begin{enumerate}
\item[$(*)_5$]  in $(*)_4$ we can replace $\bar c^k_{\beta_{m+1}}$
in the first sequence by $\bar c^k_{\alpha_{m+1}}$.
\end{enumerate}
\mn
[Why?  By $(*)_4 + (*)_2$.]   
\mn
But
\mn
\begin{enumerate}
\item[$(*)_6$]  $(\bar c^{k+1}_{\alpha_{m+1}} \char 94 
\bar e^{k+1}_{\alpha_m}) \char 94 
(\bar d_{\alpha_{m+1}} \char 94 \bar e^k_{\alpha_{m+1}})$ and
$(\bar c^{k+1}_{\beta_{m+1}} \char 94 \bar e^{k+1}_{\beta_m}) \char 94
(\bar d_{\beta_{m+1}} \char 94 \bar e^k_{\beta_{m+1}})$ realize the
same type over $B$.
\end{enumerate}
\mn
[Why?  By clause (e) of the assumption.]
\mn
\begin{enumerate}
\item[$(*)_7$]  in $(*)_4,(*)_5$ we can replace $\bar
  c^{k+1}_{\alpha_\ell} (\ell \le m)$ by $\bar c^k_{\alpha_\ell}$.
\end{enumerate}
\mn
[Why?  As $\bar c^k_\gamma \triangleleft \bar c^{k+1}_\gamma$ by clause (c)
of the assumption.]
\mn
\begin{enumerate}
\item[$(*)_8$]  $\bar c^k_{\alpha_0} \char 94 \bar d_{\alpha_0}
\char 94 \ldots \char 94 \bar c^k_{\alpha_{m-1}} \char 94 \bar
d_{\alpha_{m-1}} \char 94 \bar c^k_{\alpha_m} \char 94 \bar
d_{\alpha_m} \char 94 \bar c^k_{\alpha_{m+1}} \char 94
\bar d^m_{\alpha_{m+1}} \char 94 \bar e^k_{\beta_{m+1}}$ and

$\bar c^k_{\beta_0} \char 94 \bar d_{\beta_0} \char 94 \ldots \char 94
\bar c^k_{\beta_{m-1}} \char 94 \bar d_{\beta_{m-1}} \char 94 \bar
c^k_{\beta_m} \char 94 \bar d_{\beta_m} 
\char 94 \bar c^k_{\beta_{m+1}} \char 94 \bar d_{\beta_{m+1}} \char 94 \bar
e^k_{\beta_{m+1}}$ 

realize the same type over $B$.
\end{enumerate}
\mn
[Why?  By $(*)_7 + (*)_6$ and clause (f) of the assumption.]

We finish the induction step.  Hence we get the desired statement. 
\end{PROOF}
\bigskip

\centerline {$* \qquad * \qquad *$}
\bigskip

\noindent
We now return to the generic pair conjecture.  Central here is the
following definition; the best case is $\lambda = \kappa = \theta$ is
a measurable cardinal.
\begin{definition}
\label{ps.7} 
We say that the triple $(\lambda,\kappa,< \theta)$ is 
good or $T$-good \when \,:
\mn
\begin{enumerate}
\item[$(A)$]  $|T| < \theta = \text{ cf}(\theta) \le \kappa =
\text{ cf}(\kappa) \le \lambda = \lambda^{< \kappa}$,
\sn
\item[$(B)$]  $T$ is dependent,
\sn
\item[$(C)$]  if $M$ is $\kappa$-saturated of cardinality $\le
\lambda$ and $\bar d \in {}^{\theta >}{\gC}$ \then \, we
can find $B \subseteq M$ of cardinality $< \kappa$ and a strict 
$(\kappa,< \theta)$-decomposition
$\langle(\bar c_n,\bar d_n):n < \omega\rangle$ over $(M,B)$ 
such that $\bar d \trianglelefteq \bar d_0$,
\sn
\item[$(D)$]  if $M$ is $\kappa$-saturated of cardinality $\le
  \lambda,B \subseteq M$ has cardinality $< \kappa,\langle (\bar
  c_n,\bar d_n):n < \omega\rangle$ is a strict $(\kappa,<
  \theta)$-decomposition over $(M,B)$ and $\bar d \in {}^{\theta
  >}\gC$ \then \, there is a strict $(\kappa,< \theta)$-decomposition
  $\langle (\bar c^+_n,\bar d^+_n):n < \omega\rangle$ over $(M,B)$
  such that $\bar c_n \trianglelefteq \bar c^+_n,\bar d_n
  \triangleleft \bar d^+_n$ for $n < \omega$ and
$\bar d_0 \char 94 \bar d \trianglelefteq \bar d^+_0$.
\end{enumerate}
\end{definition}

\noindent
So to begin our analysis we need
\begin{observation}
\label{ps.11}
1) If $T$ is a dependent and $\lambda >|T|$ is a measurable cardinal
\then \, $(\lambda,\lambda,< \lambda)$ is $T$-good.

\noindent
2) If $T$ is dependent, $\kappa$ a supercompact cardinality and
   $\lambda = \lambda^{< \lambda} \ge \kappa$, then $(\lambda,\kappa,<
   \kappa)$ is $T$-good.
\end{observation}

\begin{PROOF}{\ref{ps.11}}
Immediate by \ref{pr.49}(2) and \ref{pr.23}(2), you may use
\ref{pr.28}, too.
\end{PROOF}

\noindent
For the rest of this section we assume, till but not including the end
that is, \ref{ps.56}.
\begin{hypothesis}
\label{ps.13}  
1) $T$ is dependent.

\noindent
2) $|T| < \theta = \cf(\theta) \le \kappa 
\le \lambda$ and $\lambda = \lambda^{< \lambda} > |T|$ 
and $(\lambda,\kappa,< \theta)$ is $T$-good.
\end{hypothesis}

\begin{claim}
\label{ps.21}  
Assume
\mn
\begin{enumerate}
\item[$(a)$]   $\delta$ is a limit ordinal $< \theta = \text{\rm cf}(\theta)$,
\sn
\item[$(b)$]  $\langle(\bar c^\alpha_n,\bar d^\alpha_n):n <
\omega\rangle$ is a strict $(\kappa,< \theta)$-decomposition over
$(M,B_\alpha)$ for each $\alpha < \delta$,
\sn
\item[$(c)$]  $\bar c^\alpha_n \trianglelefteq \bar c^\beta_n \wedge
\bar d^\alpha_n \trianglelefteq \bar d^\beta_n$ for $\alpha < \beta <
\delta,n < \omega$,
\sn
\item[$(d)$]  $B_\alpha \subseteq B_\beta$ for $\alpha < \beta < \delta$,
\sn
\item[$(e)$]  we define $\bar c^\delta_n = \cup\{\bar
c^\alpha_n:\alpha < \delta\},\bar d^\delta_n = \cup\{\bar
d^\alpha_n:\alpha < \delta\},B_\delta = \cup\{B_\alpha:\alpha <
\delta\}$.
\end{enumerate}
\mn
\Then \, $\langle(\bar c^\delta_n,\bar d^\delta_n):n < \omega\rangle$
is a strict $(\kappa,< \theta)$-decomposition over $(M,B_\delta)$.
\end{claim}

\begin{PROOF}{\ref{ps.21}}
We have to check Clauses (a)-(f) of Definition \ref{pr.35}(1). 
Clause (a) is trivial by assumption (b) of \ref{ps.21} recaling
$\delta < \theta = \cf(\theta)$.  
Clause (b) holds as $\delta < \theta = \cf(\theta)$ by
assumption (a) of \ref{ps.21} and $\langle \ell g(\bar c^\delta_n):n <
\omega\rangle$ is constant by assumptions (b),(c) and similarly
$\langle \ell g(\bar d^\delta_n):n < \omega\rangle$ is constant.  Next 
Clauses (c),(d),(e) hold by their local character and assumptions 
(b) + (c) of \ref{ps.21}.  

Lastly, proving Clause (f) is the main point, it means to show:
\mn
\begin{enumerate}
\item[$\odot_1$]   if $B_\delta \subseteq A \subseteq M$ and $|A| <
\kappa$ then for some pair $(\bar c,\bar d)$ of sequences from
$M$ we have $\langle \bar c \char 94 \bar d \rangle \char 94 \langle \bar
c^\delta_n \char 94 \bar d^\delta_n:n < \omega\rangle$ is an
indiscernible sequence over $A$.
\end{enumerate}
\mn
Toward this by induction on $\alpha < \delta$ we choose a pair $(\bar
c^*_\alpha,\bar d^*_\alpha)$ such that
\mn
\begin{enumerate}
\item[$\circledast_1$]   $\langle \bar c^*_\alpha \char 94 \bar
d^*_\alpha\rangle \char 94 \langle c^\alpha_n \char 94 \bar
d^\alpha_n:n < \omega\rangle$ is an indiscernible sequence over $A \cup
\bigcup\{\bar c^*_\beta \char 94 \bar d^*_\beta:\beta < \alpha\}$.
\end{enumerate}
\mn
[Why possible?  
We can choose $(\bar c^*_\alpha,\bar d^*_\alpha)$ because $\langle (\bar
c^\alpha_n,\bar d^\alpha_n):n < \omega\rangle$ being a strict
$(\kappa,< \theta)$-decomposition over $B_\alpha$, we can apply clause
(f) of Definition \ref{pr.35} recalling $B_\alpha$ being a
subset of $A \cup \{(\bar c^* \char 94 \bar d^*:\beta < \alpha\}$ and the
later being $\subseteq M$ and of cardinality $< \kappa$ as $\kappa$ is
regular $\ge \theta > \delta$.]

We can find $\langle(\bar c_\alpha,\bar d_\alpha):\alpha <
\delta\rangle$ though not necessarily in $M$ such that
\mn
\begin{enumerate}
\item[$\circledast_2$]  $(a) \quad \langle \bar c_\alpha \char 94
\bar d_\alpha:\alpha < \delta\rangle \char 94 \langle \bar c^\delta_n
\char 94 \bar d^\delta_n:n < \omega\rangle$ is an indiscernible
sequence over $A$
\sn
\item[${{}}$]  $(b) \quad \bar c^*_\alpha \trianglelefteq \bar
c_\alpha$ and $\bar d^*_\alpha \trianglelefteq \bar d_\alpha$ for
$\alpha < \delta$ and, of course, 

\hskip25pt  $\ell g(\bar c^*_\alpha) = \ell g(\bar c^\delta_\alpha),
\ell g(\bar d_\alpha) = \ell g(\bar d^\delta_\alpha)$.
\end{enumerate}
\mn
[Why?  For this by using the saturation of ${\frak C}$, it is enough
to prove that: if $n < m < \omega,\alpha_0 < \ldots < \alpha_{n-1} <
\alpha_n < \delta$ then the sequences 
$\bar c^*_{\alpha_0} \char 94 \bar d^*_{\alpha_0}
\char 94 \ldots \char 94 \bar c^*_{\alpha_{n-1}} \char 94 
\bar d^*_{\alpha_{n-1}}$ 
$\bar c^{\alpha_0}_0 \char 94 \bar d^{\alpha_0}_0 \char 94 \bar
c^{\alpha_1}_1 \char 94 \bar d^{\alpha_1}_1 \char 94 \ldots \char 94
\bar c^{\alpha_{n-1}}_{n-1} \char 94 \bar d^{\alpha_{n-1}}_{n-1}$
realize the same type over
$A \cup \{\bar c^{\alpha_n}_k \char 94 \bar d^{\alpha_n}_k:
k \in [n,m)\}$.  This is proved by
induction on $n < \omega$ and is straight.]

Hence we can find a pair $(\bar c',\bar d')$ such that:
\mn
\begin{enumerate}
\item[$\circledast_3$]  $\langle 
\bar c' \char 94 \bar d'\rangle \char 94 \langle \bar c_\alpha
\char 94 \bar d_\alpha:\alpha < \delta\rangle \char 94 \langle \bar
c^\delta_n \char 94 \bar d^\delta_n:n < \omega\rangle$ is an indiscernible
sequence over $A$.
\end{enumerate}
\mn
Lastly, we choose $(\bar c'',\bar d'')$ such that
\mn
\begin{enumerate}
\item[$\circledast_4$]  $(\bar c'',\bar d'')$ is a pair of sequences
from $M$ such that $\bar c'' \char 94 \bar d''$ realizes 
$\tp(\bar c' \char 94 \bar d',A \cup 
\bigcup\{\bar c^*_\alpha \char 94 \bar d^*_\alpha:\alpha < \delta\})$,
of course with $\ell g(\bar c'') = \ell g(\bar c^\delta_0),\ell g(\bar
d'') = \ell g(\bar d^\delta_0)$; equivalently there is an 
automorphism of ${\frak C}$ which is the
identity on $A \cup \bigcup\{\bar c^*_\alpha \char 94 \bar
d^*_\alpha:\alpha < \delta\}$ mapping $\bar c' \char 94 \bar d'$ to
$\bar c'' \char 94 \bar d''$.
\end{enumerate}
\mn
We shall prove that $(\bar c'',\bar d'')$ is as required.
Now to prove that $(\bar c'',\bar d'')$ is as required in
$\odot_1$ it suffices to prove, for each $\alpha < \delta$ that
\mn
\begin{enumerate}
\item[$\odot_2$]  $\langle (\bar c'' \restriction \ell g(\bar
c^\alpha_0)) \char 94 (\bar d'' \restriction \ell g(\bar
d^\alpha_0)) \rangle \char 94 \langle \bar c^\alpha_n \char 94 
\bar d^\alpha_n:n < \omega\rangle$

\hskip15pt  is an indiscernible sequence over $A$.
\end{enumerate}
\mn
[Why?  As $\odot_1$ is a ``local" demand, i.e. it says that $\bar c''
\char 94 \bar d''$ is a sequence realizing an appropriate type $q$
(and is from $M$) and for this it suffices to check every finite
subtype so $\odot_2$ suffices.]

Now $\odot_2$ follows by $\circledast_6$ below.
Let $(\bar c^*_{\beta,\gamma},\bar d^*_{\beta,\gamma}) =
(\bar c^*_\beta \rest \ell g(\bar c^\gamma_0),
\bar d^*_\beta \rest \ell g(d^\gamma_0))$ and $(\bar c''_\gamma,\bar
d''_\gamma) = (\bar c'' \rest \ell g(\bar c^\gamma_0),
\bar d'' \rest \ell g(\bar d^\gamma_0))$ for $\gamma < \delta$ 
and $\beta < \delta$; we use only $\beta \in [\gamma,\delta)$.

Now
\mn
\begin{enumerate}
\item[$\circledast_5$]   $\langle \bar c''_\alpha \char 94 \bar
  d''_\alpha \rangle \char 94 \langle \bar c^*_{\beta,\alpha} \char 94
\bar d^*_{\beta,\alpha}:\beta \in [\alpha,\delta)\rangle$ is 
a strict$_0$ decomposition over $(A,A)$.
\end{enumerate}
\mn
[Why?  By $\circledast_3$, this holds for 
$\langle(\bar c' \restriction \ell g(\bar c^\alpha_0)) \char 94
(\bar d' \restriction \ell g(\bar c^\alpha_0))\rangle 
\char 94 \langle \bar c^*_{\beta,\alpha} \char 94 
\bar d^*_{\beta,\alpha}:\beta \in [\alpha,\delta)\rangle$ and 
we use preservation by automorphism of
${\frak C}$, i.e. use $\circledast_4$.]
\mn
\begin{enumerate}
\item[$\circledast_6$]   For $\alpha < \delta$ and $i \le \omega$ the
sequence $\langle \bar c''_\alpha \char 94 \bar d''_\alpha \rangle
\char 94 \langle \bar c^*_{\beta,\alpha} \char 94 
\bar d^*_{\beta,\alpha}:\beta \in [\alpha,\delta)\rangle 
\char 94 \langle c^\alpha_n \char 94 \bar d^\alpha_n:n
< i\rangle$ is a strict$_0$ decomposition over $(A,A)$.
\end{enumerate}
\mn
[Why?  We prove this by induction on $i$, noticing that for $i = \omega$ we get
the desired conclusion; also for $i = \omega$ the inductive step is trivial
and for $i=0$ use $\circledast_5$.  So assume $i=n+1$, let
\mn
\begin{enumerate}
\item[$(*)_1$]   $\bullet \quad A_1 = A \cup \text{ Rang}(\bar c''_\alpha) 
\cup \text{ Rang}(\bar d''_\alpha)$,
\sn
\item[${{}}$]  $\bullet \quad A_2 = A_{2,\delta}$ where
\sn
\item[${{}}$]  $\bullet \quad$ for 
$\gamma \in [\alpha,\delta]$ we let $A_{2,\gamma} = 
\cup\{\text{Rang}(\bar c^*_{\beta,\alpha}) \cup 
\text{ Rang}(\bar d^*_{\beta,\alpha}):\beta \in [\alpha,\gamma)\}$,
\sn
\item[${{}}$]  $\bullet \quad A_3 = \cup\{\bar c^\alpha_\ell
\char 94 \bar d^\alpha_\ell:\ell < n\}$.
\end{enumerate}
\mn
Clearly
\mn
\begin{enumerate}
\item[$(*)_2$]  $(a) \quad A_1,A_{2,\delta}$ are $\subseteq M$ and
  $A_3 \subseteq \gC$,
\sn
\item[${{}}$]  $(b) \quad A_{2,\gamma}$ is $\subseteq$-increasing for
  $\gamma \in [\alpha,\delta)$.
\end{enumerate}
\mn
Let $\gamma \in (\alpha,\delta)$ be a successor ordinal and we shall
prove
\mn
\begin{enumerate}
\item[$(*)_3$]  $(c^*_{\gamma,\alpha} \char 94 \bar d^*_{\gamma,\alpha})
\char 94 \ldots \char 94 (\bar c^*_{\gamma +n,\alpha} \char 94 \bar
d^*_{\gamma +n,\alpha})$ and $(\bar c^\alpha_0 \char 94 \bar
d^\alpha_0) \char 94 \ldots \char 94 (\bar c^\alpha_n \char 94 \bar
d^\alpha_n)$ realize the same type over $A_1 \cup A_{2,\gamma}$.
\end{enumerate}
\mn
As $\delta$ is a limit ordinal this suffices by \ref{pu.16} 
(with $n = 2,I_0$ a singleton, $I_1$ isomorphic to $\delta$, $I_2$
isomorphic to $\omega$).  Now $(\bar c^*_{\gamma,\alpha} 
\char 94 \bar d^*_{\gamma,\alpha}) \char 94 \ldots
\char 94 (\bar c^*_{\gamma +n-1,\alpha} \char 94 \bar d^*_{\gamma
+n-1,\alpha})$ and $(\bar c^\alpha_0 \char 94 \bar d^\alpha_0) \char 94
\ldots \char 94 (\bar c^\alpha_{n-1} \char 94 \bar d^\alpha_{n-1})$
realize the same type over $A_1 \cup A_{2,\gamma}$ by the induction
hypothesis.  

Next by Clause (b) of the assumption 
$\tp(\bar c^\alpha_n,M \cup \bigcup\{\bar c^\alpha_m 
\char 94 \bar d^\alpha_m:m < n\})$ does not split over $B_\alpha$ 
hence $\tp(\bar c^\alpha_n,A \cup A_1 \cup A_2 \cup A_3)$ 
does not split over $B_\alpha \subseteq A$ by $(*)_1$.

Hence by the induction hypothesis
$\langle \bar c''_\alpha \char 94 \bar d''_\alpha \rangle 
\char 94 \langle \bar c^*_{\beta,\alpha} \char 94 d^*_{\beta,\alpha}:\beta \in
[\alpha,\delta)\rangle \char 94 \langle (\bar c^\alpha_k \char 94 
\bar d^\alpha_k):k < n\rangle$ is an indiscernible sequence over $A \cup
\bar c^\alpha_n$, hence
\mn
\begin{enumerate}
\item[$\bullet$]  $(\bar c^*_{\gamma,\alpha} \char 94 \bar
  d^*_{\gamma,\alpha}) \char 94 \ldots \char 94 (\bar c^*_{\gamma +n-1,\alpha}
\char 94 \bar d^*_{\gamma +n-1,\alpha}) \char 94 \bar c^\alpha_n$ and
  $(\bar c^\alpha_0 \char 94 \bar d^\alpha_0) \char 94 \ldots \char 94
  (\bar c^*_{\gamma +n-1,\alpha} \char 94 \bar d^*_{\gamma
  +n-1,\alpha}) \char 94 \bar c^\alpha_n$ realize the same type over
  $A_1 \cup A_{2,\gamma}$.
\end{enumerate}
\mn
But by $\circledast_1$
clearly $\bar c^\alpha_n,\bar c^*_{\gamma +n,\alpha}$ realize the
same type over $A_1 \cup A_{2,\gamma}$ hence
\mn
\begin{enumerate}
\item[$\bullet$]  $(c^*_{\gamma,\alpha} \char 94 
\bar d^*_{\gamma,\alpha}) \char 94 \ldots \char 94 (\bar c^*_{\gamma
  +n-1,\alpha} \char 94 \bar d^*_{\gamma +n-1,\alpha}) \char 94 \bar c_{\gamma
+n,\alpha}$ and $(\bar c^\alpha_0 \char 94 \bar d^\alpha_0) \char 94
\ldots \char 94 (\bar c^\alpha_{n-1} \char 94 \bar d^\alpha_{n-1})
\char 94 \bar c^\alpha_n$ realize the same type over $A_1 \cup
A_{2,\gamma}$.
\end{enumerate}
\mn
We can choose $\bar d'$ in $\gC$ such that
$(\bar c^*_{\gamma,\alpha} \char 94 \bar d^*_{\gamma,\alpha}) \char 94
\ldots \char 94 
(\bar c^*_{\gamma+n-1,\alpha} \char 94 \bar d^*_{\gamma +n-1,\alpha})
\char 94 (\bar c^*_{\gamma +n,\alpha} \char 94 \bar d')$ and
$(\bar c^\alpha_0 \char 94 \bar d^\alpha_0) \char 94 \ldots \char 94
(\bar c^\alpha_{n-1} \char 94 \bar d^\alpha_{n-1}) \char 94
(\bar c^\alpha_n \char 94 \bar d^\alpha_n)$ realize the same type
over $A_1 \cup A_{2,\gamma}$ so (to prove $\circledast_6$) 
it suffices to prove that 
$\bar d^*_{\gamma +n,\alpha},\bar d'$ realize the same type over $A_1
\cup A_{2,\gamma +n} \cup \bar c^*_{\gamma +n,\alpha}$.  

Recall that $\langle \bar c^*_{\beta,\alpha} \char 94 \bar
d^*_{\beta,\alpha}:\beta \in [\alpha,\delta)\rangle \char 94 \langle
\bar c^\alpha_n \char 94 \bar d^\alpha_n:n < \omega\rangle$ is 
an indiscernible sequence, hence
the sequences $\bar c^*_{\gamma +n-1,\alpha} \char 94 \bar
d^*_{\gamma +n-1,\alpha} \char 94 \bar c^*_{\gamma+n,\alpha} 
\char 94 \bar d^*_{\gamma +n,\alpha}$ and
$\bar c^*_{\gamma +n -1,\alpha} \char 94 \bar d^*_{\gamma +n-1,\alpha} 
\char 94 \bar c^\alpha_n \char 94 \bar d^\alpha_n$ realize the same type.

By the two previous sentences and 
the transitivity of the equality of types $\bar d^*_{\gamma
+n,\alpha},\bar d'$ realize the same type over 
$(\bar c^*_{\gamma +n-1,\alpha} \char 94 \bar d^*_{\gamma +n-1,\alpha})
\char 94 \bar c^*_{\gamma +n,\alpha}$, but by Clause (e)$^-$ of
Definition \ref{pr.35} which apply by Clause (b) of the assumption and
$\circledast_5$ above we have $\tp(\bar d^*_{\gamma +n,\alpha},
\bar c^*_{\gamma +n,\alpha} + \bar d^*_{\gamma
+n-1,\alpha}) \vdash \text{ tp}(\bar d^*_{\gamma +n,\alpha},A_1 +
A_{2,\gamma +n} + \bar c^*_{\gamma +n,\alpha})$ so we are done.]
\end{PROOF}

\begin{definition}
\label{ps.25}
1) We say that the strict $(\kappa,< \theta)$-decompositions $\langle
(\bar c'_\varepsilon,\bar d'_\varepsilon):\varepsilon < 
\delta\rangle,\langle (\bar c''_\varepsilon,\bar
d''_\varepsilon):\varepsilon < \delta \rangle$ over $M$ are equivalent over $B
\in [M]^{< \lambda}$ \when \, for some automorphism $f$ of $M$ over $B$ for
every $n$ and $\varepsilon_0 < \ldots < \varepsilon_{n-1} < \delta,f$
maps the type $\tp((\bar c'_{\varepsilon_\delta} \char 94 \bar
d'_{\varepsilon_\delta}) \char 94 \ldots \char 94 (\bar
c'_{\varepsilon_{n-1}} \char 94 \bar d'_{\varepsilon_{n-1}}),M)$ to
the type $\tp((\bar c''_{\varepsilon_0} \char 94 \bar
d''_{\varepsilon_0}) \char 94 \ldots \char 94 (\bar
c''_{\varepsilon_{n-1}} \char 94 \bar d''_{\varepsilon_{n-1}}),M)$ and
$\ell g(\bar c'_{\varepsilon_0}) = \ell g(\bar
c''_{\varepsilon_0}),\ell g(\bar d'_{\varepsilon_0}) = \ell g(\bar
d''_{\varepsilon_0})$. 

\noindent
2) In part (1) we say ``weakly equivalent over $B$" \when \, for every
   $\zeta < \kappa$ and $\bar b' \in {}^\zeta M$ there is $\bar b''
   \in {}^\zeta M$ and vice versa and elementary mapping $f$ such
   that: $f \supseteq \id_B,f(\bar b') = \bar b''$ and $f(\bar
   c'_\varepsilon) = f(\bar c''_\varepsilon),f(\bar d'_\varepsilon) =
   f(\bar d''_\varepsilon)$ for $\varepsilon < \delta$.

\noindent
3) If $B = \emptyset$ then we may omit it.
\end{definition}

\begin{claim}
\label{ps.23}
1) If $M$ is $\kappa$-saturated of cardinality $\lambda$ 
and $B \in [M]^{< \lambda}$ \then \, the number of
strict $(\lambda,\kappa,< \theta)$-decompositions $\langle (\bar c_n,\bar
d_n):n < \omega \rangle$ over $(M,B)$ such that $\ell g(\bar d_n) > 0$
up to weak equivalence or when $\ell g(\bar d_n) = 0$ up to equivalence
over $B$ is $\le \lambda$, see \ref{ps.25}.

\noindent
2) For $M,B$ as above, two strict $(\kappa,<\theta)$-decompositions
   are equivalent \when \, they are weakly equivalent above $(M,B)$
   and $\lambda = \kappa$.

\noindent
3) if $\langle \bar c_n,\bar d_n):n < \omega\rangle$ is a 
$(\kappa,< \theta)$-decomposition over $(M,B)$ so $B \in [M]^\kappa$ and $C =
\cup\{\rang(\bar c_n \char 94 \bar d_n):n < \omega\} \cup B$ \then
 \, $M_{[C]}$ is a sequence homogeneous model. 
\end{claim}

\begin{PROOF}{\ref{ps.23}}
1) First, if $\ell g(\bar d_n)=0$, 
by \ref{3k.0.7}(4), that is \cite[5.26]{Sh:783}, however in the
present case $\lambda$ is measurable hence strongly inaccessible so
\ref{ps.23}(1) is easy.  That is, fixing $A,B$, also the number of
$\alpha < \lambda$ and $p(\bar x_\alpha) \in \bold S(M)$ not splitting
over $B$ is $\le \lambda$ by \ref{3k.0.7}.  The case for weakly
equivalent holds as $\lambda$ is strongly inaccessible.

\noindent
2) Let $\langle (\bar c'_n,\bar d'_n):n < \omega),\langle (\bar
c''_n,\bar d''_n):n < n\rangle$ be two strict $(\kappa,<
\theta)$-decompositions over $M$.
As a $\lambda$-sequence-homogeneous model of cardinality $\lambda$
is determined up to isomorphisms by the set of complete types of
 finite tuples in it (by \cite{KM67} or see \cite[\S2]{Sh:88r}) by
   part (3) it suffices to show:
\mn
\begin{enumerate}
\item[$(*)_1$]  for every $\bar b' \in {}^{\omega >} M$ there are
  $\bar b'' \in {}^{\omega >}M$ and an elementary mapping $f$ such
  that $f(\bar b') = \bar b'',f \supseteq \id_B$ and $\varepsilon <
  \delta \Rightarrow f(\bar c'_\varepsilon) = \bar c''_\varepsilon$,
\sn
\item[$(*)_2$]  for every $\bar b'' \in {}^{\omega >} M$ there are
$\bar b' \in {}^{\omega >}M$ there $f$ as above.
\end{enumerate}
\mn
By symmetry it suffices to prove $(*)_1$.  Given $\bar b'$ let $B^* =
B \cup \bar b$.  Let $\delta < \kappa$ be a finite ordinal.  We choose
$(\bar c^*_\varepsilon,\bar d^*_\varepsilon)$, a pair of sequences from
$M$ such that $\langle \bar c'_\varepsilon \char 94 \bar
d'_\varepsilon \rangle \char 94 \langle \bar c'_n \char 94 \bar d'_n:n
< \omega\rangle$ is an indiscernible sequence over $B' \cup
\bigcup\{\bar c^*_\zeta \char 94 \bar d^*_\zeta:\zeta < \varepsilon\}$
 by \ref{pr.35}(f) with $\cup\{\bar c'_\zeta \char 94
\bar d'_\zeta:\zeta < \varepsilon\} \cup B',M,\langle (\bar c'_\varepsilon,\bar
d'_\varepsilon):\varepsilon < \delta\rangle,(\bar c'_\varepsilon,\bar
d'_\varepsilon)$ here standing for 
$A,M,\langle (\bar c_t,\bar d_t):t \in J\rangle,(\bar
c,\bar d)$ there and similarly $(\bar c''_\varepsilon,\bar
d''_\varepsilon)$ for $\varepsilon < \delta$ with $B'' := B$.

Now as $M$ is $\kappa$-saturated we can find $\bar b'' \in 
{}^{\ell g(\bar b')}M$ such that the sequences $\bar b' \char 94 \bar c' \char
  94 \bar d' \char 94 (\ldots \char 94 \bar c'_\varepsilon \char 94
  \ldots)_{\varepsilon < \delta}$ and $\bar b'' \char 94 \bar c''
\char 94 \bar d'' \char 94 (\ldots \char 94 \bar c''_\varepsilon
  \char 94 \ldots)_{\varepsilon < \delta}$ realize the same type over $A$.  
We can finish by \ref{pr.35}(1)(e) as in proof of \ref{ps.21}.

\noindent
3) Let $g$ be the identity mapping on $\cup\{\Rang(\bar c_n \char 94
   \bar d_n):n < \omega\} \cup B$.

Let
\mn
\begin{enumerate}
\item[$(*)_1$]  $\cF$ be the set of $f$ such that:
\sn
\begin{enumerate}
\item[$\bullet$]  $f$ is an elementary mapping
\sn
\item[$\bullet$]  $\Dom(f) \subseteq M$ has cardinality $< \kappa$
\sn
\item[$\bullet$]  $\Rang(f) \subseteq M$
\sn
\item[$\bullet$]  $f \cup g$ is an elementary mapping.
\end{enumerate}
\end{enumerate}
\mn
It suffices to prove $(*)_1 - (*)_4$ below:
\mn
\begin{enumerate}
\item[$(*)_1$]  $\cF \ne \emptyset$.
\end{enumerate}
\mn
[Why?  As $f = \id_B$ belongs to $\cF$.]
\mn
\begin{enumerate}
\item[$(*)_2$]  $f \in \cF$ \Iff \, $f^{-1} \in \cF$
\sn
\item[$(*)_3$]  $\cF$ is closed under union of increasing chains of
  length $< \lambda$.
\end{enumerate}
\mn
[Why?  Just check.]
\mn
\begin{enumerate}
\item[$(*)_4$]  if $f \in \cF,\ell \in \{1,2\}$ and $a_\ell \in M$ then
  for some $a_{3-\ell} \in M$ we have $f \cup \{(a_1,a_2)\} \in \cF$.
\end{enumerate}
\mn
Why?  By symmetry, \wilog \, $\ell=1$. Let $\delta$ be a limit ordinal
$< \kappa$.  We choose $\bar c'_\varepsilon,\bar d'_\varepsilon$
sequences from $M$ of length $\ell g(\bar c_0),\ell g(\bar d_0)$
respectively by induction on $\varepsilon < \delta$, by applying
\ref{pr.35}(f) to $M,\langle (\bar c_n,\bar d_n):n < \omega\rangle$
 and $B_\varepsilon = B \cup \Dom(f) \cup \{a\} \cup \Rang(f) \cup
 \bigcup\{\bar c'_\zeta \char 94 \bar d'_\zeta:\zeta < \varepsilon\}$
 and get $\bar c_\varepsilon,\bar d_\varepsilon$ from $M$ as there.
\mn
\begin{enumerate}
\item[$(*)_{4.1}$]  $\langle \bar c'_\varepsilon \char 94 \bar
  d'_\varepsilon:\varepsilon < \delta \rangle \char 94 \langle
(\bar c_n \char 94 \bar d_n):n < \omega\rangle$ is an indiscernible sequence
  over $B_0$.
\end{enumerate}
\mn
This implies that, letting $B_* = \cup\{\bar c'_\varepsilon \char 94
\bar d'_\varepsilon:\varepsilon < \delta\}$:
\mn
\begin{enumerate}
\item[$(*)_{4.2}$]  $f \cup g \cup \id_{B_*}$ is an elementary mapping.
\end{enumerate}
\mn
As $M$ is $\kappa$-saturated
\mn
\begin{enumerate}
\item[$(*)_{4.3}$]  there is $a_2 \in M$ such that
$f \cup \id_{\bar c \char 94 \bar d} \cup \{(a_1,a_2)\}$ is an
  elementary mapping.
\end{enumerate}
\mn
So we can prove as in the proof of \ref{ps.21}
\mn
\begin{enumerate}
\item[$(*)_{4.4}$]  $f \cup \{(a_1,a_2)\} \in \cF$.
\end{enumerate}
\mn
So we are done.
\end{PROOF}

Clearly \ref{ps.21} is a step forward.  Now we prove
 the generic pair conjecture; instead of assuming that the
cardinality $\lambda$ is measurable we can restrict $T$.

\noindent
Toward this
\begin{definition}
\label{ps.28}  
We say that the triple $\bold m = (M,N,{\cA})$ is a 
$(\lambda,\kappa,< \theta)$-system \when \,
($\lambda \ge \kappa \ge \theta = \text{ cf}(\theta) > |T|$ and) it
satisfies clauses (a)-(d) below, and say the triple is a full
$(\lambda,\kappa,< \theta)$-system when it satisfies clauses
(a)-(i) below where:
\mn
\begin{enumerate}
\item[$\boxplus$]   $(a) \quad M$ is $\kappa$-saturated of 
cardinality $\lambda$,
\sn
\item[${{}}$]   $(b) \quad M \prec N \prec {\gC}$ and $N$ has
cardinality $\lambda$,
\sn
\item[${{}}$]   $(c) \quad {\cA}_{\bold m}$ is a set of 
cardinality $\lambda$ of objects $\bold p$ such that:
\sn
\begin{enumerate}
\item[${{}}$]  $(\alpha) \quad \bold p$ is of the form 
$\langle(\bar c_\varepsilon,\bar d_\varepsilon):
\varepsilon < \lambda^+\rangle = \langle(\bar c_\varepsilon[\bold p],
\bar d_\varepsilon[\bold p]):\varepsilon < \lambda^+\rangle$,
\sn
\item[${{}}$]   $(\beta) \quad \bar c_0 \char 94 \bar d_0 \subseteq N$,
\sn
\item[${{}}$]   $(\gamma) \quad \langle(\bar c_\varepsilon,
\bar d_\varepsilon):\varepsilon < \lambda^+ \rangle$ is a 
strict $(\kappa,< \theta)$-decomposition over $M$,
\end{enumerate}
\sn
\item[${{}}$]   $(d) \quad {\cA}_{\bold m}$ is 
partially ordered by: $\bold p \le_{\bold m} \bold q$ iff
\sn
\item[${{}}$]  $\qquad \bullet \quad \bar c_0[\bold p] \trianglelefteq
\bar c_0[\bold q]$ and $\bar d_0[\bold p] \trianglelefteq \bar d_0[\bold q]$,
\sn
\item[${{}}$]  $\qquad \bullet \quad$ for every large enough
$\varepsilon < \lambda^+,\bar c_\varepsilon[\bold p] \trianglelefteq
\bar c_\varepsilon[\bold q],\bar d_\varepsilon[\bold p]
\trianglelefteq \bar d_\varepsilon[\bold q]$,
\sn
\item[${{}}$]  $(e) \quad {\cA}_{\bold m}$ is closed under union of $< \theta$
increasing chains of length $< \theta$,
\sn
\item[${{}}$]   $(f) \quad$ if $\bold p \in {\cA}_{\bold m}$ and 
$\bar d \in {}^{\theta >} N$ then for some $\bold q \in \cA$ above
$\bold p$ we have 

\hskip25pt $\Rang(\bar d) \subseteq \Rang(\bar d_0[\bold q])$,
\sn
\item[${{}}$]  $(g) \quad N$ is saturated,
\sn
\item[${{}}$]  $(h) \quad$  up to really equivalence every possible $\bold p$
in Clause (c) occurs (i.e. is 

\hskip25pt represented in $\cA$), where $\bold p,\bold q$ are really
equivalent \when \, they are

\hskip25pt as in (c) and $(\bar c_{\bold p,0},\bar d_{\bold p,0}) =
(\bar c_{\bold q,0},\bar d_{\bold q,0})$
\sn
\item[${{}}$]  $(i) \quad$ there is $\bold q_2 \in \cA_{\bold m}$ such that
$\bold p_2 \le_{\bold m} \bold q_2$ and $\bold q_2,\bold q_1$ are

\hskip25pt equivalent as witnessed by $f$, \when \,:
\sn
\begin{enumerate}
\item[${{}}$]  $(\alpha) \quad \bold p_1 \le_{\bold m} \bold q_1$ so
both are from $\cA_{\bold m}$,
\sn
\item[${{}}$]  $(\beta) \quad \bold p_2 \in \cA_{\bold m}$ 
is really equivalent to $\bold p_1$,
\sn
\item[${{}}$]  $(\gamma) \quad f \in \aut(M)$ maps $\bold p_1$ to $\bold p_2$.
\end{enumerate}
\end{enumerate}
\end{definition}

\begin{definition}
\label{ps.35}  
1) Let $\BP = \BP_{\lambda,\kappa,< \theta}$ be the set of
   $(\lambda,\kappa,< \theta)$-systems.

\noindent
2) If $(M_\ell,N_\ell,{\cA}_\ell)$ is a 
$(\lambda,\kappa,< \theta)$-system for $\ell=1,2$ we say
$(M_2,N_2,{\cA}_2)$ is above $(M_1,N_1,{\cA}_1)$
or $(M_1,N_1,{\cA}_1) \le_{\BP_{\lambda,\kappa,< \theta}} (M_2,N_2,{\cA}_2)$
\when \, $M_1 = M_2,N_1 \prec N_2$ and ${\cA}_1 \subseteq {\cA}_2$.

\noindent
3) We may write just $\BP,\le_{\BP}$ \when \, 
$(\lambda,\kappa,< \theta)$ is clear from the context.

\noindent
4) If $\bold m \in \BP$ and $\bold p \in \cA_{\bold m}$, we say that
   $B$ is a base of $\bold p$ \when \, $B \subseteq M_{\bold m}$ has
   cardinality $< \kappa$ and $\bold p$ is a strict $(\kappa,<
   \theta)$-decomposition over $(M,B)$.
\end{definition}

\begin{claim}
\label{ps.42}  
Assume $\kappa = \lambda = \lambda^{< \lambda}$ (see Definition \ref{ps.7}). 

\noindent
1) If $M \prec {\gC}$ is $\kappa$-saturated of cardinality
$\lambda$ \then \, there is a pair $(N,{\cA})$ such that 
$(M,N,{\cA})$ is a $(\lambda,\kappa,< \theta)$-system.

\noindent
2) If $(M,N_1,{\cA}_1)$ is a $(\lambda,\kappa,< \theta)$-system 
and $N_2 \prec {\gC}$ and $\|N_2\| = \lambda$ \then \, there is a
pair $(N_3,{\cA}_3)$ such that $(M,N_3,{\cA}_3)$ is a
$(\lambda,\kappa,< \theta)$-system above $(M,N_1,{\cA}_1)$ and $N_2
\prec N_3$.

\noindent
3) If $\langle (M,N_\varepsilon,{\cA}_\varepsilon):\varepsilon <
\delta\rangle$ is an increasing sequence of $(\lambda,\kappa,< \theta)$-systems
and $\delta$ is a limit ordinal $< \lambda^+$ \then \, the union, 
$(M,\bigcup\limits_{\alpha < \delta} N_\alpha,\bigcup\limits_{\alpha < \delta} 
{\cA}_\alpha)$ is a $(\lambda,\kappa,< \theta)$-system which is a least 
upper bound of
$\{(M,N_\varepsilon,{\cA}_\varepsilon):\varepsilon < \delta\}$.

\noindent
4) If in part (3) we have {\rm cf}$(\delta) = \lambda$ and each is
full, \then \,  so is the union.
\end{claim}

\begin{PROOF}{\ref{ps.42}}
1) Let $\bold m = (M,M,\emptyset)$ and check.

\noindent
2) Let $N_3 \prec \gC$ be such that $N_1 \cup N_2 \prec N_3,N_3$ is
   saturated of cardinality $\lambda$ and use $(M,N_3,\cA_1)$.

\noindent
3) Easy.

\noindent
4) Easy.
\end{PROOF}

\begin{claim}
\label{ps.43}
1) We have:
\mn
\begin{enumerate}
\item[$(a)$]  $\le_{\BP}$ is a partial order on $\BP$,
\sn
\item[$(b)$]   any $\le_{\BP}$-increasing sequence of
length $< \lambda^+$ has an upper bound,
\sn
\item[$(b)^+$]  moreover if $\langle
(M,N_\alpha,\cA_\alpha):\alpha < \delta\rangle$ is
$\le_{\BP}$-increasing \then \, $(M,\cup\{N_\alpha:\alpha <
\delta\},\cup\{\cA_\alpha:\alpha < \delta\})$ is a least 
$\le_{\BP}$-upper bound,
\sn
\item[$(c)$]  $\BP$ is not empty, moreover for every $M \in
\EC_{\lambda,\lambda}(T)$ there is $\bold m \in \BP$ such that
$M_{\bold m} = M = N_{\bold m},\cA_{\bold m} = \emptyset$.
\end{enumerate}
\mn
2) If $\bold m$ is a $(\lambda,\kappa,<\theta)$-system \then \,:
\mn
\begin{enumerate}
\item[$(a)$]   $\le_{\bold m}$ is a partial order of $\cA_{\bold m}$,
\sn
\item[$(b)$]  if also $\bold n$ is a $(\lambda,\kappa,< \theta)$-system and
$\bold m \le_{\BP} \bold n$ then
$\le_{\bold m} = \le_{\bold n} \rest \cA_{\bold m}$,
\sn
\item[$(c)$]  if $\delta$ as a limit ordinal
$< \theta$ and $\bold p_\alpha \in \cA_{\bold m}$ for $\alpha <
\delta$ is $\le_{\bold m}$-increasing with $\alpha$ \then \, there is
$\bold n \in \BP$ such that $\bold m \le_{\BP} \bold n$ and
$\le_{\bold n}$-upper bound of $\{\bold p_\alpha:\alpha < \delta\}$,
\sn
\item[$(d)$]  if $\bold p \in \cA_{\bold m}$
and $\bar d \in {}^{\theta >} N$ \then \, there are a
$(\lambda,\kappa,< \theta)$-system $\bold n$
satisfying $\bold m \le_{\BP} \bold n$ and $\bold q \in \cA_{\bold n}$
such that $\bold p \le_{\bold n} \bold q$ and $\bar d \subseteq
\Rang(\bar d^{\bold q}_0)$,
\sn
\item[$(e)$]  if $\bold p_1,\bold p_2 \in \cA_{\bold m}$ are
equivalent and $\bold p_1 \le_{\bold m} \bold q_2$ \then \, for some
$\bold n \in \BP$ we have $\bold m \le_{\BP} \bold n$ and for some
$\bold q_2 \in \cA_{\bold n}$ equivalent to $\bold q_1$ we have $\bold
p_1 \le_n \bold q_2$.
\end{enumerate}
\mn
\noindent
3) If $(M,N,\cA)$ is a $(\lambda,\kappa,< \theta)$-system \then \, there
 is a full $(\lambda,\kappa,< \theta)$-system above it.
\end{claim}

\begin{PROOF}{\ref{ps.43}}
1) Clause (a) holds easily by checking \ref{ps.35}(2).  Clauses
(b),(b)$^+$ holds easily by \ref{ps.42}(3), and for Clause (c) use
 $(M,M,\emptyset)$.

\noindent
2) Why?  \underline{Clauses (a),(b)}: Easy.

\noindent
\underline{Clause (c)}: For $\alpha < \beta < \delta$ let 
$\zeta = \zeta(\alpha,\beta) < \lambda^+$ be such that:
\mn
\begin{enumerate}
\item[$\bullet$]   if $\varepsilon \in [\zeta,\lambda^+)$ then $\bar
c^{\bold p_\alpha}_\varepsilon \trianglelefteq \bar c^{\bold
p_\beta}_\varepsilon,\bar d^{\bold p_\alpha}_\varepsilon
\trianglelefteq \bar d^{\bold p_\beta}_\varepsilon$.
\end{enumerate}
\mn
Let $\zeta(*) = \sup\{\zeta(\alpha,\beta) +2:\alpha < \beta <
\delta\}$, so necessarily $\zeta(*) < \lambda^+$.  Clearly if $\zeta
\in [\zeta(*),\lambda^+) \cup \{0\}$ 
then $(\bar c^{\bold p_\alpha}_\zeta:\alpha
< \delta\rangle$ is $\trianglelefteq$-increasing and also the sequence
 $\langle \bar d^{\bold p_\alpha}_\zeta:
\alpha < \delta\rangle$ is $\trianglelefteq$-increasing.

We choose $\bar c_0 = \cup\{\bar c^{\bold p_\alpha}_0:\alpha <
\delta\},\bar d_0 = \cup\{\bar d^{\bold p_\alpha}_0:\alpha < \delta\}$
and for $\varepsilon < \lambda^+$ let $\bar c_{1 + \varepsilon} =
\cup\{\bar c^{\bold p_\alpha}_{\zeta(*)+\varepsilon}:\alpha <
\delta\}$ and $\bar d_{1 + \varepsilon} = \cup\{\bar d^{\bold
p_\alpha}_{\zeta(*) + \varepsilon}:\alpha < \delta\}$.

Now easily
\mn
\begin{enumerate}
\item[$\bullet$]  $\bold p = \langle (\bar c_\varepsilon,\bar
d_\varepsilon):\varepsilon < \lambda^+\rangle$ is well defined and
satisfies the demand in \ref{ps.28}$(c)(\alpha)+(\beta)+(\gamma)$.
\end{enumerate}
\mn
[Why?  By Claim \ref{ps.21}.]

So we define $\bold n$ as $(M,N_{\bold m},\{\bold p\} \cup \cA_{\bold
m})$, easily
\mn
\begin{enumerate}
\item[$\bullet$]  $\bold n \in \BP$ and $\bold m \le \bold n$.
\end{enumerate}
\mn
Lastly, $\alpha < \delta \Rightarrow \bold p_\alpha \le_{\bold n}
\bold p$ as witnessed by $\zeta(*) \times \omega$.  So we are done
proving clause (c).

\noindent
\underline{Clause (d)}:

Easy.  By Clause (D) of Definition \ref{ps.7} of
``$(\lambda,\kappa,< \theta)$ is $T$-good" which holds by Hypothesis 
\ref{ps.13}.

\noindent
\underline{Clause (e)}:

By the definition of ``equivalence" in Definition \ref{ps.25} (and
\ref{ps.28}(h)).

\noindent
3) We fix a $(\lambda,\kappa,< \theta)$-system $(M,N_*,\cA_*)$.
We now shall choose $\bold m_i \in \BP,\le_{\BP}$-increasing by induction on
$i \le \lambda$ and
\mn
\begin{enumerate}
\item[$\odot_1$]  for $i=0:\bold m_0$ is $(M,N_*,\cA_*)$
\sn
\item[$\odot_2$]  for limit $i < \lambda,\bold m_i$ is a
$\le_{\BP}$-upper bound of $\langle \bold m_j:j < i\rangle$, in fact
the union.
\end{enumerate}
\mn
[Why?  Possible by clause $(b)^+$ of part (1).]

By bookkeeping
\mn
\begin{enumerate}
\item[$\odot_3$]  if $\langle \bold p_\alpha:\alpha < \delta\rangle$
is $\le_{\bold m_j}$-increasing, $\delta$ a limit ordinal $< \theta$
and $j < \lambda$, 
then for some $i \in (j,\lambda)$ the sequence has an upper bound by
$\le_{\bold m_{i+1}}$.
\end{enumerate}
\mn
[Why?  By Clause (c) of part (2).]
\mn
\begin{enumerate} 
\item[$\odot_4$]   if $\bold p \in \cA_{\bold m_j}$ and $\bar d \in
{}^{\theta >}N_{\bold m_j}$ and $j < \lambda$ 
\then \, for some $i \in (j,\lambda)$,
there is $\bold q \in \cA_{\bold m_{i+1}}$ such that $\bold p
\le_{\bold m_{i+1}} \bold q$ and $\Rang(\bar d) \subseteq \Rang(\bar
d^{\bold q}_0)$. 
\end{enumerate}
\mn
[Why?  By Clause (d) of part (2).]
\mn
\begin{enumerate} 
\item[$\odot_5$]   if $j < \lambda$ and $\bold p_1,\bold q_1,\bold
  p_2,f$ satisfy Clauses $(\alpha)-(\varepsilon)$ of \ref{ps.28}(i)
  with $\bold m_j$ here standing for $\bold m$ there \then \, for some $i
  \in (j,\lambda)$ there is $\bold q_2 \in \cA_{\bold m_{i+1}}$ such
  that $\bold p_2 \le_{\bold m_{i+1}} \bold q_2$ and $\bold q_1,\bold
  q_2$ are equivalent as witnessed by $f$.
\end{enumerate}
\mn
[Why?  Just think.]

So we can carry the induction.  Now $\bold m_\lambda$ is as required.
\end{PROOF}

\begin{theorem}
\label{ps.49}  
Assume $(\lambda,\kappa,< \theta) = (\lambda,\lambda,< \lambda)$ so is a
$T$-good triple, see \ref{ps.7}, \ref{ps.13}.

\noindent
1) If $\lambda^+ = 2^\lambda,\langle M_\alpha:\alpha
<\lambda^+\rangle$ is an $\prec$-increasing continuous sequence of
members of $\EC_\lambda(T),M_\alpha$ saturated if $\alpha$ is non-limit
and $M = \cup\{M_\alpha:\alpha < \lambda^+\}$ is saturated \then \,
for some club (= closed unbounded subset) 
$E$ of $\lambda^+$ for any $\alpha < \beta < \delta \in
E$ and $\alpha,\beta$ are non-limit or are from
$S^{\lambda^+}_\lambda$, the pairs $(M_\delta,M_\alpha)$ and
$(M_\delta,M_\beta)$ are isomorphic.

\noindent
2) If $\bold m_\ell = (M_\ell,N_\ell,\cA_\ell)$ is a full
$(\lambda,\lambda,< \lambda)$-system for $\ell=1,2$ \then \,
$(N_1,M_1) \cong (N_2,M_2)$ that is there is an isomorphism $f$ from
$N_1$ onto $N_2$ mapping $M_1$ onto $M_2$.
\end{theorem}

\begin{PROOF}{\ref{ps.49}}
1) By part (2), noting that $M_\alpha,M_\beta$ are saturated and
   recallihng \ref{ps.43} and its proof.

\noindent
2) We define the set AP of approximation:
\mn
\begin{enumerate}
\item[$(*)_1$]   AP is the set of triples $\bold h = (\bold
p_1,B_1,\bold p_2,B_2,f) = (p_1[\bold h],B_2[\bold h],p_2[\bold
h],B_2[\bold h],f[\bold h])$ satisfying:
\sn
\begin{enumerate}
\item[$(a)$]  $\bold p_\ell \in \cA_\ell$ for $\ell=1,2$,
\sn
\item[$(b)$]  $B_\ell \in [M_\ell]^{< \lambda}$ is a base for $\bold
p_\ell$, see \ref{ps.35}(4),
\sn
\item[$(c)$]  $f$ is an elementary mapping which maps $B_1$
onto $B_2$ such that $\Dom(f) = B_1$,
\sn
\item[$(d)$]  there is an isomorphism $f^+$ from $M_1$ onto $M_2$
extending $f$ such that: if $\alpha_0 < \ldots < \alpha_{n-1} < \lambda^+$
\then \, $f^+$ maps $\tp((\bar c_{\bold p_1,\alpha_1} \char 94 \bar
d_{\bold p_1,\alpha_1}) \char 94 \ldots \char 94 (\bar c_{\bold
p_1,\alpha_{n-1}} \char 94 \bar d_{\bold p_1,\alpha_{n-1}}),M_1)$ onto
$\tp((\bar c_{\bold p_2,\alpha_0} \char 94 \bar d_{\bold
p_2,\alpha_0}) \char 94 \ldots \char 94 (\bar c_{\bold
p_2,\alpha_{n-1}} \char 94 \bar d_{\bold p_2,\alpha_{n-1}}),M_2)$.
\end{enumerate}
\item[$(*)_2$]  we define the two-place relation $\le_{\AP}$ by: 
$\bold h_1 \le_{\AP} \bold h_2$ \Iff \,
\sn
\begin{enumerate}
\item[$(a)$]  both are in AP,
\sn
\item[$(b)$]  $\bold p_\ell[\bold h_1] \le_{\bold m_\ell} \bold
p_\ell[\bold h_2]$ for $\ell=1,2$,
\sn
\item[$(c)$]  $B_\ell[\bold h_1] \subseteq B_\ell[\bold h_2]$,
\sn
\item[$(d)$]  $f_1[\bold h_1] \subseteq f[\bold h_2]$.
\end{enumerate}
\end{enumerate}
\mn
Obviously
\mn
\begin{enumerate}
\item[$(*)_3$]  $\le_{\AP}$ partially ordered AP.
\end{enumerate}
\mn
Also
\mn
\begin{enumerate}
\item[$(*)_4$]  if $\delta < \lambda$ is a limit ordinal and $\langle
\bold h_\alpha:\alpha < \delta\rangle$ is $\le_{\AP}$-increasing,
\then \, this sequence has an $\le_{\AP}$-upper bound.
\end{enumerate}
\mn
[Why?  As in the proof of \ref{ps.43}(2)(c).]
\mn
\begin{enumerate}
\item[$(*)_5$]  if $\bold h \in \AP$ and $a \in M_1$
then there is $\bold h' \in \AP$ such that $\bold h \le_{\AP} \bold
h'_1$ and $a \in B_{\bold p_\ell[\bold h'_1]}$. 
\end{enumerate}
\mn
[Why?  Let $a_1 = a$ and let $f^+ \supseteq f_{\bold h}$ be as in
  $(*)_1(d)$ and let $a_2 \in M_2$ realize 
$f^+(\tp(a_1,B_{\bold p_1[\bold h]}))$.

Let 
$B'_\ell = (B_{\bold p_\ell[\bold h]} \cup \{a_\ell\})$ and let
$\bold p'_\ell \in \cA_{\bold m_\ell}$ be defined by:
$(\bar c_\varepsilon[\bold p'_\ell],\bar d_\varepsilon[\bold p'_\ell])
= (\bar c_\varepsilon[\bold p_\ell[\bold h_1],\bar d_\varepsilon[\bold
p_\ell[\bold h])$.

Lastly, let $\bold h' = (\bold p'_1,B'_1,\bold p'_2,B'_2,f_{\bold h} \cup
\{(a_1,a_2)\})$.]
\mn
\begin{enumerate}
\item[$(*)_6$]  If $\bold h \in \AP$ and $a \in M_2$
\then \, there is $\bold h' \in \AP$ such that $\bold h \le_{\AP}
\bold h'$ and $a \in B_2[\bold h']$.
\end{enumerate}
\mn
[Why?  Like $(*)_5$.]
\mn
\begin{enumerate}
\item[$(*)_7$]   if $\bold h \in \AP$ and $d \in N_1$ then for some
$\bold h' \in \AP$ we have $\bold h \le_{\AP} \bold h'$ and $d \in
\Rang(\bar d_0[\bold p_1[\bold h']])$.
\end{enumerate}
\mn
[Why?  There is $\bold q_1 \in \cA_{\bold m_1}$ such that $\bold
p_1[\bold h] \le_{\bold m_1} \bold q_1$ and $d \in \bar d_0[\bold
p_1[\bold h']$.]

Let $f' \supseteq f_{\bold h}$ be as in $(*)_1(d)$ so in particular 
an isomorphism from $M_1$ onto $M_2$.

Now by clause (i) of Definition \ref{ps.28} there is $\bold q_2 \in
\cA_{\bold m_2}$ such that $f'$ maps $\bold q_1$ to $\bold q_2$.

The rest should be clear.]
\mn
\begin{enumerate}
\item[$(*)_8$]  if $\bold h \in \AP$ and $d \in N_2$ \then \, for some
$\bold h' \in \AP$ we have $\bold h \le_{\AP} \bold h'$ and $d \in
\Rang(\bar d_0[\bold p_1(\bold h')])$.
\end{enumerate}
\mn
[Why?  Like $(*)_7$.]

Together
\mn
\begin{enumerate}
\item[$\odot$]  there  is a sequence $\langle \bold h_i:i <
\lambda\rangle$ such that (for $\ell=1,2$)
\sn
\begin{enumerate}
\item[$(a)$]  it is $\le_{\AP}$-increasing
\sn
\item[$(b)_\ell$]  if $a \in M_\ell$ then $a \in B_\ell[\bold h_i]$
for some $i$
\sn
\item[$(c)_\ell$]  if $d \in N_\ell$ then $d \in \Rang(\bar d_0[\bold
p_\ell(\bold h_i)])$ for some $i < \lambda$.
\end{enumerate}
\end{enumerate}
\mn
\relax From this sequence we can ``read" an isomorphism as required, say
$g(a_1) = a_2$ \Iff \, for some $i$ and $\varepsilon < \ell g(\bar
d_0[\bold p_1[\bold h_i])$ we have $a_1 = 
(\bar d_0[\bold p_1[\bold h_i]])_\varepsilon,d_2 
= (\bar d_0[\bold p_2[\bold h_i]])_\varepsilon$.]
\end{PROOF}

\noindent
Another form, not assuming Hypothesis \ref{ps.13}, is
\begin{conclusion}
\label{ps.56}  
Assume $(\lambda,\lambda,< \lambda)$ is $T$-good, e.g. $\lambda > |T|$
 is a measurable cardinality and $\lambda = \lambda^{< \lambda}$.  
Then for some club $\bold F$ we have:
\mn
\begin{enumerate}
\item[$(A)$]   $(a) \quad \bold F$ is as in \cite[3.3]{Sh:88r}, i.e.
\sn
\begin{enumerate}
\item[$(\alpha)$]  $\bold F$ is a function with domain
$\{\bar M:\bar M$ has the form $\langle M_i:i \le \beta\rangle$, a
$\prec$-increasing continuous sequence such that $M_i$ is models of
$T$ of cardinality $\lambda$ with universe an ordinal $\in
[\lambda,\lambda^+)$ and if $i$ is non-limit then $M_i$ is saturated$\}$,
\sn
\item[$(\beta)$]  ${\bold F}(\bar M)$ is such that $\bar M
\char 94 \langle \bold F(\bar M)\rangle \in \Dom(\bold F)$
\end{enumerate}
\sn
\item[$(B)$]   $\bar M = \langle M_\alpha:\alpha <
\lambda^+\rangle$ obeys $\bold F$ which means that 

\hskip25pt $\lambda^+ = \sup\{\alpha:\bold F
(\bar M \restriction (\alpha +1) \prec M_{\alpha +1}\}$ \then \,

\hskip25pt for some club $E$ of $\lambda^+$ we have:
\sn
\item[${{}}$]   $(a) \quad \cf(\alpha) = \lambda \Rightarrow M_\alpha$ is
saturated,
\sn
\item[${{}}$]  $(b) \quad$ if 
$M_{\alpha_\ell}$ is saturated, $\cf(\delta_\ell) =
\lambda$ and $\alpha_\ell < \delta_\ell \in E$ for $\ell=1,2$ 

\hskip25pt then $(M_{\delta_1},M_{\alpha_1}) \cong (M_{\delta_2},M_{\beta_2})$.
\end{enumerate}
\end{conclusion}
\newpage

\bibliographystyle{alphacolon}
\bibliography{lista,listb,listx,listf,liste,listy}

\end{document}